\documentclass[a4paper,12pt]{article}
\usepackage[utf8x]{inputenc}
\usepackage{arabtex}
\usepackage{afterpage}
\usepackage{xcolor}
\usepackage{relsize}
\usepackage{mathrsfs}
\usepackage[all]{xy}
\usepackage{amsmath,amssymb,amsfonts,amscd, graphicx, epstopdf, latexsym, verbatim, multirow, color}
\usepackage{epsfig}
\usepackage{geometry} 
\usepackage{adjustbox}

\newcommand{\F}{\mathcal F}
\newcommand{\T}{\mathcal T}
\newcommand{\Q}{\mathbf Q}
\newcommand{\Z}{\mathbf Z}

\def\xor{\underline{\lor}}
 \renewcommand{\R}{\mathbf R}

\newcommand{\psl}{\mathrm{PSL}_2(\Z)}
\newcommand{\pgl}{\mathrm{PGL}_2(\Z)}

\newtheorem{theorem}{Theorem}
\newtheorem{notatdef}[theorem]{DEFINITION-NOTATION}
\newtheorem{corollary}[theorem]{Corollary}
\newtheorem{lemma}[theorem]{Lemma}
\newtheorem{proposition}[theorem]{Proposition}

\usepackage{color} 

\newcommand\nt[1]{\textcolor{red}{{ #1}}}

\newcommand{\Jimm}{{\adjustbox{scale=.8, raise=1ex, trim=0px 0px 0px 7px}{\RL{j}}}}
\newcommand{\jimm}{{\adjustbox{scale=.5, raise=1ex, trim=0px 0px 0px 7px}{\RL{j}}}}
\newcommand{\waw}{{\raisebox{1mm}{\RL{w}}}}

\renewcommand{\Delta}{{\rm Discr}}

\newcommand{\sherh}[1]{\fboxsep=0pt\setlength{\fboxrule}{1pt}
\begin{center}
   \fbox{\colorbox{green}{
         \begin{minipage}[t]{13cm}
            #1
         \end{minipage}
      }
   }
\end{center}}
\newcommand{\sherhh}[1]{\fboxsep=0pt\setlength{\fboxrule}{1pt}
\begin{center}
   \fbox{\colorbox{yellow}{
         \begin{minipage}[t]{13cm}
            #1
         \end{minipage}
      }
   }
\end{center}}

\newcommand{\sherhhh}[1]{\fboxsep=0pt\setlength{\fboxrule}{1pt}
\begin{center}
   \fbox{\colorbox{red}{
         \begin{minipage}[t]{13cm}
            #1
         \end{minipage}
      }
   }
\end{center}}
\newcommand{\unut}[1]{#1}

\renewcommand{\sherh}[1]{}\renewcommand{\sherhh}[1]{}\renewcommand{\sherhhh}[1]{}\renewcommand{\nt}[1]{}

\usepackage{hyperref} 

\begin{document}

\title{Jimm, a Fundamental Involution}
\author{A. Muhammed Uluda\u{g}$^*$, Hakan Ayral\footnote{
{Galatasaray University, Department of Mathematics,}
{\c{C}{\i}ra\u{g}an Cad. No. 36, 34357 Be\c{s}ikta\c{s}}
{\.{I}stanbul, Turkey}}}

\maketitle

\centerline{ \it Dedicated to Yılmaz Akyıldız, who shared our enthusiasm about jimm}


\begin{abstract}
We study the  involution of the real line induced by the outer automorphism of the extended modular group PGL(2,Z). 
This `{modular}' involution is discontinuous at rationals but satisfies a surprising collection of functional equations. It preserves the set of real quadratic irrationalities mapping them in a highly non-obvious way to each other. It commutes with the Galois action on real quadratic irrationals.

More generally, it preserves set-wise the orbits of the modular group, thereby inducing an involution of the moduli space of real rank-two lattices. It induces a duality of Beatty partitions of the set of positive integers. 

This involution conjugates (though not topologically) the Gauss' continued fraction map to an intermittent dynamical system on the unit interval with an infinite invariant measure. \unut{The transfer operator (resp. the functional equation) naturally associated to this dynamical system is closely related to the Mayer transfer operator (resp. the Lewis' functional equation).}

We give a description of this involution as the boundary action of a certain automorphism of the infinite trivalent tree. We prove that its derivative exists and vanishes almost everywhere. It is conjectured that algebraic numbers of degree at least three are mapped to transcendental numbers under this involution.
\end{abstract}


\section{Introduction}\label{sec:introduction}
It is (it seems not very well-) known that  the group $\pgl$ has an involutive outer automorphism, which was discovered by Dyer in the late 70's \cite{dyer1978automorphism}. 
It would be very strange if this automorphism had no manifestations in myriad contexts where $\pgl$ or its subgroups play a major role. 
Our aim in this paper is to elucidate one of these manifestations, which appears to have remained in obscurity until now. This may be because  the conventional arithmetic, algebraic and geometric structures are not `{respected}' by this involution (i.e. it sends parabolics to hyperbolics) which makes its study all the more appealing.  Being thus exotic in several ways, we denote it - and some other involutions it gives rise to - by the letter\footnote{This is the fifth letter of the arabic alphabet in the hij\^a’\^i order. Latex preamble commands for a latin-compatible typography of $\Jimm$ are given at the end of the paper.} $\Jimm$ (read as: ``jimm"), hoping that this notation will help to keep track of its manifestations. 

Let $\widehat{\R}:=\R\cup \{\infty\}$. The manifestation in question of $\Jimm$ is a map $\Jimm_\R: \widehat{\R}\rightarrow \widehat{\R}$.
Denoting the continued fractions in the usual way
$$
[n_0,n_1,n_2,\dots]=n_0+\cfrac{1}{n_1+\cfrac{1}{n_2+\cfrac{1}{\dots}}},
$$
one has, for an irrational number $[n_0,n_1,n_2,\dots]$ with $n_0, n_1, \dots \geq 2$,
\begin{equation}
\Jimm_\R([n_0,n_1,n_2,\dots])=[1_{n_0-1},2,1_{n_1-2},2,1_{n_2-2},\dots]
\end{equation}
where $1_k$ is the sequence ${1,1,\dots, 1}$ of length $k$.
This formula remains valid for $n_0, n_1, \dots \geq 1$, if the emerging $1_{-1}$'s are eliminated  in accordance with the rule $[\dots m, 1_{-1},n,\dots]=[\dots m+n-1,\dots]$ and $1_0$ with the rule 
$[\dots m, 1_{0},n,\dots]=[\dots m,n,\dots]$. See page \pageref{examples} below for some examples.

It is possible to extend this definition of $\Jimm_\R$ to all of $\widehat{\R}$. If we ignore rationals and the noble numbers (i.e. numbers in the $\pgl$-orbit of the golden section, $\Phi:=(1+\sqrt{5})/{2}$), then $\Jimm_\R$ becomes an involution.
It is well-defined and continuous at irrationals,  but two-valued and discontinuous at rationals (Theorem~\ref{jimmcontini}).  Although it is impossible to draw its graph, we will give (see page \pageref{jimmplot}) a boxed-graph to indicate where the graph lies. One can choose one among the two values at rational arguments so that the function becomes upper semicontinuous.  The amount of jump of $\Jimm_\R$ at $q\in \Q$ provides a canonical (signed) measure of complexity of a rational number.

\medskip\noindent
{\bf Guide for notation}. In what follows,\\
$\Jimm$ denotes Dyers' outer automorphism of $\pgl$,\\
$\Jimm_\F$ denotes the automorphism of the tree $|\F|$ induced by $\Jimm$,\\
$\Jimm_{\partial\F}$ denotes the homeomorphism of the boundary $\partial\F$ induced by $\Jimm_\F$,  \\
$\Jimm_\R$ denotes the involution of $\R$ induced by  $\Jimm_{\partial\F}$, and\\
$\Jimm_\Q$ denotes the involution of $\Q\setminus\{0,\infty\}$ induced by  $\Jimm_{\F}$.\\
However, we reserve the right to drop the subscript and simply write $\Jimm$ when we think that confusion won't arise.

\subsubsection*{Functional equations} 
The involution $\Jimm_\R$ shares the privileged status of the fundamental involutions generating the extended modular group 
$\pgl$,
$$
K:x\rightarrow 1-x, \quad U:x\rightarrow 1/x, \quad V:x\rightarrow -x,
$$ 
as it interacts in a very harmonious way with them. 
Indeed, if we ignore its values at rational points then $\Jimm_\R$ satisfies the following set of functional equations:

\begin{center}
\fbox{\fbox{\begin{minipage}{13cm}

\vspace{2mm}
$$
\mbox{(FE:0)} \quad \Jimm(\Jimm(x))=x
$$

$$
\mbox{(FE:I)}\quad \Jimm U=U\Jimm \iff 
\Jimm\bigl(\frac{1}{x}\bigr)=\frac{1}{\Jimm(x)}
$$

$$
\mbox{(FE:II)} \quad\Jimm V= UV\Jimm\iff 
\Jimm(-x)=-\frac{1}{\Jimm(x)}
$$

\vspace{2mm}
$$
\mbox{(FE:III)} \quad \Jimm K=K\Jimm
\iff \Jimm(1-x)=1-\Jimm(x)
$$

\vspace{1mm}
\end{minipage}}}
\end{center}

The following set of functional equations are derived from the above ones:

\begin{center}
\fbox{\fbox{\begin{minipage}{13cm}

$$
\mbox{(FE:II$'$)} \quad\Jimm UV= V\Jimm\iff 
\Jimm(-\frac{1}{x})=-{\Jimm(x)}
$$
\vspace{2mm}
$$
\mbox{(FE:IV)} \quad \Jimm KV=KUV \Jimm
\iff \Jimm(1+x)=1+\frac{1}{\Jimm(x)}
$$

\vspace{2mm}
$$
\mbox{(FE:V)} \quad 
(\Jimm_\R M \Jimm_\R) (x)=(\Jimm M)(x) \quad \forall M\in \pgl$$

\vspace{2mm}
$$
\mbox{(FE:VI)} \quad 
\Jimm_\R (M x)=\Jimm M \Jimm_\R(x) \quad \forall M\in \pgl
$$

\vspace{1mm}
\end{minipage}}}
\end{center}

\bigskip
Equations (FE:I) and (FE:III) states that $\Jimm$ is covariant  with the operators $U$ and $K$,
whereas equation (FE:II) is a kind of lax-covariance\nt{\footnote{For a discussion of the functions covariant with the modular group in the literature, see the concluding remarks at page \pageref{concluding}.}.} Since $U$, $V$ and $K$ generate the group
$\pgl$, we may say that $\Jimm$ is lax-covariant with the $\pgl$-action on $\R$.
Of course, if we consider $\Jimm$ as an operator then covariance is simply  a relation of commutativity. Relation (FE:IV) is not independent from the rest as it can be easily deduced  from (FE:II) and (FE:III). The final relation is the most general form of the functional equations and says that $\Jimm_\R$ conjugates the M\"obius transformation $M\in \pgl$ to 
the M\"obius transformation $\Jimm(M)\in \pgl$. It is readily deduced from (FE:I-III) by using the involutivity of $\Jimm$. As an instance of (FE:V), $\Jimm$ conjugates the translation $1+x$ to $1+1/x$ and the transformation $1/(1+x)$ to $x/(1+x)$. (FE:VI) is got from (FE:V) by 
setting $x\to\Jimm_\R(x)$. It says that $\Jimm_\R$ is a kind of modular function.

\sherhh{(FE:II$'$) is reminiscent of the momentous theta function identity 
$\theta(-1/z)=\sqrt{z/i}\,\theta(x)$ and indeed the function 
$f(x):=^4\!\!\!\sqrt{x}\,\Jimm(x)$ do satisfy the identity 
$f(-1/z)=-\sqrt{z/i} f(x)$.--\nt{this is not very significant since $f$ is not periodic.}
}

Before attempting to play with them, beware that these functional equations are not consistent on the set of rational numbers. \nt{must explain why it is not consistent.} The first thing to try are the noble numbers. They are sent to rationals under $\Jimm$ and $\Jimm$ is 2-to-1 on this set. On the other hand, there is a related involution $\Jimm_\Q$ on the set of positive rationals satisfying some  functional equations. It must be stressed that $\Jimm_\Q$ is not the restriction of $\Jimm_\R$ to ${\Q}$. See page \pageref{jimmq}.

The privileged status of $\Jimm$ is vindicated by the fact that it preserves the ``{real-multiplication locus}", i.e. the set of real quadratic irrationalities. It does so in a highly non-trivial manner, though it preserves setwise the $\pgl$-orbits of real quadratic irrationalities. More generally $\Jimm$ preserves setwise the $\pgl$-orbits on $\widehat{\R}$ thereby inducing an involution of the moduli space of real rank-2 lattices, $\widehat{\R}/\pgl$. Moduli of the lattices represented by the numbers $[\overline{1_k, k+2}]$ are fixed under this involution. For a precise description of all fixed points, see Proposition~\ref{fixed}.

One could use the functional equations (FE:I)-(FE:III) to directly define and study the involution $\Jimm_\R$. However, the elusive\nt{\footnote{We must admit that it was not easy to firmly establish the well-definedness of $\Jimm_\R$, one end of it always appears to be loose, and the safest way for us to think about it have been to see it as acting on the boundary of $\pgl$. This boundary can be viewed as the set of equality classes of convergent semiregular continued fractions, see \cite{kunming}.}} nature of $\Jimm_\R$ is best understood by considering it as a homeomorphism of the boundary of the Farey tree, induced by an automorphism of the tree. This automorphism is the one which twists all but one vertex of the tree.

Finally, the following two-variable consequence of the functional equations is noteworthy:
$$
\boxed{\frac{1}{x}+\frac{1}{y}=1\iff \frac{1}{\Jimm(x)}+\frac{1}{\Jimm(y)}=1}
$$
Hence $\Jimm$ sends harmonic pairs of numbers to harmonic pairs. See page \pageref{twovariable} for a complete set of two-variable functional equations. 

Recall that, if the pair $(x,y)$ satisfies ${1}/{x}+{1}/{y}=1$, then the so-called Beatty sequences 
$
(\lfloor nx\rfloor)_{n\geq 1}$, $(\lfloor ny\rfloor)_{n\geq 1}
$
gives a partition of the set of positive integers (Rayleigh's theorem).  
Therefore every Beatty partition of positive integers admits a $\Jimm$-dual partition
$
(\lfloor n\Jimm(x)\rfloor)_{n\geq 1}$,  $(\lfloor n\Jimm(y)\rfloor)_{n\geq 1}.
$
\sherh{The difference $\lfloor (n+1)x\rfloor-\lfloor nx\rfloor$ is a characteristic Sturmian word in the alphabet $\{\lfloor x\rfloor+1,\lfloor x\rfloor\}$.}
\subsubsection*{Some examples.} \label{examples}
Here is a list of assorted values of $\Jimm$. Recall that $\Phi$ denotes the golden section.
\begin{align}\label{jimmphi}
\Phi=[\overline{1}]\implies\Jimm(\Phi)=\infty=\Jimm\bigl(-{1}/{\Phi}\bigr),
\end{align}
where by $\overline{v}$ we denote the infinite sequence $v,v,v, \dots$, for any  finite sequence  $v$. 
From \ref{jimmphi} by using (FE:II) we find 
\begin{align}\label{jimmminusphi}
\Jimm(-\Phi)=\Jimm\bigl({1}/{\Phi}\bigr)=0.
\end{align}
Repeated application of (FE:IV) gives
\begin{align}\label{jimmenplusphi}
\Jimm(n+\Phi)=\Jimm([n+1,\overline{1}])={F_{n+1}}/{F_{n}}
\end{align}
where $F_n$ denotes the $n$th Fibonacci number. One has
$$
\Jimm([1_n,2,\overline{1}])=[n+1,\infty]=n+1
$$
For the number $\sqrt{2}$ we have something that looks simple
$$
\sqrt{2}=[1,\overline{2}]\implies\Jimm(\sqrt{2})=[\overline{2}]=1+\sqrt{2}
$$
but this is not typical as the next example illustrates:
$$
\Jimm({(3+5\sqrt{2})}/{7})=
\Jimm([1,\overline{2, 3, 1, 1, 2, 1, 1, 1}])=
[2,\overline{2,1,4,5}]
=\frac{-3+2\sqrt{95}}{7} 
$$
In a similar vein, consider the examples
$$
\sqrt{11}=[3;\overline{3, 6}] 
\implies \Jimm(\sqrt{11})=[\overline{1,1,2,1,2,1,1}]=
\frac{15+\sqrt{901}}{26},
$$
$$
\Jimm(-\sqrt{11})=-1/\Jimm(\sqrt{11})=-\frac{26}{15+\sqrt{901}}=\frac{15-\sqrt{901}}{26}.
$$
The last example hints at the following result
\begin{theorem}
The involution $\Jimm$ commutes with the conjugation of real quadratic irrationals; i.e. for every real quadratic irrational $\alpha$ one has
$
\Jimm({\alpha}^*)={\Jimm(\alpha)^*}.
$
\end{theorem}
{\it Proof.} Every real quadratic irrational $\alpha$ is the fixed point of some $M\in \psl$, the Galois conjugate $\alpha^*$  being the other root of the equation $Mx=x$. But then $\Jimm(M\alpha)=\Jimm(M)\Jimm(\alpha)=\Jimm(\alpha)$, i.e. $\Jimm(\alpha)$ is a fixed point of the equation $\Jimm(M)$, the other root being $\Jimm(\alpha)^*$. Finally
$$
M\alpha=\alpha\implies M\alpha^*=\alpha^*\implies \Jimm(\alpha^*)=\Jimm(M\alpha^*)=\Jimm(M)\Jimm(\alpha^*),
$$
so $\Jimm(\alpha^*)$ is also a fixed point of $\Jimm(M)$, i.e. it must coincide with $\Jimm(\alpha)^*$. \hfill $\Box$.

\sherh{{\bf Jimm and the absolute Galois group}
 Today we understand that this result is better expressed in terms of the modular group: The automorphism group of the space
 $\mathbb{P}^1\backslash\{0,1,\infty\}$ is the symmetric group  $\Sigma_3$, and the quotient  
 ${\mathcal M}:=\mathbb{P}^1\backslash\{0,1,\infty\}/\Sigma_3$ is the modular curve (in fact an orbifold) ${\mathbb H}/PSL(2,\Z)$. Hence the fundamental group $\pi_1({\mathcal M})$ is the modular group $PSL(2,\Z)$, and the group ${\pi_1}(\mathbb{P}^1\backslash\{0,1,\infty\}\simeq F_2$
 is a subgroup of index 6 in this group. This gives rise to a more fundamental outer representation 
 \[\Gamma_{\mathbb{Q}}\to \mathrm{Out}(\widehat{\pi_1}({\mathcal M}))=\mathrm{Out}(\widehat{PSL(2,\Z)})\]
 
 Now, as you now the space $\mathbb{P}^1\backslash\{0,1,\infty\}$ has a further symmetry, which is the complex conjugation. 
 The quotient  
 ${\mathcal M}:=\mathbb{P}^1\backslash\{0,1,\infty\}/\Sigma_3$ is the modular tile ${\mathbb H}/PGL(2,\Z)$. This would give
 rise to an even more fundamental outer representation 
 \[\Gamma_{\mathbb{Q}}\to \mathrm{Out}(\widehat{\pi_1}({\mathcal M}))=\mathrm{Out}(\widehat{PGL(2,\Z)})\]
 And now the question is: whether if the element $\Jimm\in\mathrm{Out}(\widehat{PGL(2,\Z)})$ lies in the image of this representation, and if it does; if what is the element $\Gamma_{\mathbb{Q}}$ whose image is $\Jimm$.
 
 By above we know that it commutes with a big chunk. By the case, in case this action is not well-defined, one might extend the field 
 $\overline{\Q}$ by $\Jimm(\overline{\Q})$. Then $\Jimm$ acts on this thing for sure. However, we know that it is not a field automorphism when we view it as a map on $\R$.
}

To finish, let us give the $\Jimm$-transform of a  non-quadratic algebraic number, 
\begin{eqnarray*}
\Jimm(^3\!\sqrt{2})=
\Jimm([1; 3, 1, 5, 1, 1, 4, 1, 1, 8, 1, 14, 1, 10, 2, 1, 4, 12, 2, 3, 2, 1, 3, 4, 1, \dots])\\
=[2,1,3,1,1,1,4,1,1,4,1_6,3,1_{12},3,1_8,2,3,1,1,2,1_{10},2,2,1,2,3,1,2,1,1,3,\dots ]
\\
=2.784731558662723\dots
\end{eqnarray*}
and the transforms of two familiar transcendental numbers:
\begin{eqnarray*}
\Jimm(\pi)=
\Jimm([3, 7, 15, 1, 292, 1, 1, 1, 2, 1, 3,  \dots])=
[1_2, 2, 1_5, 2, 1_{13}, 3, 1_{290}, 5, 3,  \dots]\\
=1.7237707925480276079699326494931025145558144289232\dots
\end{eqnarray*}
\begin{eqnarray*}
\Jimm(e)=
\Jimm([2,1,2,1,1,4,1,1,6,1,1,8,\dots])= 
[1,3,4,1,1,4,1,1,1,1,\dots, \overline{4,1_{2n}}]\\
=1.3105752928466255215822495496939143349712038085627\dots
\end{eqnarray*}
We have been unable to relate these numbers to other numbers of mathematics.

\subsubsection*{Algebraicity and transcendence}
We made some numerical experiments  on the algebraicity of a few numbers $\Jimm(x)$ where $x$ is an algebraic number of degree $>2$, with a special emphasis on $\Jimm(^3\!\sqrt{2})$ (see our forthcoming paper). 
It is very likely that these numbers are transcendental. Furthermore, as we prove in Lemma~\ref{density}, if $X$ is a uniformly distributed random variable over the unit interval, then almost everywhere the average of continued fraction entries of $\Jimm(X)$ tends to 1, i.e. $\Jimm(X)$ does not obey the Gauss-Kuzmin statistics. As it is  widely believed and experimentally affirmed that the algebraic numbers of degree$\geq 3$ do obey the Gauss-Kuzmin statistics, one can state with confidence the following conjecture:

\bigskip\noindent
{\bf Conjecture.} If $x\in \R$ is algebraic of degree $>2$, then $\Jimm(x)$ is transcendental.

\bigskip
One philosophy concerning the continued fraction expansions of algebraic numbers is that it should not be possible to approximate them too closely by the rationals. On the other hand, since by the above remark on continued fraction entries, the convergents of the $\Jimm$-transform of a uniformly distributed number converges almost surely very slowly,  they can not be too closely approximated  by the rationals, and this seem to provide some (rather weak) evidence against the conjecture. 

It might be possible to tackle this conjecture with the methods recently introduced by Adamczewski and Y. Bugeaud, see 
\cite{Bugeaud}.

\unut{
\subsubsection*{Jimm-conjugates} 
Equation (FE:V) states that $\Jimm$ conjugates elements of $\pgl$ to elements of $\pgl$. In the next section we will see that $\Jimm$ conjugates the Gauss continued fraction map to something that makes sense. What about the other familiar functions of analysis? This question is probably not just a dull curiosity.

\bigskip\noindent
{\bf Conjecture.} If $M$ is an element of 
$\mathrm{PGL}_2(\Q)\setminus \pgl$, then the conjugation
$\Jimm M \Jimm$ assumes transcendental values at algebraic arguments of degree at least 3. 
}
\unut{
\bigskip
Here $\mathrm{PGL}_2(\Q)$ denotes the group of projective two-by-two {\it non-singular} integral matrices.

\nt{What about the conjugates $\Phi+x$ and $\Phi x$ etc.}}

\subsubsection*{Dynamics}
Denote by ${\mathbb T}_G$ the Gauss map $[0,1]\rightarrow [0,1]$, sending $x$ to the fractional part of $1/x$ and denote by 
${\mathbb T_F}$ the Farey map $[0,1]\rightarrow [0,1]$, defined as
\begin{equation}
{\mathbb T}_F(x)= 
\begin{cases}
{x\over 1-x}, &0\leq x\leq 1/2\\
{1-x\over x}, &1/2 < x\leq 1
               \end{cases}
\end{equation}
The involution $\Jimm$ sends the unit interval onto itself, and the $\Jimm$-conjugate of ${\mathbb T_F}$ is  ${\mathbb T_F}$ itself with the two branches being permuted. On the other hand, the involution $\Jimm$ conjugates ${\mathbb T}_G$ (but not topologically) to a self-map ${\mathbb T}_\jimm$ of the unit interval, defined by 
$$
{\mathbb T}_{\jimm}:
[0,1_k,n_{k+1},n_{k+2}, \dots] \in [0,1] \to 
[0,n_{k+1}-1,n_{k+2}, \dots] \in [0,1],
$$
where it is assumed that $n_{k+1}>1$ and $0\leq k <\infty$.
\unut{\begin{figure}[h]
\begin{center}\label{tjimm}
\noindent{\includegraphics[scale=.25]{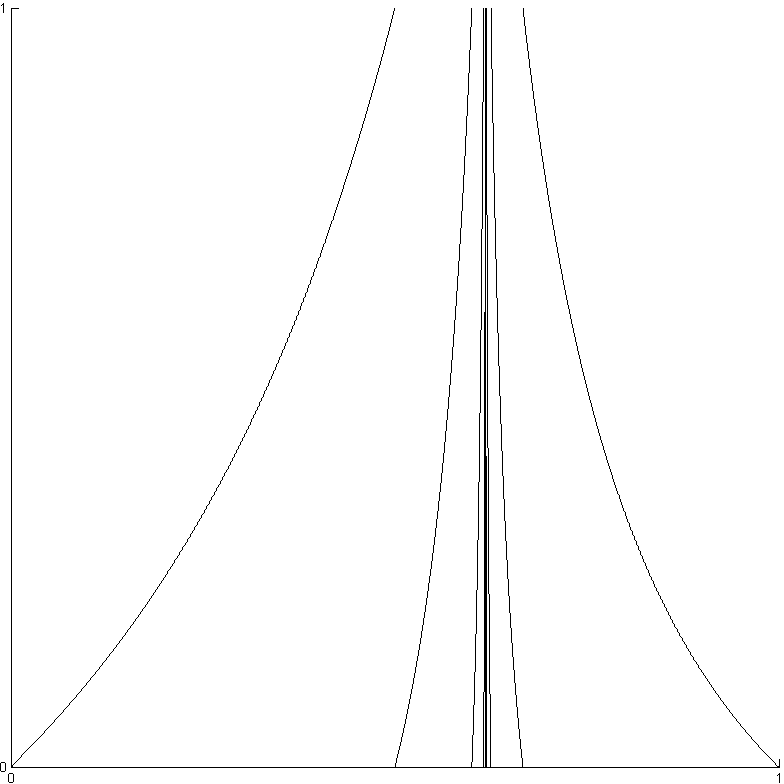}} \quad 
\noindent{\includegraphics[scale=.25]{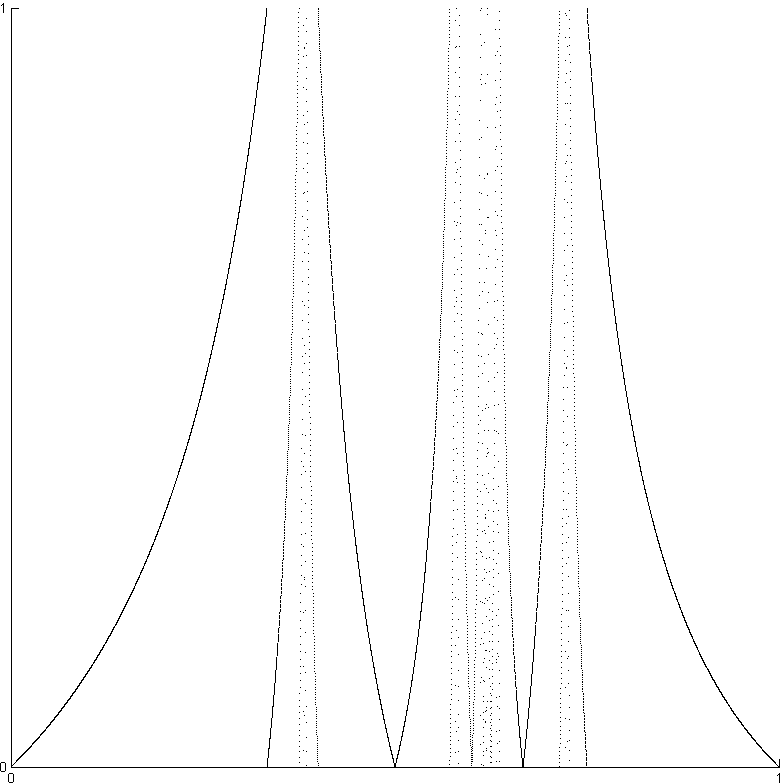}}
\end{center}
\vspace{-.3cm}
\caption{\small The graph of the map ${\mathbb T}_{jimm}$ and its first iterate.}
\end{figure}}
(This map appears also in two recent papers \cite{bonanno2014thermodynamic}, \cite{isola2014continued}, where the name ``Fibonacci map" was coined.) We discovered by trial and error that the resulting dynamical system have the infinite invariant measure $1/x(x+1)$, which can be readily verified \nt{(it can be also deduced from the correspondence (\ref{corr}) below). }
\unut{The map ${\mathbb T}_{\jimm}$ gives rise to a transfer operator 
(recall that $F_n$ is the $n$th Fibonacci number) 
\begin{equation}\label{jimmtransfer}
({\mathscr L}_s^{\jimm}\psi)(y)=\sum_{k=1}^\infty \frac{1}{(F_{k+1}y+F_{k})^{2s}}\psi\left(\frac{F_{k}y+F_{k-1}}{F_{k+1}y+F_{k}}\right)
\end{equation}
eigenfunctions of which satisfies the three-term functional equation
\begin{equation}\label{jimmfe}
\psi(y)=
\frac{1}{y^{2s}}\psi\left(\frac{y+1}{y}\right)+\frac{1}{\lambda}\frac{1}{(y+1)^{2s}}\psi\left(\frac{y}{y+1}\right)
\end{equation}
for the eigenvalue $\lambda$. It is straightforward to check that the solutions of this functional equation for $\lambda=1$
corresponds in a one-to-one manner to solutions of Lewis' functional equations \cite{zagier2001new}),
\begin{eqnarray}\label{lewisfe}
\phi(y)=\phi(y+1)+\frac{1}{(y+1)^{2s}}\phi\left(\frac{1}{y+1}\right) 
\end{eqnarray}
for the fixed points of the Mayer transfer operators~\cite{mayer1990thermodynamic}
\begin{eqnarray}\label{mayerto}
({\mathscr L}_s\phi)(y)=\sum_{k=1}^\infty \frac{1}{(y+n)^{2s}}\phi\left(\frac{1}{y+n}\right),
\end{eqnarray}
under the transformation 
\begin{equation}\label{corr}
\phi(y)=y^{-2s}\psi(1/y) \leftrightarrow \psi(y)=y^{-2s}\phi(1/y).
\end{equation}\nt{find the version of this with eigenvalues $\lambda$ other then 1.}
The correspondence (\ref{corr}) works also for the fixed functions of the operators (\ref{mayerto}) and (\ref{jimmtransfer}). 
Zagier and Lewis \cite{lewis2001period} established a correspondence between the Maass wave forms for the extended modular group and the solutions of the functional equation (\ref{lewisfe}).
Furthermore, according to a result of Mayer~\cite{mayer1990thermodynamic}, the  determinant of the operator $1-{\mathscr L}_s$, 
when made to act on the space of holomorphic functions on the disc $\{|z-1|\leq 3/2\}$, 
exists in the Fredholm sense and equals the 
Selberg zeta function of $\pgl$. Hence, our story is related to the Selberg zeta, although we don't expect an exact statement of Mayer's result to hold for the operator ${\mathscr L}_{s}^{\jimm}$. On the other hand, there is an alternative way of 
introducing some zeta analogues, by using the operator  ${\mathscr L}_{s}^{\jimm}$. 
Since $\zeta_{H}(2s,x)={\mathscr L}_s{\mathbf 1}(x)$ is the Hurwitz zeta function, the function
$$
\zeta_{\jimm}(2s, y):={\mathscr L}_{s}^{\jimm} {\mathbf 1}(y)=
\sum_{k=1}^\infty 
\frac{1}{(F_{k+1}y+F_k)^{2s}}
$$
arises as an analogue of the Hurwitz zeta and satisfies the functional equation
$$
\zeta_{\jimm}(s, y)=y^{-s}\zeta_{\jimm}(s, 1+1/y)+(y+1)^{-s}
$$
It reduces to the so-called ``{Fibonacci zeta}"  when $y=0$ 
$$
\zeta_{fib}(s):=\sum_{k=1}^\infty 
\frac{1}{F_k^{s}}=\zeta_{\jimm}(s, 0),\quad
$$
studied for its own sake in the literature \cite{murty2013fibonacci}. The $\zeta_{\jimm}$-values at $y=1$ and $y=2$ 
are also related to the Fibonacci zeta, via
$$
\zeta_{\jimm}(s, 1)=2+\zeta_{fib}(s), \quad \zeta_{\jimm}(s, 2)=2+2^{-s}+\zeta_{fib}(s),
$$
whereas $\zeta_{\jimm}$-value at $y=2$ is related to the so-called ``Lucas zeta" \cite{kamano2013analytic}:
$$
\zeta_{\jimm}(s, 2)=3^{-s}+\zeta_{luc}(s), \quad \zeta_{luc}(s):=\sum_{k=1}^\infty 
\frac{1}{(F_k+F_{k+2})^{s}}.
$$}
For more details about the dynamical system of ${\mathbb T}_{\jimm}$, its siblings and the associated zeta functions, see our forthcoming paper \cite{dynamicaljimm}. 

\unut{The involution $\Jimm$ is induced by an automorphism of the abstract Farey tree (denoted $|\F|$ in what follows).
This connection calls for a systematic study of dynamical properties of the conjugates of the Gauss map by automorphisms of 
$|\F|$ that fix an edge $I$, denoted $Aut_I(|\F|)$. This latter group  is an uncountable non-abelian profinite group (we give two quite concrete descriptions of its elements in this paper, in Theorems 2 and 3).  In the same vein, it is of interest to know about the dynamical properties which distinguish $Aut_I(|\F|)$-orbits of dynamical maps on the unit interval. }

\sherh{Is it true that the the set of piecewise-$\pgl$ maps (with rational break points) is preserved under the $Aut_I(|\F|)$-action ?}

\unut{It might also be of interest to study the $\Jimm$-conjugates of other euclidean dynamical systems \cite{vallee2006euclideans} (or of any other dynamical system with discontinuities at rationals for that matter) on $[0,1]$.}

\unut{\subsubsection*{Combinatorics and arithmetic}
Being an automorphism of $\pgl$, 
the automorphism $\Jimm$ acts on the system of subgroups 
$\mathbf{Sub}^*(\pgl)$ and on the system of conjugacy classes of subgroups $\mathbf{Sub}(\pgl)$.  We call the quotient
$\pgl\backslash {\mathcal H}$ the {\it modular tile} and denote by 
$\mathcal T$. This is an orbifold with boundary and can also be described as the quotient of the modular curve $\psl\backslash {\mathcal H}$ under the complex conjugation. The action of $\Jimm$ on subgroups (respectively conjugacy classes of subgroups) induces an action on the  base-pointed system of coverings $\mathbf{Cov}^*({\mathcal T})$
(respectively the system of coverings $\mathbf{Cov}({\mathcal T})$) of the modular tile. These coverings are surfaces possibly with boundary and with an orbifold structure. This action does not respect the genera of these surfaces. It does respect the rank of their fundamental groups. Jones and Thornton studied these actions to some depth in terms of the language of maps. One result they obtained (and besides Dyers', this is the only other result in the literature, about an action of $\Jimm$ that we are aware of) is

\medskip\noindent
{\bf Theorem} (Jones and Thornton \cite{jones1986automorphisms})  
{\it Let $G$ be a congruence subgroup of $\pgl$ such that $\Jimm(G)$ is also a congruence subgroup. Then $G\geq \Gamma(600)$, where $\Gamma(N)$ the level-$N$ principal congruence subgroup of $\psl$.}

\medskip  
Since these covering systems are naturally equivalent to systems of (half-) ribbon graphs, we have the {\it combinatorial action} of $\Jimm$ on these graphs. This action can be related to an action on combinatorial objects called necklaces,  bracelets,  Lyndon words, to the indefinite integral binary quadratic forms of Gauss and also to an action on geodesics on the modular curve. 
This latter action does not respect nor preserve the geodesic length. It does preserve the primitivity. The $\Jimm$-action on real quadratic irrationals studied in this paper to some depth, is related to these actions but we detail these connections elsewhere. }

\sherh{
\subsubsection*{Group Theory}
I have a vague feeling that $\Jimm$ is the outer automorphism group of $Aut(T)$ where $T$ is the Thompson group ; i.e. $\Jimm$ acts on $Aut(T)$ and together they define a subgroup of the homeomorphism group of the Farey tree. There will be a new set of functional equations concerning the interaction of $\Jimm$ with the elements of Thompson's groups.

But this might definitely be the case if $Aut(T)$ is the group of piecewise-$\pgl$ homeomorphisms of the circle. In any case we definitely have $T={\mathrm PPSL}_2(\Z) < {\mathrm{PPGL}_2(\Z)$, and this latter group appears to be acted upon by $\Jimm$. If you can find a presentation for ${\mathrm{PPGL}_2(\Z)}$ then it might be possible to pin down a presentation of the extension of $\mathrm{PPGL}_2(\Z)$ by $\Jimm$.}}

\unut{\subsubsection*{Numerics} 
There are many numerical excursions to be made around $\Jimm$. One important challenge which we failed until now, is to recognize 
(i.e. express by some algebraic/analytic procedure other then via $\Jimm$ itself) the
$\Jimm$-transform of some number which is not rational and not a quadratic irrational. A second curiosity is to study the statistics of $\Jimm$-transforms of some numbers (such as $\pi$ or $^3\!\sqrt{2}$) and see if they are ``{normal}" in some proper sense.
If you want to do some numerical experiments yourself, we plan to make our matlab and pyton codes and some graphics available via  a tablet computer version of this paper and on our webpage.}

\section{Modular group and its automorphism group}
The {\it modular group} is the projective group $\psl$ of two by two unimodular integral matrices \cite{koruoglu}. It acts on the upper half plane by M\"obius transformations. It is the free product of its subgroups generated by $S(z)=-1/z$ and $L(z)= (z-1)/z$, respectively of orders 2 and 3. Thus
$$
\psl =\langle S, L \,|\, S^2=L^3=1\rangle
$$
The projective group $\pgl$ consists of two by two integral matrices of determinant $\pm 1$. \unut{Its M\"obius action on the sphere does not preserve the upper half plane, i.e. if $\det(M)=-1$ then $M$ exchanges the upper and lower half planes. There is a modification of this action which preserves the upper half plane, where $M\in \pgl$ acts as 
$$
M=
\left[\begin{matrix}
p&q\\
r&s
\end{matrix}\right]
\leadsto 
\begin{cases} 
\frac{pz+q}{rz+s}, & if \quad \det(M)=1\\
\frac{p\overline{z}+q}{r\overline{z}+s},& if \quad \det(M)=-1
\end{cases}
$$
This representation of $\pgl$ is called the {\it extended modular group}.
Here we are interested in the $\pgl$-action on the boundary circle of the upper half plane, so we do not need to distinguish between these actions. But we shall keep calling $\pgl$ ``{the extended modular group}".

Note that, whenever we have an action of $\pgl$, we may compose it with $\Jimm$ and get another $\pgl$-action. 
Although it is essentially different from the first action, this latter action will have exactly the same invariants as the first one.}

\medskip
Below is a list of elements of $\pgl$ that will be used in the text.
{\scriptsize $$
\begin{array}{|c|c|c|c|c|} \hline &&\\[-1ex]
\quad \mbox{\bf Symbol}\quad &\quad \mbox{\bf M\"obius map}\quad &\quad \mbox{\bf Matrix form} \quad \\[2ex]\hline &&\\[-2ex]
S=UV&-1/x&\left[\begin{matrix} 0&1\\-1&0 \end{matrix}\right]\\[2ex]\hline &&\\[-2ex]
L=KU &1-1/x&\left[\begin{matrix} 1&-1\\1 &0 \end{matrix}\right]\\[2ex]\hline &&\\[-2ex]
V&-x&\left[\begin{matrix} -1&0 \\ 0&1 \end{matrix}\right]\\[2ex]\hline &&\\[-2ex]
T=LS=KV&x+1&\left[\begin{matrix} 1&1\\0&1 \end{matrix}\right]\\[2ex]\hline &&\\[-2ex]
\widetilde{T}=TU=KS&1+1/x&\left[\begin{matrix} 1&1\\1&0 \end{matrix}\right]\\[2ex]\hline &&\\[-2ex]
U=SV&1/x&\left[\begin{matrix} 0&1\\ 1 &0 \end{matrix}\right]\\[2ex]\hline &&\\[-2ex]
K&1-x&\left[\begin{matrix} -1&1\\ 0 &1 \end{matrix}\right]\\[2ex]\hline
\end{array}
$$}

\medskip
The fact that $Out(\pgl)\simeq \Z/2\Z$ has a story with a twist:
Hua and Reiner \cite{hua1952automorphisms} determined the automorphism groups of various projective groups in 1952, claiming that $\pgl$ has no outer automorphisms. The error was corrected by Dyer \cite{dyer1978automorphism} in 1978. This automorphism also appears in the work of Djokovic and Miller \cite{djokovic} from about the same time. Dyer also proved that the automorphism tower of $\pgl$ stops here; i.e. $Aut(Aut(\pgl))\simeq Aut(\pgl)$. Note that $\pgl\simeq Aut(\psl)$.

\sherh{Note that for us the story do not end here; embeddings of $\psl$ and $\pgl$ inside groups other then their automorphism groups also of interest for us. Non-normal embeddings.}

\medskip
Below is a list of presentations of $\pgl$ in terms of several sets of generators, and the outer automorphism $\Jimm$ of $\pgl$ defined in terms of these generators.

\medskip
{\scriptsize $$
\begin{array}{|l|l|} \hline &\\[-1ex]
\mbox{\bf Presentation of $\pgl$}\quad &\mbox{\bf The automorphism $\Jimm$}\quad  \\[2ex]\hline &\\[-.5ex]
\langle V,U,K\, |\, V^2=U^2=K^2=(VU)^2=(KU)^3=1 \rangle &(V,U,K) \rightarrow (UV,U,K)\\[2ex]\hline &\\[-.5ex]
\langle V,U,L\, |\, V^2=U^2=(LU)^2=(VU)^2=L^3=1 \rangle&(V,U,L) \rightarrow (UV, U, L)\\[2ex]\hline &\\[-.5ex]
\langle S,V,K\, |\, V^2=S^2=K^2=(SV)^2=(KSV)^3=1 \rangle\quad&(S,V,K) \rightarrow (V,S,K)\\[2ex]\hline &\\[-.5ex]
\langle T, U\, |\, { U^2= (T^{-2}UTU)^2=(UTUT^{-1})^3=1 } \rangle&(T,U) \rightarrow (\widetilde{T},U)\\[2ex]\hline &\\[-.5ex]
\langle T, \widetilde{T}\, |\, (T^{-1}\widetilde{T})^2= (T^{-3}\widetilde{T}^2)^2=
(T^{-2}\widetilde{T}^2)^3=1   \rangle&(T,\widetilde{T}) \rightarrow (\widetilde{T},T)\\[2ex]\hline 
\end{array}
$$}
\sherhh{With  $P:=\widetilde{T}^{-1}$ we get 
$$
\langle T, P\, |\,  (TP)^2= (T^3P^2)^2=
(T^2P^2)^3=1
\rangle \quad (T,P) \rightarrow (P^{-1},T^{-1}) 
$$
Not perfectly symmetric. To be verified.
}
\sherh{\begin{center}\label{jimmplot}
\noindent{\includegraphics[scale=1]{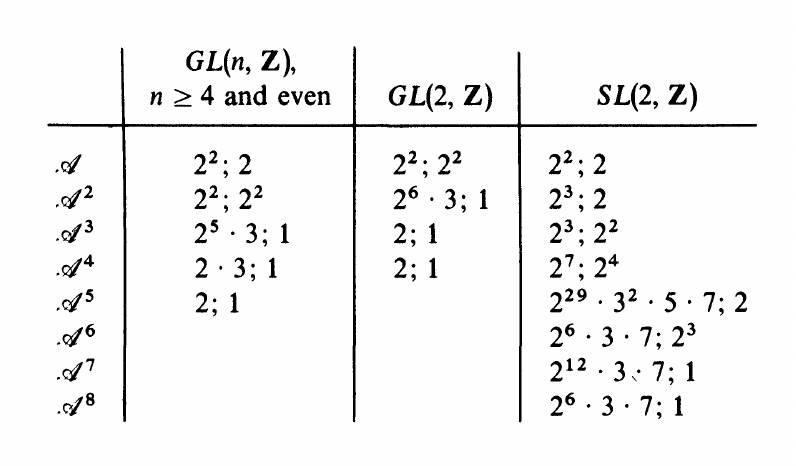}}\\
{\small {\bf Table.} Automorphism towers (from Dyers' paper \cite{dyer1978automorphism}). }\\
\end{center}}

\medskip\noindent
From the first presentation we see that Klein's viergruppe and the symmetric group on three letters are subgroups of $\pgl$:
$$
\Z/2\Z\times \Z/2\Z  \simeq\langle V,U\,|\, V^2=U^2=1 \rangle < \pgl, 
$$ 
$$
\Sigma_3\simeq \langle U,K\,|\, K^2=U^2=(UK)^3=1 \rangle < \pgl.
$$
In fact $\pgl$ is an amalgamated product of these subgroups \cite{dicks}:
$$
\pgl\simeq 
\langle V,U\rangle *_{\langle U\rangle } \langle V,K\rangle 
\simeq
(\Z/2\Z\times \Z/2\Z)  *_{\Z/2\Z} \Sigma_3
$$
In this description, we see that the outer automorphism of $\pgl$ originates from the following automorphism of Klein's viergruppe:
$$
\Jimm: (I, U, V, UV) \rightarrow (I, U, UV, V)
$$
We also see that $\Jimm$ leaves invariant the subgroup 
$$
\langle U, K \rangle \simeq \left\{z, \frac{1}{z}, 1-z, 1-\frac{1}{z},\frac{1}{1-z}, \frac{z}{z-1} \right\} \simeq\Sigma_3.
$$ 
It is easy to find the $\Jimm$ of an element of $\pgl$ given as a word in one of the presentations listed above. On the other hand, there seems to be no algorithm to compute the $\Jimm$ of an element of $\pgl$ given in the matrix form, other then actually expressing the matrix in terms of one of the presentations above, finding the $\Jimm$-transform, and then computing the matrix. 

The trace $\Jimm(M)$ has no simple expression in terms of the trace of $M$. It seems that $tr(\Jimm(M))$ is a novel and subtle class invariant of 
$M\in \pgl$ and of the binary quadratic form obtained by homogenization from the fixed point equation of $M$.

\nt{How strong is this invariant? Does $(tr(M), tr(\Jimm(M)))$ characterize the class? If yes might it be possible to express the Gauss product directly in terms of this pair?}

It is easy to see that the translation $2+x\in \psl$ is sent to the transformation $(2x+1)/(x+1)\in \psl$ under $\Jimm$, 
i.e. it may send a parabolic element inside $\psl$ to a hyperbolic element of inside $\psl\subset \pgl$.
Here is a selecta of other examples illustrating the effect of $\Jimm$ on matrices:
{\scriptsize 
$$
\begin{array}{|cc|cc|}
\hline&&&\\ [-1ex]
M&\Jimm(M)&M&\Jimm(M)\\[1ex]\hline &&&\\[-1ex]
\left[\begin{matrix}1& 0\\1& 1\end{matrix}\right] &
\left[\begin{matrix}0& 1\\1& 1\end{matrix}\right] &
\left[\begin{matrix}16& 1\\15& 1\end{matrix}\right] &
\left[\begin{matrix}987& 1597\\377& 610\end{matrix}\right]\\ [2ex]
\left[\begin{matrix}15& 1\\14& 1\end{matrix}\right] &
\left[\begin{matrix}610& 987\\233& 377\end{matrix}\right]&
\left[\begin{matrix}29& 2\\14& 1\end{matrix}\right] &
\left[\begin{matrix}843& 1364\\610& 987\end{matrix}\right]\\ [2ex]
\left[\begin{matrix}14& 1\\13& 1\end{matrix}\right] &
\left[\begin{matrix}377& 610\\144& 233\end{matrix}\right]&
\left[\begin{matrix}41& 3\\27& 2\end{matrix}\right] &
\left[\begin{matrix}665& 1076\\144& 233\end{matrix}\right]\\ [2ex]
\left[\begin{matrix}27& 2\\13& 1\end{matrix}\right] &
\left[\begin{matrix}521& 843\\377& 610\end{matrix}\right]& 
\left[\begin{matrix}40& 3\\13& 1\end{matrix}\right] &
\left[\begin{matrix}898& 1453\\521& 843\end{matrix}\right]\\ [2ex]
\hline
\end{array}
$$}
\unut{\paragraph{Presentation of $Aut(\pgl)$}
$$
\langle V,U,K, \Jimm\, |\, V^2=U^2=K^2=\Jimm^2=(VU)^2=(KU)^3=
(\Jimm V)^2U=(\Jimm K)^2=(\Jimm U)^2=1 \rangle 
$$
One can eliminate $U=(\Jimm V)^2$ from the presentation. Simplification gives
$$
\langle V,K, \Jimm\, |\, V^2=K^2=\Jimm^2=(\Jimm V)^4
=(K\Jimm V\Jimm V)^3=(\Jimm K)^2=1 \rangle 
$$

\paragraph{The largest Jimm-invariant subgroup.} 
The subgroup $\psl$ is not invariant under $\Jimm$, its image is the 
subgroup 
$$
\Jimm(\psl)\simeq \langle V, L \,|\, V^2=L^3=1 \rangle \simeq \langle -z, 1-1/z\rangle 
$$
Note that $\Jimm$ sends elements of order three in $\psl$ to elements of order three in $\psl$: since the determinant is multiplicative, there are in fact no elements of order three in $\pgl\setminus \psl$.

The largest $\Jimm$-invariant subgroup of $\psl$ is the index-2 subgroup $$
\psl\cap\Jimm \psl=\langle R, SRS \rangle = \Z/3\Z * \Z/3\Z
$$
This subgroup is the kernel  
$$
0\longrightarrow \langle R, SRS \rangle\longrightarrow \psl\longrightarrow \Z/2\Z\longrightarrow0
$$
Hence this is the subgroup of $\psl$ which consists of words in $L$ and $S$ with an even number of $S$'s; in particular this subgroup does not contain any elliptic elements of order 2.
Thus $\Gamma(2)$ and therefore all congruence subgroups 
$\Gamma(2N)$ are contained in it. On the other hand, since $T^{2N+1}\in\Gamma(2N+1)$, the subgroups $\Gamma(2N+1)$ are never contained in it.}

\section{The Farey tree and its boundary} 
\unut{\subsection{The bipartite Farey tree $\F$.}
(This subsection is taken to a large extent from a joint project of M. Uluda\u g with A. Zeytin on Thompson's groups.)}
From $\psl$ we construct the {\it bipartite Farey tree $\F$} 
tree on which $\psl$ acts, as follows.
The edges of $\F$ consists of the elements of the modular group; i.e.
$E(\F)=\psl$. The set of vertices of $\F$ are the left cosets of the subgroups $\langle S\rangle$ and $\langle L\rangle$.  Two distinct vertices $v$ and $v'$ are joined by an edge if and only if the intersection $v \cap v'$ is non-empty and in this case the edge between the two vertices is the only element in the intersection. Since cosets of a subgroup are disjoint, $\F$ is a bipartite graph. Cosets of $\langle S\rangle$ are always 2-valent vertices and the cosets of 
$\langle L\rangle$ are always 3-valent. 
The edges incident to the vertex $\{W,WL,WL^{2}\}$ are $ W, \, WL$ and $ WL^{2} $, and these edges inherit a natural cyclic ordering which we fix for all vertices as $(W, WL, WL^2)$. This endows $\F$ with the structure of a ribbon graph. It is connected since $\psl$ is generated by $S$ and $L$ and is circuit-free since it is freely generated by these elements.
Hence $\F$ is an infinite bipartite tree with a ribbon structure.

\sherh{
{\bf Exercise.} Show the edges corresponding to some elliptic elements on the tree. Show also an element $M$ and the conjugate of that element $NMN^{-1}$ for some $N$.
}

$M\in \psl$ acts on $\F$ from the left by ribbon graph automorphisms 
by sending the edge labeled $W$ to the edge labeled  $MW$. This action is free on $E(\F)$ but not free on the set of vertices: the vertex $\{W,WL,WL^{2}\}$ is left invariant by the order-3 subgroup $\langle WLW^{-1}\rangle$ and the vertex $\{W,WS\}$ is left invariant by the order-2 subgroup $\langle WSW^{-1}\rangle$.

\sherh{

\begin{center}
\noindent{\includegraphics[angle=90, scale=1]{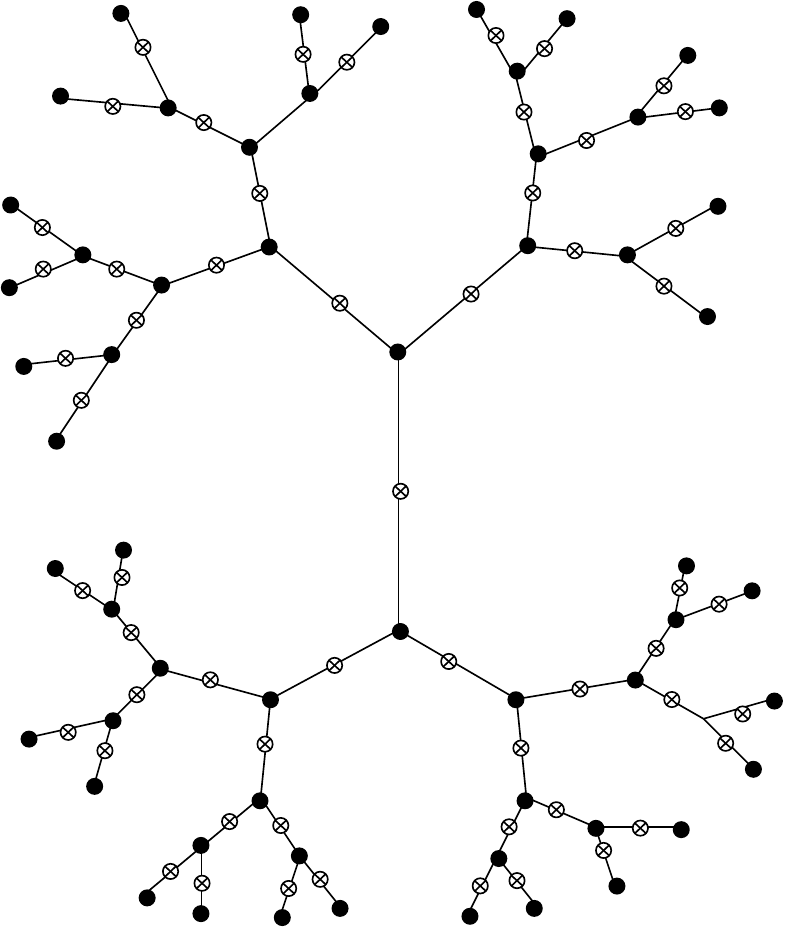}}\\
{\small {\bf Figure.} The tree. }\\
\end{center}

}

\paragraph{The boundary of $\F$.} 
\label{sec:boundary/of/F}
A {\it path} on a graph $\mathcal G$ is a sequence of edges $e_1,e_2,\dots,e_k$ of $\mathcal G$ such that $e_i$ and $e_{i+1}$ meet at a vertex, for each $1\leq i< k-1$.  Since the edges of $\F$ are labeled by reduced words in the letters $L$ and $S$, a path in $\F$ is a sequence of reduced words $(W_i)$ in $L$ and $S$, such that $W_i^{-1}W_{i+1}\in \{L, L^2, S\}$ for every $i$.
Since  $\F$ is a tree, there is a unique non-back\-tracking path through any two edges.

An {\it end} of $\F$ is an equivalence class of infinite (but not bi-infinite) non\--back\-tracking paths in $\F$, where eventually coinciding paths are considered as equivalent. In other words, an end of $\F$ is the equivalence class of an infinite sequence of finite reduced words $(W_i)$ in $L$ and $S$ with $W_i^{-1}W_{i+1}\in \{L, S\}$ for every $i$, where sequences with coinciding tails are equivalent. 

The set of ends of $\F$ is denoted $\partial \F$. 
The action of $\psl$ on $\F$ extends to an action on the set $\partial \F$, the element $M\in \psl$ sending the path $(W_i)$ to the path $(MW_i)$. 

Given an edge $e$ of $\F$ and an end $b$ of $\F$, there is a unique path in the class $b$ which starts at $e$. Hence for any edge $e$, we may  identify the set $\partial \F$ with the set of infinite non-backtracking paths that start at $e$. We denote this latter set by $\partial\F_e$  and endow it with the product topology. This topology is generated by the open sets, called {\it Farey intervals} ${\mathcal O}_{e'}$, which are defined to be the set of infinite paths starting at $e$ and passing through $e$. The space $\partial\F_e$ is also endowed with a natural cyclic ordering induced by the ribbon structure of $\F$. Hence $\partial\F_e$ is a cyclically ordered topological space.  Given  a second edge $e'$ of $\F$, the spaces $\partial\F_e$ and $\partial\F_{e'}$ are  homeomorphic under the order-preserving map which pre-composes with the unique path joining $e$ to $e'$. 
The action of the modular group on the set $\partial\F$ induces an action  of $\psl$ by order-preserving homeomorphisms of the topological space $\partial\F_e$, for any choice of a base edge $e$.

This construction of the Farey tree and that of its boundary below, was inspired by Qaiser Mushtaq's work \cite{mushtaq1}, \cite{mushtaq2} on coset diagrams, and also by a paper of Manin and Marcolli \cite{marcolli}. It can be established for a large class of trees with a planar structure, see \cite{tugce} and \cite{northshield}.

\paragraph{The continued fraction map.} 
$\partial\F_e$ is homeomorphic to the Cantor set, i.e. it is an uncountable, compact, totally disconnected, Hausdorff topological space.
By exploiting its cyclic order structure, we ``smash the holes" of this Cantor set to obtain the continuum, as follows. Define a {\it rational end} of $\F$ to be an eventually left-turn or eventually right-turn path. Now introduce the equivalence relation $\sim_*$ on $\partial\F$ as: left- and right- rational paths which bifurcate from the same vertex are equivalent. On $\partial\F_e$ this equivalence sets equal those points which are not separated by (with respect to the order relation) a third point.

\begin{center}
\noindent{\includegraphics[scale=.8]{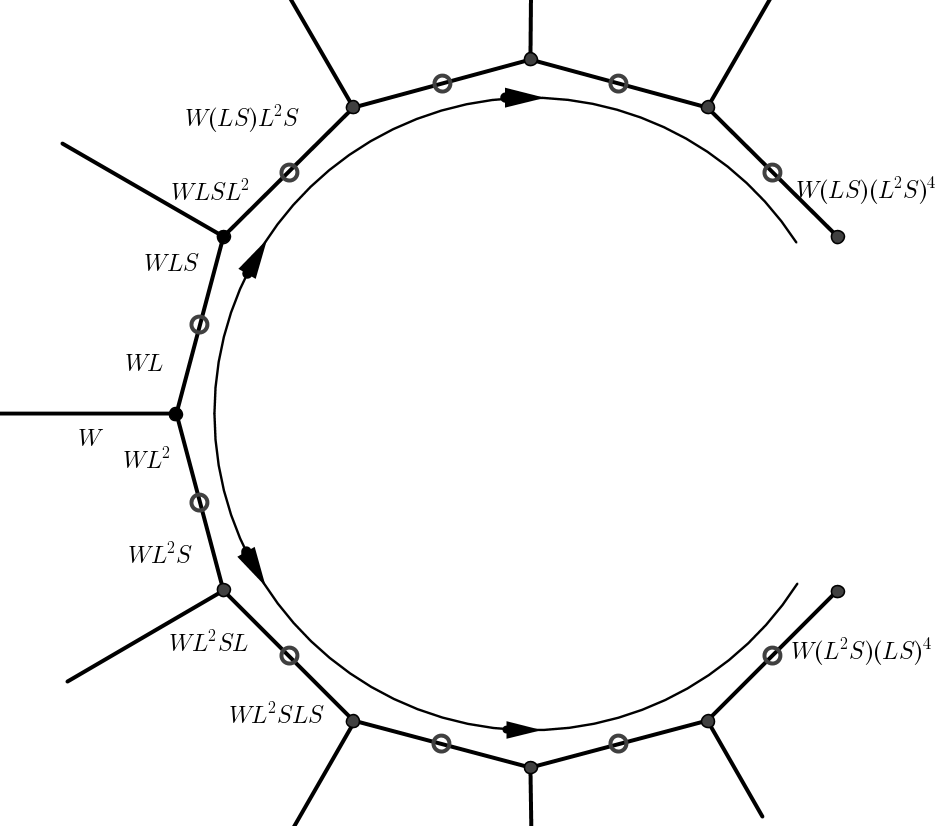}}\\
{\small {\bf Figure.} A pair of rational ends.}\\
\end{center}

On the quotient space $\partial\F_e/\!\sim_*$ there is the quotient topology induced by the topology on $\partial\F_e$ such that the projection map 
\begin{equation}\label{projectionmap}
 \partial\F_e\longrightarrow \partial\F_e/\!\sim_* 
 \end{equation}
  is continuous. The quotient space\footnote{The space $S^1_e:=\partial\F_e/\!\sim_*$ is called the circle boundary of $\F$, Northshield \cite{northshield} gave a much general construction for planar graphs, in the context of potential theory.} $\partial\F_e/\!\sim_*$ is a cyclically ordered topological space under the order relation inherited from $\partial\F_e$.  We shall denote this quotient space by $S^1_e$.  The equivalence relation is preserved under the canonical homeomorphisms $\partial\F_e \longrightarrow \partial\F_{e'}$ and is also respected by the $\psl$-action. Therefore we have the commutative diagram 
\begin{displaymath}
    \xymatrix{
\partial\F_e  \ar@{->}[r] \ar@{->}[d]  &  \partial\F_{e'}\ar@{->}[d]\\ 
S^1_{e}            \ar@{->}[r]                       &  S^1_{e'}\\
 }
\end{displaymath}
where the horizontal arrows are order-preserving homeomorphisms and the vertical arrows are projections. 
Moreover, $\psl$ acts by homeomorphisms on $S^1_{e}$, for any $e$.

Now, $\F$ comes equipped with a distinguished edge, the one marked $I$, the identity element of the modular group. Hence all spaces $S^1_{e}$ are canonically homeomorphic to $S^1_{I}$. 
\sherhh{One may prove this fact by using a topological characterization of the circle. However, below we shall prove this by specifying a canonical homeomorphism. See Theorem 1 below.}

Any element of $S^1_I$ can be represented by an infinite word in $L$ and $S$. Regrouping occurrences of $LS$ and $L^2S$, any such word $x$ of $S^1_I$ can be written in one of the following forms:
\begin{align}\label{reps}
	x&=(LS)^{n_0}(L^2S)^{n_1}(LS)^{n_2}(L^2S)^{n_3}(LS)^{n_4} \cdots \mbox{ or  }\\
	x&=S(LS)^{n_0}(L^2S)^{n_1}(LS)^{n_2}(L^2S)^{n_3}(LS)^{n_4} \cdots,
\end{align}

\noindent where $n_0, n_1 \dots \geq 0$. Since our paths do not have any backtracking we have $n_0\geq0$ and $n_{i}>0$ for $i = 1,2,\cdots$. 
The pairs of words 
\begin{align*}
	(LS)^{n_0}\cdots (LS)^{n_{k}+1} (L^2S)^{\infty} \mbox{ and  } &(LS)^{n_0}\cdots (LS)^{n_{k}}(L^2S) (LS)^{\infty}, & (k\mbox{ even}) \\
	(LS)^{n_0}\cdots (L^2S)^{n_{k}+1} (LS)^{\infty} \mbox{ and  } &(LS)^{n_0}\cdots (L^2S)^{n_{k}}(LS) (L^2S)^{\infty}, & (k\mbox{ odd})  
\end{align*}
correspond to pairs of rational ends and represent the same element of $S_I^1$. For irrational ends this representation is unique. 

The $\psl$-action on $S^1_I$ is then the pre-composition of the infinite word by the word in $L,S$ representing the element of $\psl$. In this picture it is readily seen that this action respects the equivalence relation
$\sim_*$

Set $T:=LS$, so that $T(x)=1+x$. Note that
\begin{eqnarray*}
(LS)^n.(L^2S)^m.(LS)^k(x)=(LS)^n.S.[S.(L^2S)^m .S]. S . (LS)^k(x)\\
=(LS)^n.S.[SL^2]^m. S . (LS)^k(x)=(LS)^n.S.[LS]^{-m}. S . (LS)^k(x)\\
=(x+n)\circ (-1/x) \circ (x-m) \circ (-1/x)\circ (x+k) =
n+\cfrac{1}{m+\cfrac{1}{k+x}}
\end{eqnarray*}
Accordingly, define the {\it continued fraction map} $cfm: S^1_I\rightarrow \hat {\mathbf R} $ by
$$
cfm(x)=\left\{ \begin{array}{rl}
[n_0,n_1,n_2,\dots]& \mbox{ if } x=(LS)^{n_0}(L^2S)^{n_1}(LS)^{n_2}(L^2S)^{n_3}(LS)^{n_4}\dots\\  
-1/[n_0,n_1,n_2,\dots]& \mbox{ if } x=S(LS)^{n_0}(L^2S)^{n_1}(LS)^{n_2}(L^2S)^{n_3}(LS)^{n_4}\dots
\end{array}\right.
$$\\[-10mm]
\begin{theorem}
	The continued fraction map $cfm$ is a homeomorphism.
	\label{thm:cfm}
\end{theorem}

\noindent{\it Proof.}
To each rational end $[n_0,n_1,n_2,\dots,n_k,\infty]$ we associate the rational number $[n_0,n_1,n_2,\dots,n_k]$. Likewise for the rational ends in the negative sector. This is an order preserving bijection between the set of equivalent pairs of rational ends and 
$\Q\cup\{\infty\}$.  Now observe that an infinite path is then no other than a Dedekind cut and conversely every cut determines a unique infinite path, see  \cite{uludag2012binary} for details.\hfill $\Box$

\bigskip
As a consequence of this result, we see that the continued fraction map conjugates the $\psl$-action on $S_I^1$ to its action on  
 $\widehat {\mathbf R}$ by M\"obius transformations.
Furthermore, there is a bijection between $\widehat {\Q}$ and the set of pairs of equivalent rational ends. 
Any pair of equivalent rational ends determines a unique {\it rational horocycle}, a bi-infinite left-turning (or right-turning) path.  
Hence, there is a bijection between $\widehat {\Q}$ and the set of rational horocycles. 

\nt{a figure here}

\begin{lemma} \label{correspondences} Given the base edge $I$ of $\F$,\\
(i) There is a natural bijection between the 3-valent vertices of $\F$ and $\widehat {\Q}\setminus\{0,\infty\}$.\\
(ii) There is a bijection between the 2-valent vertices and the Farey intervals $[p/q, r/s]$ with $ps-qr=1$.
\end{lemma}
{\it Proof.}
(i) On each rational horocycle there lies unique trivalent vertex which is closest to the base edge $I$. If we exclude the two horocycles on which $I$ lies, this gives a bijection between the set of trivalent vertices of $\F$ and the set of horocycles. The continued fraction map sends
the two horocycles through $I$ to $0$ and $\infty$. Hence, there is a natural correspondence between the set of trivalent vertices of $\F$ and 
$\widehat {\Q}\setminus \{0, \infty\}$. This correspondence sends the trivalent vertex $\{W, WL, WL^2\}$ to $W(1)\in \Q$, where it is assumed that $W$ is a reduced word which ends with an $S$. (ii) Every 2-valent vertex $\{W, WS\}$ lies exactly on two horocycles, and the set of paths based at $I$ and through $\{W, WS\}$ is sent to the interval $[W(0), WS(0)]$ under $cfm$. \hfill $\Box$
\sherhh{According to the lemma, there is a correspondence between the one-third of $\psl$ (i.e. the set of cosets of $\langle L\rangle$ = the set of trivalent vertices) and $\widehat {\Q}\setminus \{0, \infty\}$.

On the other hand, note that, given a point at infinity on the unit circle, the torus $\widehat{\R}^1\times \widehat{\R}^1$ can be viewed as the set of intervals on the circle, and it is possible to define a map 
$$
\Jimm: \widehat{\R}^1\times \widehat{\R}^1\to \widehat{\R}^1\times \widehat{\R}^1,
$$
sending intervals to intervals. If $\alpha<\beta$, then 
$$
[\alpha, \beta] \to [\lim_{x\to \alpha^+}\Jimm(x),  \lim_{x\to \beta^-}\Jimm(x) ]
$$
and
$$
[\beta, \alpha] \to [\lim_{x\to \beta^+}\Jimm(x),  \lim_{x\to \alpha^-}\Jimm(x) ]
$$
Note that this map is well-defined at each point (and this holds true if the above trick is applied to a function with only jump discontinuities).
}

\sherh{One may use this lemma to introduce a modular-group structure directly on the set of rationals (i.e. by concatenating paths).}

\paragraph{Periodic paths and the real multiplication set.}
Let $\gamma:=(W_1, W_2, \dots, W_n)$ be a finite path in $\F$. Then the {\it periodization} of $\gamma$ is the path $\gamma^\omega$ defined as
\begin{eqnarray*}
W_1, W_2, \dots, W_n,   W_nW_1^{-1}W_2,     W_nW_1^{-1}W_3,    \dots,    W_nW_{1}^{-1}W_n,  \qquad\qquad\qquad \\ 
\quad (W_nW_{1}^{-1})^2W_2, \dots, (W_nW_{1}^{-1})^2W_n, (W_nW_{1}^{-1})^3W_2, \dots
\end{eqnarray*}
In plain words, $\gamma^\omega$ is the path obtained by concatenating an infinite number of copies of a path representing $W_1^{-1}W_n$, starting at the edge $W_1$. If $W_1^{-1}W_n$ is elliptic, then $\gamma^\omega$ is an infinitely backtracking finite path. If not, 
$\gamma^\omega$ is actually infinite and represents an end of $\F$. We call these {\it periodic ends of $\F$}. 
Thus we have the periodization map
$$
per:\,\,   \{\,finite\, \, hyperbolic \, \, paths\,  \} \longrightarrow \{\, ends\, \, of\, \, \F\,  \}
$$
whose image consists of periodic ends. The map $per$ is not one-to one but its restriction to the set of primitive paths is one-to-one. 
The modular group action on $\partial\F$ preserves the set of periodic ends and the periodization map is $\psl$-equivariant.

Given an edge $e$ of $\F$ and a periodic end $b$ of $\F$, there is a unique path in the class $b$ which starts at $e$. This way the set of periodic ends of $\F$ is identified with the set of eventually periodic paths based at $e$. This set is dense in $\partial\F_e$  and preserved under the canonical homeomorphisms 
between the spaces $\partial\F_e$ and $\partial\F_{e'}$. Every periodic end has a unique $\psl$-translate, which is a purely periodic path based at $e$.
Finally, the set of periodic ends descends to a well-defined subset of $S^1_e$. The image of this set under the continued fraction map consists of the set of eventually periodic continued fractions, i.e. the set of real quadratic irrationalities (the ``real-multiplication set"). These are precisely the fixed points of the $\psl$-action on $S^1_e$; 
a hyperbolic element $M=(LS)^{n_0}(L^2S)^{n_1}\cdots (LS)^{n_{k}}\in \psl$ fixing the numbers represented by the infinite words
\begin{eqnarray*}
(LS)^{n_0}(L^2S)^{n_1}\cdots (LS)^{n_{k}}(LS)^{n_0}(L^2S)^{n_1}\cdots (LS)^{n_{k}}\dots, \mbox{ and}\\
(LS)^{-n_k}\cdots (L^2S)^{-n_1}(LS)^{-n_{k}}(LS)^{-n_k}\cdots (L^2S)^{-n_1}(LS)^{-n_{k}}\dots. 
\end{eqnarray*}
\sherh{{\bf Exercice.} Show that the fixed points are purely periodic if and only if $M$ is a cyclically reduced word.}
\subsection{The automorphism group of $\F$}
Recall that $\psl$ acts on  $\F$ by ribbon graph automorphisms and on $\partial \F$ by homeomorphisms respecting the pairs of rational ends thereby acting by homeomorphisms on $S^1$. \unut{Any non-identity element of $\psl$ is either of finite order, stabilizes a pair of rational ends,  or it stabilizes  a pair of eventually periodic ends of $\F$; in which case it is respectively called {\it elliptic}, 
{\it parabolic} or {\it hyperbolic}.  }

\sherh{\begin{lemma}
Elliptic elements in $Aut(\F)\simeq \psl$ are conjugate to either $S$, $L$ or $L^{-1}$. Parabolic elements are conjugate to a translation $T^n$ for some $n\in \Z$.
\end{lemma}}

Now let us forget about the ribbon structure of  $\F$.
This gives an abstract graph which we denote  by $|\F|$.
The automorphism group\footnote{The group $Aut(|\F|)$ is naturally isomorphic to the automorphism group of the abstract trivalent tree obtained from $|\F|$ by forgetting the vertices of degree 2.} $Aut(|\F|)$  of $|\F|$  is much bigger then $Aut(\F)$. It is uncountable and contains $Aut(\F)$ as a non-normal subgroup.
It is not compact but locally compact under its natural topology. 
In this topology, a neighborhood base of an automorphism
$\gamma$ consists those elements of $Aut(|\F|)$ which agree with $\gamma$ on finite subtrees.

The map
$$
Aut(|\F|) \times \partial \F \to \partial \F
$$
is continuous and $Aut(|\F|)$ also acts by homeomorphisms on 
$\partial \F$. However, automorphisms of $|\F|$ does not respect the ribbon structure on $\F$ in general and does not induce a well-defined homeomorphism of  $S_I^1$ in general (for example, as in the case of $\Jimm_{\partial \F}$, an automorphism of $|\F|$ may send a pair of rational ends to distinct irrational ends.).

\sherh{Note that these problematic points are at most countable in number.}

\unut{A bi-infinite path without back-tracking in $\F$ is called an {\it oriented geodesic} of $\F$; 
the set of geodesics of $\F$ is in one-to-one correspondence with the ``{Cantor torus}"
 $\partial \F\times \partial \F$, minus the diagonal. 
 In other words, every pair of distinct ends determine a unique oriented geodesic. }

\sherh{{\bf Remark.} $Aut_e(\F)$ is a proper subgroup of $Homeo(\partial\F)$. It is a rather small subgroup.
As an example, elements of $PGL_2(\Q)$ extends to homeomorphisms of the boundary, and they are not induced from the automorphisms, even those in $\psl$. Simply because they are continuous on the circle, and the elements of $Aut_e(\F)$ rarely induce continuous homeomorphisms on the boundary. 

Note that one can obtain many homeomorphisms of the boundary, by permuting the branches of the tree is some way. This is equivalent to permuting some disconnected components of the boundary in some way. This seem to correspond to a sort of completion of Thompson's group. This subgroup of homeomorphisms include those induced by the automorphisms of the tree. It might be equal to the whole group of homeomorphisms of the boundary. 

In particular, $x\in S^1_I \rightarrow 2x\in S^1_I$ lifts to a homeomorphism of $\partial\F$ but it is not a twist.\\

The lemma below is a variation of a lemma from \cite{cartier}.
\begin{lemma}
For  an automorphism  $\gamma$ of $|\F|$, there are two possibilities:
\begin{itemize} 
\item $\gamma$ does not stabilize any vertex.  
In this case there exists an oriented geodesic  $(e_i)$ 
and an $n\in \mathbf N$ such that $\gamma(e_i)=e_{i+n}$ for all $i$.
Hence, $\gamma$ stabilizes a pair of ends of $|\F|$.
Moreover, $\gamma$ is conjugate to a translation $T^n$ of  $\F$ for some $n\in\Z$. 
\item $\gamma$ stabilizes a vertex.  In this case there are two possibilities:
\begin{itemize} 
\item  $\gamma$ does not stabilize an edge. Then $\gamma$ is conjugate to a rotation $S$ or $L$ of $\F$. 
\item  $\gamma$ stabilizes  an edge. 
\end{itemize}
\end{itemize}
\end{lemma}

Now we elaborate on the last case. Evidently, if $\gamma$ stabilizes an edge $e$ 
 then it also stabilizes  the edge that meet $e$ at a vertex of degree 2. 
 
{\bf Exercice.} Show that the $PGL_2(\Q)$ action on the circle can be lifted to an action on the Cantor space $\partial\F$, by homeomorphisms. 
 }
 
 Fix an edge $e$ of $\F$ and denote by $Aut_e(|\F|)$  the group of automorphisms of $|\F|$ that stabilize  $e$. For any pair $e$, $e'$ of edges, $Aut_e(|\F|)$ and $Aut_{e'}(|\F|)$ are conjugate subgroups.
 
For $n>0$ let $|\F_n|$ be the finite subtree of $|\F|$ containing vertices of distance $\leq n$ from $e$. 
Then $|\F_{n}|\subset |\F_{n+1}|$ forms an injective system with respect to inclusion and $Aut_e(|\F_{n+1}|) \rightarrow Aut_e(|\F_{n}|)$ forms a projective system,
and one has 
$$
Aut_e(|\F|)=\lim_{\longleftarrow} Aut_e(|\F_{n}|).
$$
Hence $Aut_e(|\F|)$ is a profinite group.  Note that any edge $e$ splits $\F$ into two components and each one of these components are preserved by all elements of $Aut_e(|\F|)$.
Hence, 
$$
Aut_e(|\F|) \simeq Aut({\mathcal T}) \times Aut({\mathcal T}),
$$
where $\T$ is the rooted infinite binary tree. The group  $Aut({\mathcal T})$ 
might be described as a certain wreath product of an infinite number of copies of $\Z/2\Z$ (see Nekrashevych \cite{similar}), but we prefer and we shall present below two alternative descriptions which appears to be more intuitive.

\sherh{There is an injection $Aut(\F)\subset {\mathcal Homeo}(\partial\F)$, which is not surjective. For example,  there is a natural embedding (by lifting) $PGL_2(\Q)\to {\mathcal Homeo}(\partial\F)$.}

\sherh{If we omit the vertices of degree 2, then there is an automorphism which exchanges the two rooted trees, so that we have
$$
Aut_e(|\F|) \simeq (Aut({\mathcal T}) \times Aut({\mathcal T}))\ltimes \Z/2\Z.
$$
It can be shown that every element of finite order is in fact of order $2^n$ for some $n$. It must be possible to adopt some facts from  Nekrashevych's work, for example concerning automata}

\sherh{
\nt{The assertion below is almost but not entirely correct because to be precise one must assume that $T_n$ is bipartite, which brings forth some annoying details. Anyway, this is not needed for the sequel.}
Note that 
$$
Aut_e(|\F_n|) \simeq Aut({\mathcal T_{n+1}}) \times Aut({\mathcal T_n}),
$$
where ${\mathcal T_n}$ is the rooted binary tree of depth $n$.
}

\paragraph{A description of tree automorphisms as shuffles.} Let us turn back to the ribbon graph $\F$, which by construction comes with a base edge, the one labeled with $I\in\psl$. Given any vertex $v$, this base edge permits us to speak about the full subtree (i.e. Farey branch) attached to $\F$ at $v$.

Denote by $V_\bullet$ the set of vertices of degree three of $\F$.

The ribbon structure of $\F$ serves as a sort of coordinate system to describe all automorphisms of $|\F|$, as follows.

\sherh{{\bf Discussion.}
Note that the following data on $|\F|$ are equivalent:
\begin{itemize}
\item A planar embedding of  $|\F|$,
\item A faithful action of $\psl$ on  $|\F|$,
\item A ribbon graph structure on  $|\F|$,
\item A cyclic ordering of the edges of  $|\F|$,
\end{itemize}
There are uncountably many ways of endowing $|\F|$ with any of this data. 
Any choice of this data, together with the choice of a base edge, may be called a {\it frame} of $|\F|$. 
Thus our $\F$ comes equipped with a canonical frame by construction. The notion of a finitary automorphism depends on the frame. 

\medskip\noindent$\bullet$
There is yet another way of constructing $\F$ with its canonical frame, given as a cyclic ordering of its set of edges.

\medskip\noindent
$\bullet$ Is it true that, ``any maximal subgroup of $Aut(|\F|)$ acting without fixed edges is isomorphic (conjugate) to $\psl$''? 

 \medskip\noindent
$\bullet$ $\T$  can be constructed from the free monoid on two letters $\langle a,b\rangle^+$ by labeling its edges with the words in $\langle a,b\rangle^+$, such that any word $w$ has two children $wa$ and $wb$. The root vertex is labelled with the empty word, the unit element of the monoid.
The vertices of $\T$ can also be labeled with the elements of $\langle a,b\rangle^+$, except the root vertex - the vertex connecting the edge labeled $w$ to its children is labeled $w$.

\medskip\noindent
$\bullet$ Study the map sending the edge marked by $g$ to the edge marked by $g^{-1} $. \\
}

Given a vertex $v=\{W, WL, WL^2\}$ of type $V_\bullet$ of $\F$, the {\it shuffle} $\sigma_v$ is the automorphism of $|\F|$ which is defined as:
\begin{eqnarray}\label{twist}
\sigma_v: \mbox{edge labeled } M \longrightarrow 
\mbox{edge labeled} \left\{
\begin{array}{ll}
M, & \mbox{ if } M\neq WLX\\ 
WL^2X,& \mbox{ if } M=WLX\\
WLX,& \mbox{ if } M=WL^2X
\end{array}
\right.
\end{eqnarray}
where $M$ and $X$ are assumed to be reduced words in $S$ and $L$. Thus $\sigma_v$ is the identity away from the Farey branches at $v$, whereas it exchanges the two Farey branches at $v$.
Note that $\sigma_v^2=I$, i.e. the shuffle $\sigma_v$ is involutive.

\sherh{{\bf Exercice 1.} $\sigma_v$ and $\sigma_{v'}$ are conjugates in $Aut(|\F|)$ by an element of $\psl$. Similarly for $\theta_v$ and $\theta_{v'}$.
\\
{\bf Exercice 2.} Describe $\sigma_v$ as a map on the real line.
}

The automorphism $\sigma_v$ is obtained by permuting the two branches attached at $v$, by shuffling these branches one above the other. 
Beware that $\sigma_v$ is {\it not} the automorphism $\theta_v$ of $|\F|$ obtained by rotating  in the physical 3-space the branches starting at $WL$ and at $WL^2$ around the vertex $v$. We call this latter automorphism a {\it twist}, see below for a precise definition.

We must also stress that the definition of $\sigma_v$ 
requires the ribbon structure of $\F$ as well as a base edge, although $\sigma_v$  is never an automorphism of $\F$. 

Evidently, $\sigma_v$ stabilizes $I$, so one has $\sigma_v\in Aut_I(|\F|)$.

Given an arbitrary (finite or infinite) set $\nu$ of vertices in $V_\bullet$, we inductively define the shuffle 
$\sigma_\nu \in Aut_I(|\F|)$ as follows:
First order the elements of $\nu$ with respect to the distance from the base edge $I$, i.e. set 
$\nu^{(1)}=(v_1^{(1)}, v_2^{(1)},\dots)$, where $v_i^{(1)}\neq v_{i+1}^{(1)}$ and such that 
$$
d(v_i^{(1)},I) \leq d(v_{i+1}^{(1)},I) \mbox{ for all } i=1,2,\dots
$$
For $i=1,2\dots$ set
$$
\sigma^{(i)}_\nu:=\sigma_{v_i^{(i)}} \mbox{ and } v_j^{(i+1)}=\sigma^{(i)}_\nu(v_j^{(i)}) \mbox{ for } j\geq i+1
$$
Since for any $j>i$, the automorphism $\sigma^{(j)}_\nu$ agrees with  $\sigma^{(i)}_\nu$ on $\F_i$,
the sequence $\sigma^{(i)}_\nu$ converges in $Aut_I(|\F|)$ and we set
$$
\sigma_\nu:=\lim_{i\rightarrow \infty} \sigma^{(i)}_\nu
$$
for the limit automorphism of  $|\F|$. There is some arbitrariness in the initial ordering of $\nu$, concerning its elements of constant distance to the edge $I$, but the limit does not depend on this. The reason is that the shuffles corresponding to those elements commute. 
$\sigma_\emptyset$ is the identity automorphism by definition. Note that $\sigma_\nu$ is not involutive in general.

\sherh{{\bf Exercice }  Describe $\sigma_\mu$ as a map on the real line for some $\mu$.

Note that, it is not necessary to drag the flagged vertices inductively; absolute flags do just as good. Similarly for the twists. However, this choice is more natural, for example if $\gamma$ is an infinite path, then $\sigma_\gamma$ is something really about the path $\gamma$. 
It transforms the path $\gamma$ in an interesting way.. Whereas absolute flags don't do this.

{\bf Exercice }  Describe $\sigma_\gamma$ as a map on the real line for some $\gamma$.

{\bf Exercice }  Study the maps 
$$
\gamma\in \partial\F\to \sigma_\gamma\in Aut(\F)
$$
$$
\gamma\in \partial\F\to \sigma_\gamma\in {\mathcal Homeo}(\partial\F) 
$$
(more generally, there are maps $\partial\F^n\to \sigma_\gamma\in Aut(\F)$).

{\bf Exercice }  Find the order of the element $\sigma_\gamma$ for some $\gamma$.

{\bf Exercice } Do the same with twists.
}

This gives us a unique opportunity to glimpse inside an infinite, non-abelian profinite group;
any  automorphism of $|\F|$ that fix the edge $I$ is in fact a shuffle:
\begin{theorem}
$Aut_I(|\F|) =\{ \sigma_\nu\, |\, \nu\subseteq V_\bullet(\F)\}$.
\end{theorem}

\noindent{\it Proof.} It suffices to show that the restrictions of shuffles to finite subtrees $|\F_n|$ 
gives the full group $Aut_I(|\F_n|)$. This is easy.
\hfill$\Box$

Beware the trade-off: in the shuffle description, elements of $Aut_I(|\F|)$
are quite visible whereas the group operation is not so direct, as attempts to compute some powers of some non-trivial elements shows.

\sherh{{\bf Exercice} Describe some elements of order $n$ 
in $Aut_I(|\F|)$.

{\bf Remark.} This description permits us to define and study many special automorphisms of $|\F|$. For example, one may shuffle along a path or along a bi-infinite geodesic. Similarly for twists. Similarly for flips.

{\bf Question.} For any subgroup $\Gamma<\psl$, there is a shuffle $\sigma_\Gamma$ and there is a twist $\theta_\Gamma$. These will define some maps on $\R$. It is a very interesting question to study these maps. Can they be of finite order? What is their variation properties? Their derivatives and integrals? Do they conjugate the Gauss map to some interesting dynamical systems? What are their invariant measures and transfer operators? Zeta functions ?

}

\sherh{
\medskip
Recall that, upon a choice of a base edge, 
the set of ends $\partial\F$ is naturally identified with the set $\partial\F_I$, 
which carries the natural product topology.\\
{\bf Lemma}
$Aut_I(|\F|)$ acts transitively on $\partial\F_I$.\\
\noindent{\it Proof.} Given two ends $x,y\in \partial\F$, one can construct a convergent sequence of automorphisms sending longer and longer initial segments
of $x$ to initial segments of $y$. \hfill$\Box$

Recall that $\psl$ also acts on $\partial\F_I$ by homeomorphisms, so 
$\partial\F_I$ has many automorphisms not coming from an automorphism of  $|\F|$
that fix $I$.}

\sherh{
A {\it finitary automorphism} is a shuffle $\sigma_\nu$, where  $\nu$  is  a finite set of vertices in $V_\bullet$. The {\it group of finitary automorphism} is the group 
$$
Aut_I^{fin}(|\F|) =\{ \sigma_\nu\, |\, \nu\subset V_\bullet(\F), \quad |\nu|<\infty\}.
$$
$Aut_I^{fin}(|\F|)$ is not a normal subgroup. 
The notion of a finitary automorphism depends on the ribbon graph structure. Finitary subgroups defined by using different ribbon structures yield conjugate subgroups.

The above must be called the {\it finitary \underline{shuffle} group}. 
There is the notion of a finitary twist-automorphism as well. The co-finitary versions makes sense as well.

The group $Aut_I^{fin}(|\F|)$ acts by homeomorphisms on the boundary $\partial \F$.
However, shuffles (except $\sigma_{V_\bullet}$) do not respect pairs of rational ends, so $Aut_I^{fin}(|\F|)$ does not act on 
$S^1_I$ by homeomorphisms . However, they do preserve the set of rational ends and thanks to this it is possible to represent them as lower semicontinuous self-maps of the circle - or as interval exchange maps.

Finitary automorphisms gives rise to several interesting subgroups of $Homeo(\partial \F)$.
For example, there is the subgroup of $Aut(|\F|)$ generated by 
$Aut(\F)\simeq \psl$ and $Aut_I^{fin}(|\F|)$.

The group generated by finitary twists and flips act on the boundary $\partial \F$ as well.  

This automorphism $\theta_v$ is a local application of the automorphism $L\mapsto L^{-1}$ of $\psl$ and can be expressed as $\sigma_\nu$ where $\nu$ is the full subtree at vertex $v$ - to be confirmed. It is possible to describe its effect on $S^1$ explicitly. It is also possible to define $\theta_v$ inductively for subsets $\nu$, as in the case of $\sigma_\nu$. This way they give rise to a group of ``{cofinitary automorphisms}'' of $\F$.

One may also consider the subgroup generated by the finitary and cofinitary automorphisms of the tree. 

Are twists are co-finitary shuffles? Not really. Cofinitary shuffle group is really related to $\Jimm$.

Finitary automorphisms are dense in the full automorphism group..}

\paragraph{The automorphism $\sigma_{V_\bullet}$.}
What happens if we shuffle every $\bullet$-vertex of $|\F|$? 
In other words, what is effect of the automorphism $\sigma_{V_\bullet}$ on 
$\partial_I F$?
If $x\in \partial_I F$ is represented by an infinite word\footnote{In what follows, we will be sloppy about the difference between a path $x\in \partial_I\F$, its equivalence class in $S^1_I$, and the real number it represents, i.e. its image $cfm(x)$ under the continued fraction map.} in $L$, $L^2$ and $S$, then according to (\ref{twist}), we see that $\sigma_{V_\bullet}$
replaces $L$ by $L^2$ and vice versa. 
In other words, if $x$ is of the form 
$$	
x=(LS)^{n_0}(L^2S)^{n_1}(LS)^{n_2}(L^2S)^{n_3}(LS)^{n_4}\dots, 
$$
then one has
$$
\sigma_{V_\bullet}(x)=
(L^2S)^{n_0}(LS)^{n_1}(L^2S)^{n_2}(LS)^{n_3}(L^2S)^{n_4}\dots
$$
In terms of the contined fraction representations, we get, for $n_0>0$,
$$
\sigma_{V_\bullet}([n_0, n_1,n_2,\dots])=[0, n_1,n_2,\dots]
$$
and for $n_0=0$
$$
\sigma_{V_\bullet}([0, n_1,n_2,\dots])=[n_1,n_2,\dots].
$$
It is readily verified that a similar formula also holds when $x$ starts with an $S$, and we get
$$
\sigma_{V_\bullet}(x)=1/x=U(x).
$$ 
Observe (keeping in mind the forthcoming parallelism between $\Jimm$ and $U$) that $U$ is an involution satisfying the equations 
\begin{equation}\label{fesforu}
U=U^{-1}, \quad US=SU, \quad
LU=UL^2.
\end{equation}
Here, $U$, $L$ and $S$ are viewed as operators acting on the boundary $\partial_I\F$.
These equations can be re-written in the form of  functional equations as below:
$$
U(Ux)=x, \quad U(-\frac{1}{x})=-\frac{1}{Ux}, \quad U(\frac{1}{1-x})=1-\frac{1}{Ux}.
$$
One may derive other functional equations from these, i.e. $UT=ULS=L^2US=L^2SU$ is written as
$$
U(1+x)=\frac{Ux}{1+Ux}.
$$
One may consider $U$ as an ``{orientation-reversing}" automorphism of the ribbon tree $\F$. 
In the same vein, $U$ is a homeomorphism of 
$\partial_I\F$ which reverses its canonical ordering. It is the sole
element of $Aut_I(|\F|)$ which respects the equivalence $\sim_*$
(unlike $\Jimm$) and hence $U$ acts by homeomorphism on $S^1_I$.

\sherh{{\bf Exercice.} Prove the last claim}

\sherh{{\bf To do.} Show that $U$ is a limit of shuffles, and insert the animated gif (or a sequence of pictures) here.}

\sherh{
\begin{center}
\noindent{\includegraphics[scale=.25]{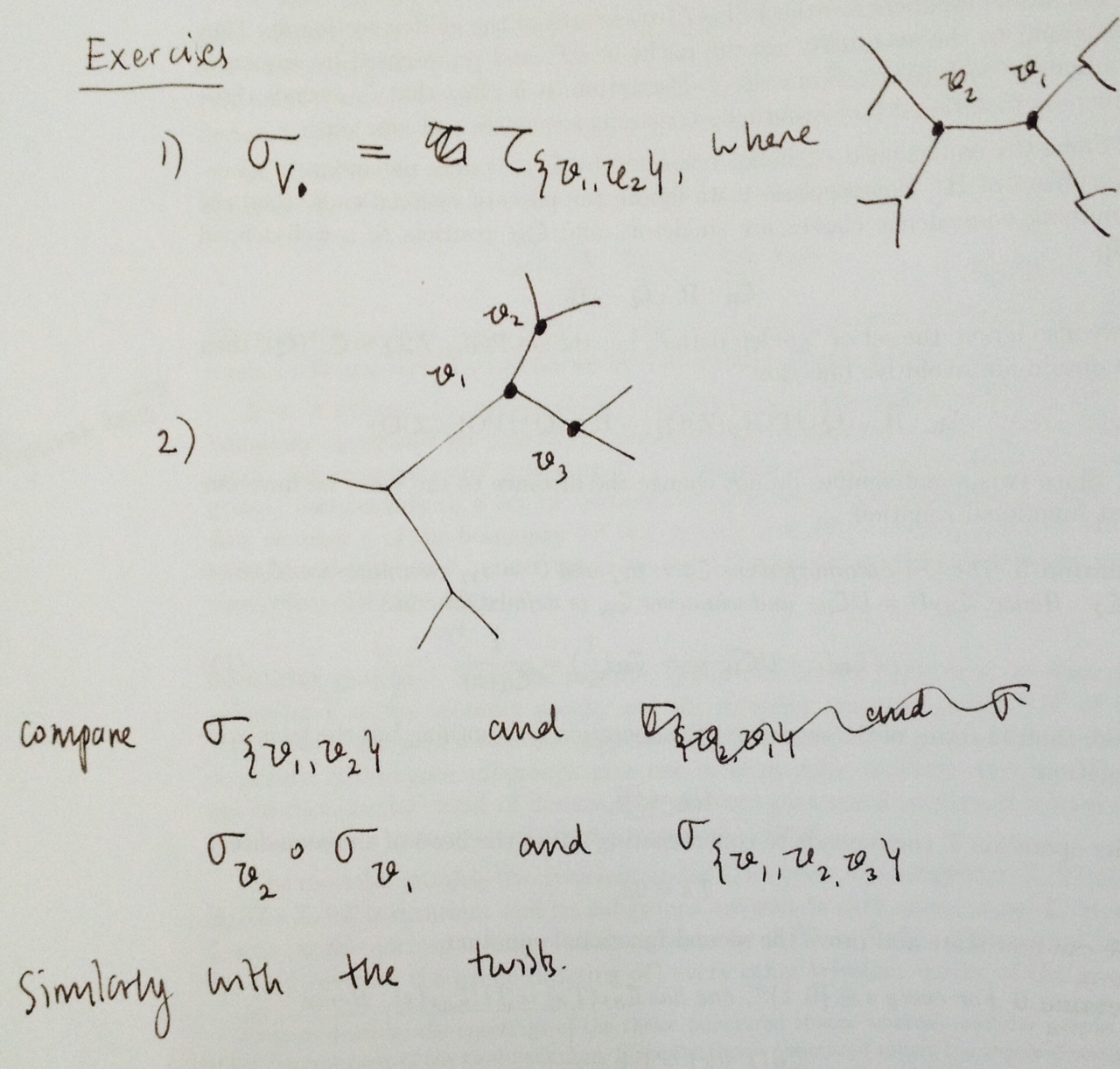}}\\
{\small {\bf Figure.}  }\\
\end{center}}

In fact, taken as a map on $E(\F)=\psl$, the involution $U$ is an automorphism of $\psl$, the one defined by 
$$
U:\left(\begin{matrix} S\\L\end{matrix}\right) \to 
\left(\begin{matrix} S\\L^2\end{matrix} \right)
$$
The automorphism group of $\psl$ is generated by the inner automorphisms and $U$: this is the group $\pgl$.

\sherh{{\bf Exercice.} Show
\begin{center}\label{jimmplot}
\noindent{\includegraphics[scale=.4]{./figures/reflect}}\\
{\small {\bf Figure.} }\\
\end{center}}

\sherh{({\bf Exercice.}
One has $\psl=\langle R, S\rangle$. The conjugation by $U$ gives $USU=S$ and $URU=R^2 \iff UR^2U=R$.
Hence, this action replaces all ``left-turns" by ``right-turns" and vice versa. In terms of the modular graphs, this action reverses the orientation of every vertex.
This action fixes (not element-wise) the congruence subgroups $\Gamma(N)$ and $\Gamma_0(N)$, since
$$
\left[\begin{matrix} 0&1\\1&0 \end{matrix}\right]\left[\begin{matrix} a&b\\c&d \end{matrix}\right]\left[\begin{matrix} 0&1\\1&0 \end{matrix}\right]
=\left[\begin{matrix} d&c\\b&a \end{matrix}\right]
$$
However, the congruence subgroups $\Gamma_1(N)$ are not fixed.

We may also consider this as an operation on the molecular space for $\pgl$. As such, its effect is nothing but a change of the base point.
}

\sherh{Let $(\nu_k)$ be a sequence of finite subsets of 
$V_\bullet(\F)$ with $\mu_k\subset \mu_{k+1}$ and $\cup \mu_k=V_\bullet(\F)$.
Then $\lim_{k\to\infty} \sigma_{\mu_k}$ is $V$. There is a similar fact for $\Jimm$. A canonical choice is to take $\nu_k=\F_{2k}$. This should suffice to describe $\Jimm$ as a limit of interval exchange maps.}

\paragraph{Twists.} Let $v\in V_\bullet$ and let 
$\nu_v$ be the set of all vertices (including $v$) on the Farey branch (with reference to the edge $I$) of $\F$ at $v$. Then the {\it twist} of $v$
 is the automorphism of $\F$ defined by 
$$
\theta_v:=\sigma_{\nu_v}.
$$
In words, $\theta_v$ is the shuffle of every vertex of the Farey branch at $v$.
As in the case of shuffles, for any subset $\mu\subset V_\bullet$, 
one may define the twist $\theta_\mu$ as a limit of 
convergent sequence of individual twists.

\begin{lemma}\label{twisttoshuffle}
Let $v$ be a vertex in $V_\bullet$ and $v'$, $v''$ its two children 
(with respect to the ancestor $I$). Put $\mu:=\{v,v',v''\}$.  
Then
$
\sigma_v=\theta_\mu.
$ 
\end{lemma}
Hence, shuffles can be expressed in terms of twists and vice versa. This proves the following result.
\begin{theorem}
$Aut_I(|\F|) =\{ \theta_\mu\, |\, \mu\subseteq V_\bullet(\F)\}$.
\end{theorem}

\begin{notatdef}\label{notatdef}
Let $v^*$ be the trivalent vertex incident to the base edge $I$, i.e. $v^*:=\{I,L,L^2\}$ and set 
$V_\bullet^*:=V_\bullet\backslash \{v^*\}$. 
We denote the special automorphism $\theta_{V_\bullet^*}$ 
by  $\Jimm_\F$.  In words, this 
is the automorphism of $|\F|$ obtained by 
twisting every trivalen vertex except the vertex $v^*$. This is the same automorphism obtained by 
shuffling every other trivalent vertex, such that the vertex $v^*$ is not shuffled.
\end{notatdef}
\sherh{The set of trivalent vertices of distance $4n$ to the base edge gives a subgroup, and these vertices are those remain non-shuffled under $\Jimm$. The set of trivalent vertices of distance $4n+2$ are shuffled. Thus, it is natural to shuffle subgroups and their cosets. }
More precisely, the following consequence of Lemma \ref{twisttoshuffle} holds:
\begin{lemma}
One has 
$
\theta_{V_\bullet^*}=\Jimm_\F=\sigma_J,
$
where $J\subset V_\bullet$ is the set of degree-3 vertices, whose distance to the base edge is an odd number. 
\end{lemma}
Obviously, $\Jimm_\F$ induce an involutive homeomorphism of $\partial_I\F$ 
which do not respect its canonical ordering.  In fact, one may say that it destroys the ordering of the boundary in the most terrible possible way. We denote this homeomorphism by $\Jimm_{\partial\F}$. Terrible as they are, we must emphasize that 
$\Jimm_{\F}$ and $\Jimm_{\partial\F}$  are perfectly well-defined mappings on their domain of definition.
They don't exhibit such things as the two-valued behavior of $\Jimm_{\R}$ at rationals. This two-valued behavior is a consequence of the fact that $\Jimm_{\partial\F}$ do not respect the equivalence relation $\sim_*$. 

To see the effect of $\Jimm_{\partial\F}$ on $x\in\partial_I\F$, assume 
$$	
x=S^\epsilon (LS)^{n_0}(L^2S)^{n_1}(LS)^{n_2}(L^2S)^{n_3}(LS)^{n_4}\dots, \, (n_0\geq 0, \, n_i>0 \mbox{ if } i>0, \,\epsilon\in\{0, 1\}).
$$
We may represent this element by a string of 0's and 1's (0 for $L$ and 1 for $L^2$):
\begin{eqnarray*}
x=S^\epsilon
\underbrace{00\dots0}_{n_0}
\underbrace{11\dots1}_{n_1}
\underbrace{00\dots0}_{n_2}
\underbrace{11\dots1}_{n_3}
\dots, \, 
\end{eqnarray*}
Then rational numbers are represented by the eventually constant strings where for any finite string $a$, the strings $a0111\dots$ and $a1000\dots$ represent the same rational number\footnote{The two representations of the number $0\in \R$ are
$1^\omega$ and $-0^\omega$ and the  two representations of $\infty\in \widehat\R$ are $0^\omega$ and $-1^\omega$.}.

Let $\phi:=(01)^\omega$ be the zig-zag path to infinity, and set $\phi^*:=\lnot \phi=(10)^\omega$.
Then 
$$
\Jimm_{\partial\F}(x)= \begin{cases} 
a \xor \phi,& x=a\in \{0,1\}^\omega,\\
S(a \xor \phi*),& x=Sa, \, a\in \{0,1\}^\omega.
\end{cases}
$$
where $\xor$ is the operation of term-wise exclusive or (XOR) on the strings of 0's and 1's. The involutivity of $\Jimm_{\partial\F}$ then stems from the reversibility of the disjunctive or: $(p\xor q) \xor q=p$.

\sherh{xoring (or performing term-wise logical operations) with a fixed infinite string defines boundary homeomorphisms 
induced by special kind of tree automorphisms, which can be called radially symmetric.}

\sherh{One may take a fixed bi-infinite string with a section, $\alpha^*|\alpha$ and xor with this 
$$
\Jimm_{\partial\F}= \begin{cases} 
a \xor \alpha,& x=a\in \{0,1\}^\omega,\\
S(a \xor \alpha*),& x=Sa, \, a\in \{0,1\}^\omega.
\end{cases}
$$
Example $\alpha^*|\alpha=(100)^*|(001)^*$. Something like this.
These boundary homeomorphisms comes from the tree automorphisms, those which shuffle the 1's on the string $\alpha^*|\alpha$. These are precisely the ``radially symmetric" automorphisms (-twists or -shuffles) of the tree.
}
\medskip\noindent
{\bf Some Examples.}
The string $0(0011)^\omega$ corresponds to the continued fraction 
$[1,2,2,2,\dots]$, which equals $\sqrt{2}$. One has
$$
\begin{array}{l|cccccccccccc|c}
\partial\F&&\multicolumn{10}{c}{\{0,1\}^\omega}&&\widehat{\R}\\
\hline
x 
&0&1&1&0&0&1&1&0&0&1&1&0\dots &\sqrt{2}\\
\phi
&0&1&0&1&0&1&0&1&0&1&0&1\dots&\Phi\\
\Jimm(x)
&0&0&1&1&0&0&1&1&0&0&1&1\dots&1+\sqrt{2}
\end{array}
$$
As for the value of $\Jimm_\R$ at $\infty$, one has
$$
\begin{array}{l|cccccccccccc|c}
\partial\F&&\multicolumn{10}{c}{\{0,1\}^\omega}&&\widehat{\R}\\
\hline
x 
&0&0&0&0&0&0&0&0&0&0&0&0\dots &\infty\\
\phi
&0&1&0&1&0&1&0&1&0&1&0&1\dots&\Phi\\
\Jimm(x)
&0&1&0&1&0&1&0&1&0&1&0&1\dots&\Phi
\end{array}
$$
The two values $\Jimm_\R$ assumes at the point $1$ are found as follows:
$$
\begin{array}{l|cccccccccccc|c}
\partial\F&&\multicolumn{10}{c}{\{0,1\}^\omega}&&\widehat{\R}\\
\hline
x 
&0&1&1&1&1&1&1&1&1&1&1&1\dots &1\\
\phi
&0&1&0&1&0&1&0&1&0&1&0&1\dots&\Phi\\
\Jimm(x)
&0&0&1&0&1&0&1&0&1&0&1&0\dots&1+\Phi
\end{array}
$$
$$
\begin{array}{l|cccccccccccc|c}
\partial\F&&\multicolumn{10}{c}{\{0,1\}^\omega}&&\widehat{\R}\\
\hline
x 
&1&0&0&0&0&0&0&0&0&0&0&0\dots &1\\
\phi
&0&1&0&1&0&1&0&1&0&1&0&1\dots&\Phi\\
\Jimm(x)
&1&1&0&1&0&1&0&1&0&1&0&1\dots&1/(1+\Phi)
\end{array}
$$
Conversely, one has $\Jimm_\R(1+\Phi)=\Jimm_\R(1/(1+\Phi))=1$, 
illustrating the two-to-oneness of $\Jimm_\R$ on the set of noble numbers. 

\bigskip
The noble numbers  are the $\psl$-translates of the golden section $\Phi$.
They correspond to the eventually zig-zag paths in $\partial \F$,
and represented by strings terminating with $(01)^\omega$. 
From the $\xor$-description, it is clear that $\Jimm_{\partial\F}$ sends those strings to rational (i.e. eventually constant) strings and vice versa.

Since the equivalence $\sim_*$ is not respected by $\Jimm_{\partial\F}$,
it does not induce a homeomorphism of  $\widehat{\R}$, not even a well-defined map. Nevertheless, 
if we ignore the pairs of rational ends, then the remaining equivalence classes are singletons and $\Jimm_{\partial\F}$ restricts to a well-defined map 
$$
\Jimm_\R:\widehat{\R}\setminus\widehat{\Q}
\to \widehat{\R}
$$
If we also ignore the set of ``{golden paths}", i.e. the set
$\pgl \widehat{\Q}=\Jimm^{-1}(\Q)$, then we obtain an involutive bijection
$$
\Jimm_\R: \widehat{\R}\setminus(\widehat{\Q}\cup \pgl \widehat{\Q})\to
\widehat{\R}\setminus(\widehat{\Q}\cup \pgl \widehat{\Q})
$$

\medskip\noindent
{\bf The functional equations.}
Since twists and shuffles do not change the distance to the base, 
we have our first functional equation:
\begin{lemma}
The $|\F|$-automorphisms $\Jimm_\F$  and $U=\sigma_{V_\bullet}$ commute, i.e. $\Jimm_\F U=U\Jimm_\F$. Hence,
$\Jimm_{\partial\F}U=U\Jimm_{\partial\F}$ and 
whenever $\Jimm_\R$ is defined, one has 
\begin{equation}\label{jimmu}
\Jimm_\R U= U\Jimm_\R \iff \Jimm_\R\bigl(\frac{1}{x}\bigr)=\frac{1}{\Jimm_\R(x)}
\end{equation}
\end{lemma}
In fact, $\Jimm_\F U$ is the automorphism of $\F$ which shuffles every other vertex, starting with the vertex $v^*$. 
In other words, it shuffles those vertices which are not shuffled by $\Jimm_\F$. We denote this automorphism by 
$\Jimm_\F^*$.
Note that, in terms of the strings, the operation $U$ is nothing but the term-wise negation:
$$
Ua=\lnot a,
$$ 
and the lemma merely states the fact that $\lnot (\phi \xor a)=\phi \xor \lnot a$. The boundary homeomorphism induced by the 
automorphism $\Jimm_\F^*$ is thus the map $\Jimm_{\partial \F}^*$ which xors with the string $\phi^*:=(10)^\omega$.

Now, consider the operation $S$. If $x=a\in \{0,1\}^\omega$, then one has $\Jimm_{\partial \F} (Sx)=$
\begin{eqnarray*}
 = S(a \xor \phi^*) 
=S(a\xor \lnot \phi) =S(\lnot(a\xor\phi))=S(U(a\xor\phi)) = SU\Jimm_{\partial \F} (x)=V\Jimm_{\partial \F} (x)
\end{eqnarray*}
The same equality holds if $x=Sa$, and we get our second functional equation:
\begin{lemma}
The $\partial_I\F$-homeomorphisms $\Jimm_{\partial\F}$  and $S$ satisfy $\Jimm_{\partial\F}S=V\Jimm_{\partial\F}$. 
Hence, whenever $\Jimm_\R$ is defined, one has 
\begin{equation}\label{jimmux}
\Jimm_\R S= V\Jimm_\R \iff \Jimm_\R\bigl(-\frac{1}{x}\bigr)=-{\Jimm_\R(x)}
\end{equation}
\end{lemma}

Now we consider the operator $L$. Suppose that $x=Sa$.
Then $Lx= LSa=0a$. Hence, noting that $\phi=(01)^\omega=0 (10)^\omega=0\phi^*$ we have
\begin{eqnarray*}
 \Jimm_{\partial \F} (Lx) = 0a \xor \phi=0a\xor 0\phi^*=0(a\xor \phi^*)= LS(a\xor \phi^*) 
 = L\Jimm_{\partial \F} (x).
\end{eqnarray*}
Another possibility is that $x=0a$.
In other words, $x$ starts with an $L$. 
Then $Lx$ starts with an $L^2$, i.e. $Lx= 1a$. Hence, 
\begin{eqnarray*}
 \Jimm_{\partial \F} (Lx) = 1a \xor \phi=1a\xor 0\phi^*=1(a\xor \phi^*),\\
  \Jimm_{\partial \F}(x)=0a \xor \phi=0a\xor 0\phi^*=0(a\xor \phi^*) 
\implies L \Jimm_{\partial \F}(x)= 1(a\xor \phi^*).
\end{eqnarray*}
Finally, if $x=1a$, then
$Lx$ starts with an $S$, i.e. $Lx= Sa$. Hence, 
\begin{eqnarray*}
 \Jimm_{\partial \F} (Lx) = S(a \xor \phi^*),\\
  \Jimm_{\partial \F}(x)=1a \xor \phi=1a\xor 0\phi^*=1(a\xor \phi^*) 
\implies L \Jimm_{\partial \F}(x)= S(a\xor \phi^*).
\end{eqnarray*}
Whence the third functional equation:
\begin{lemma}
The $\partial_I\F$-homeomorphisms $\Jimm_{\partial\F}$  and $L$ satisfy $\Jimm_{\partial\F}L=L\Jimm_{\partial\F}$. 
Hence, whenever $\Jimm_\R$ is defined, one has 
\begin{equation}\label{jimmux}
\Jimm_\R L= L\Jimm_\R \iff \Jimm_\R\bigl(1-\frac{1}{x}\bigr)=1-\frac{1}{\Jimm_\R(x)}
\end{equation}
\end{lemma}
(To be brief, the lemma holds true since $L$ rotates paths around the vertex $v^*=\{I,L,L^2\}$ and so it does not change the distance of an edge to that vertex.)

Since $U$, $S$ and $L$ generate the group $\pgl$, these three functional equations  forms a complete set of functional equations, from which the rest can be deduced. For example, 
$$
T=LS\implies \Jimm_{\partial \F}T=LV\Jimm_{\partial \F}\iff
\Jimm_\R(1+x)=1+\frac{1}{\Jimm_\R(x)}
$$

These functional equations also shows that $\Jimm_\R$ acts as the desired outer automorphism of $\pgl$. 

\sherh{Here is a direct proof for $T$:
In the positive sector, the operation $T$ concatenates a 0 to the head of the string:
$$
Ts=0s.
$$
We can now state and prove the second functional equation.
\begin{lemma} For every $s\in \{0,1\}^\omega$, one has
$\Jimm_{\partial\F}(Ts)=TU\Jimm_{\partial\F}(s)$.
Hence 
\begin{equation}\label{jimmt}
\Jimm_\R(1+x)=1+\frac{1}{\Jimm_\R(x)}.
\end{equation}
\end{lemma}
{\it Proof.} Using the fact that $\phi=0\lnot \phi$ and that 
$p\xor (\lnot q)=\lnot(p\xor q)$ we get
$$
\Jimm_{\partial\F}(Ts)=\Jimm_\F(0s)=\phi \xor 0s=(0\xor 0)((\lnot \phi)\xor s) 
=0 (\lnot (\phi \xor s))=\qquad\qquad\qquad\qquad
$$
$$
\qquad\qquad\qquad\qquad\qquad\qquad\qquad\qquad\qquad\qquad\qquad=0U\Jimm_\F(s)=TU\Jimm_{\partial\F}(s). 
\qquad \Box
$$}

\medskip\noindent
{\bf Jimm as a limit of piecewise-$\pgl$ functions.}
In $\F$, denote by  
$V_\bullet(n)$ the set of vertices of 
distance $\leq n$ to the base (excluding $v^*$), and let $V_\bullet^\prime(n)$ be the set of vertices of odd distance $\leq n$ to the base. Then one has
$$
\Jimm_\F=
\lim_{n\to\infty} \theta_{V_\bullet(n)}
=\lim_{n\to\infty} \sigma_{V_\bullet^\prime(n)}.
$$
The maps $\theta_{V_\bullet(n)}$ and $\sigma_{V_\bullet^\prime(n)}$ 
induce a sort of finitary 
projective interval exchange maps on $\widehat{\R}$
and thus $\Jimm_\R$ can be written as a limit of such functions. 
Below we draw the twists $\theta_{V_\bullet(n)}$ for $n=1\dots 9$. 
\sherh{These are obtained by xoring with the sequences $(01)^n1^\omega$ or $(01)^n0^\omega$. The non-defined values 0 and $\infty$ can be thought of those rationals having nothing to be xored with.}
\begin{center}\label{jimmiterate}
\noindent
{\includegraphics[scale=.2]{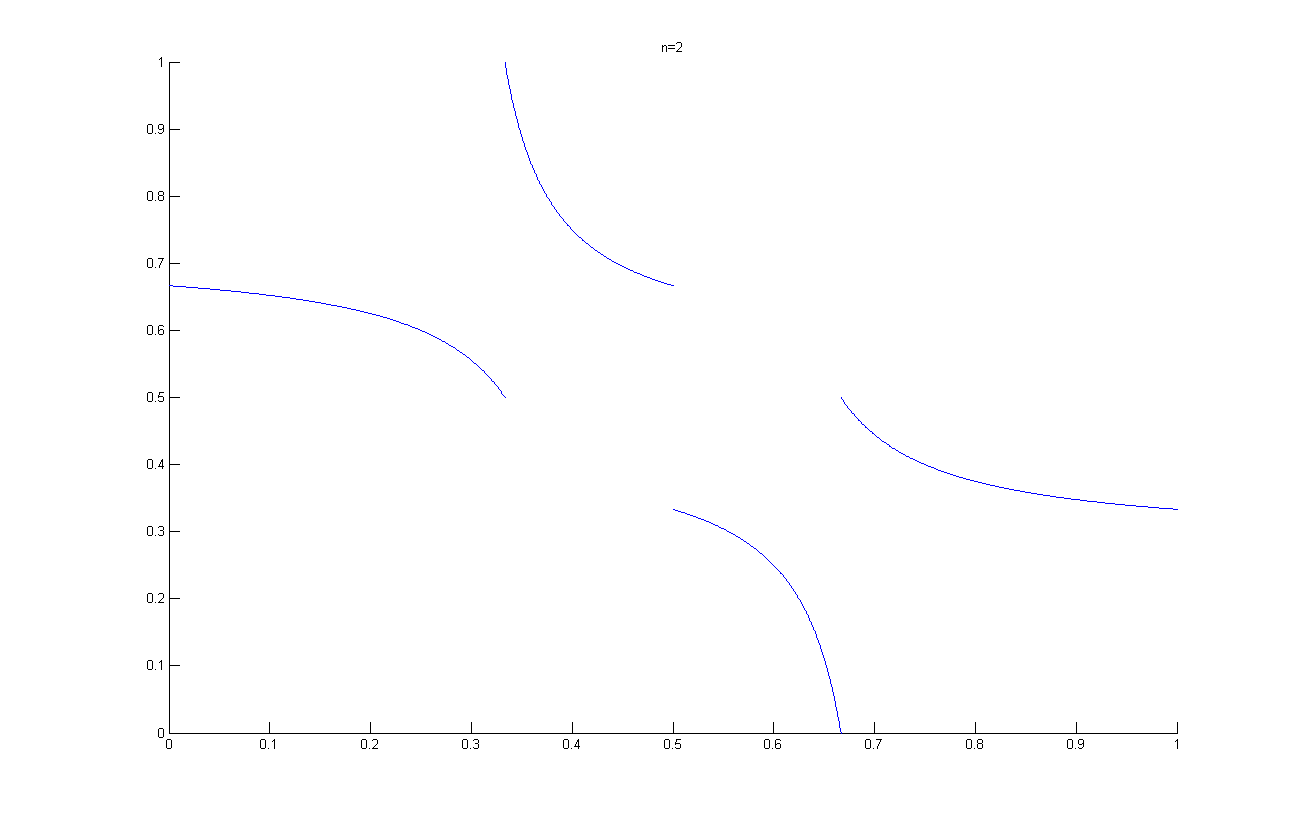}}
{\includegraphics[scale=.2]{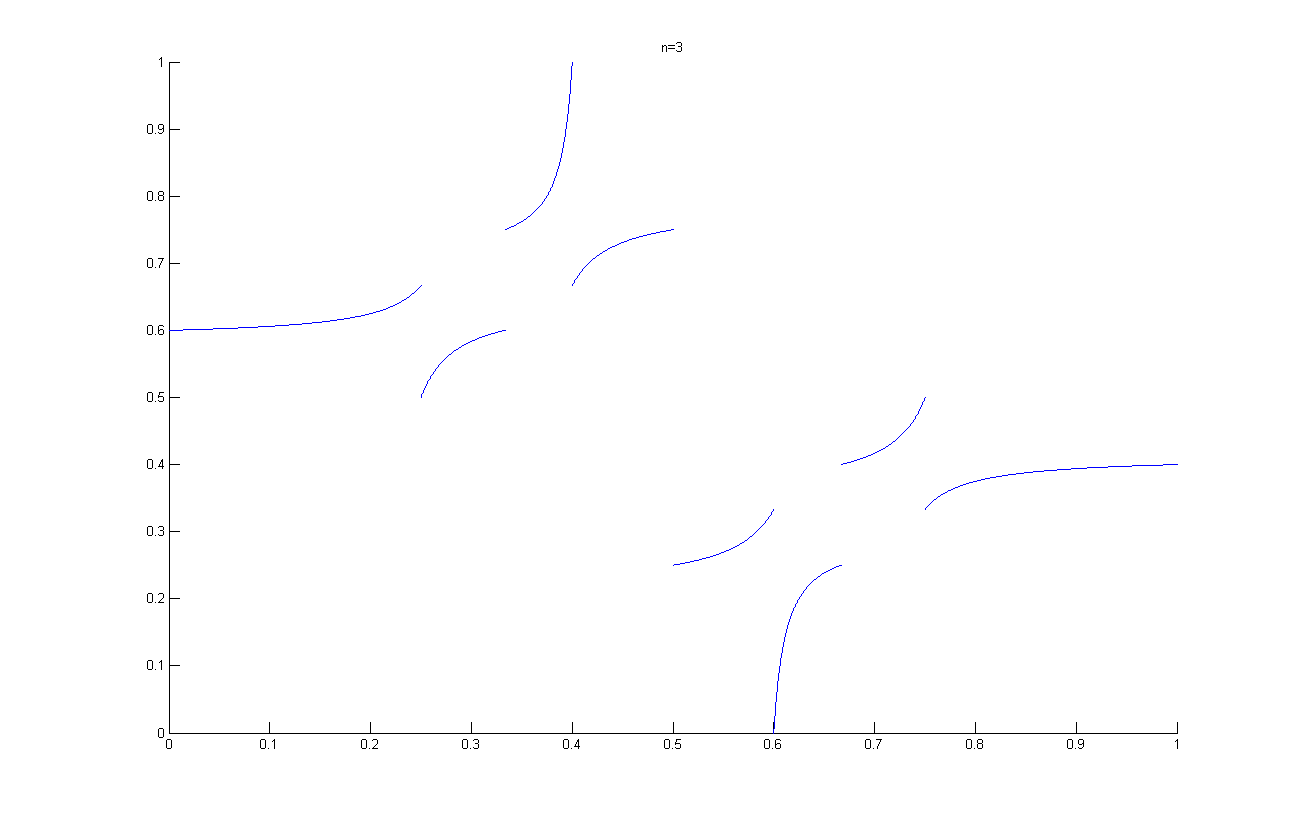}}
{\includegraphics[scale=.2]{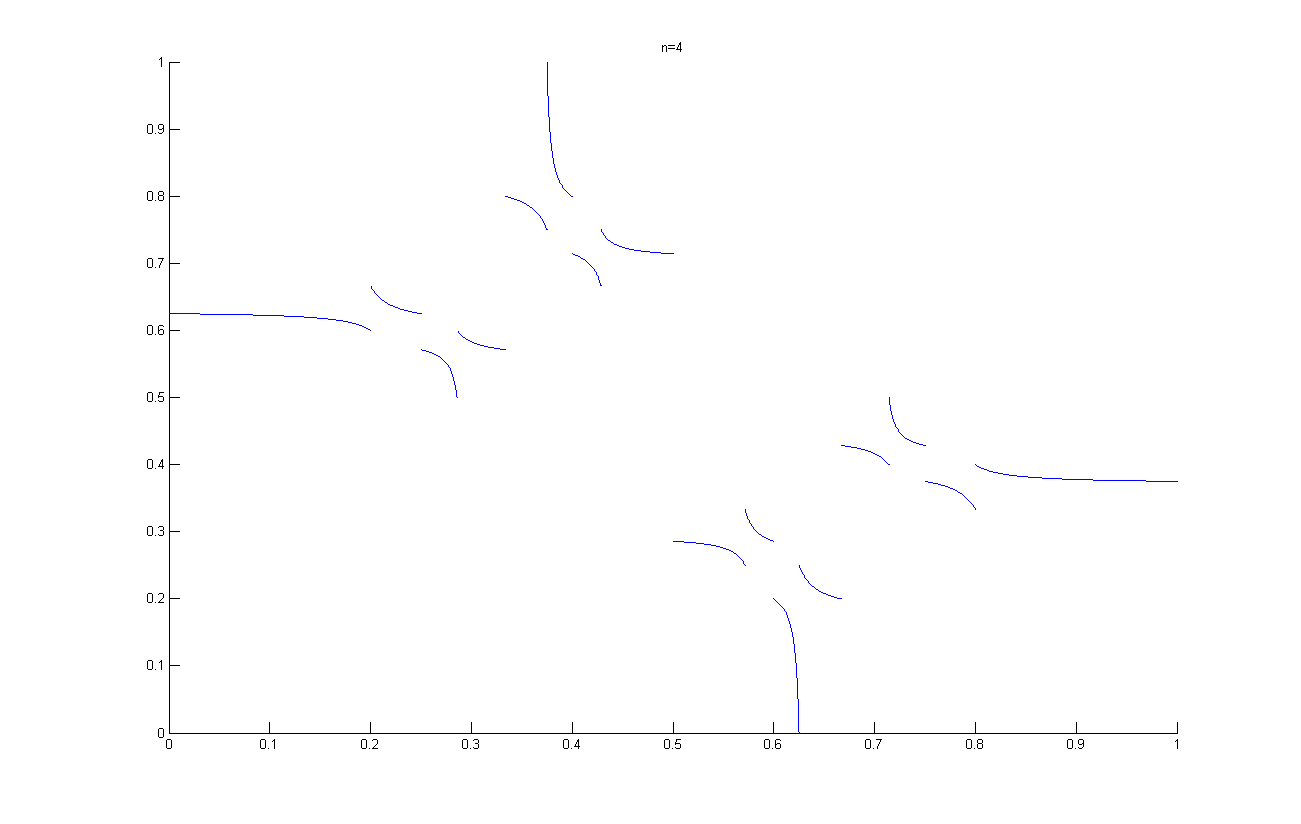}}
{\includegraphics[scale=.2]{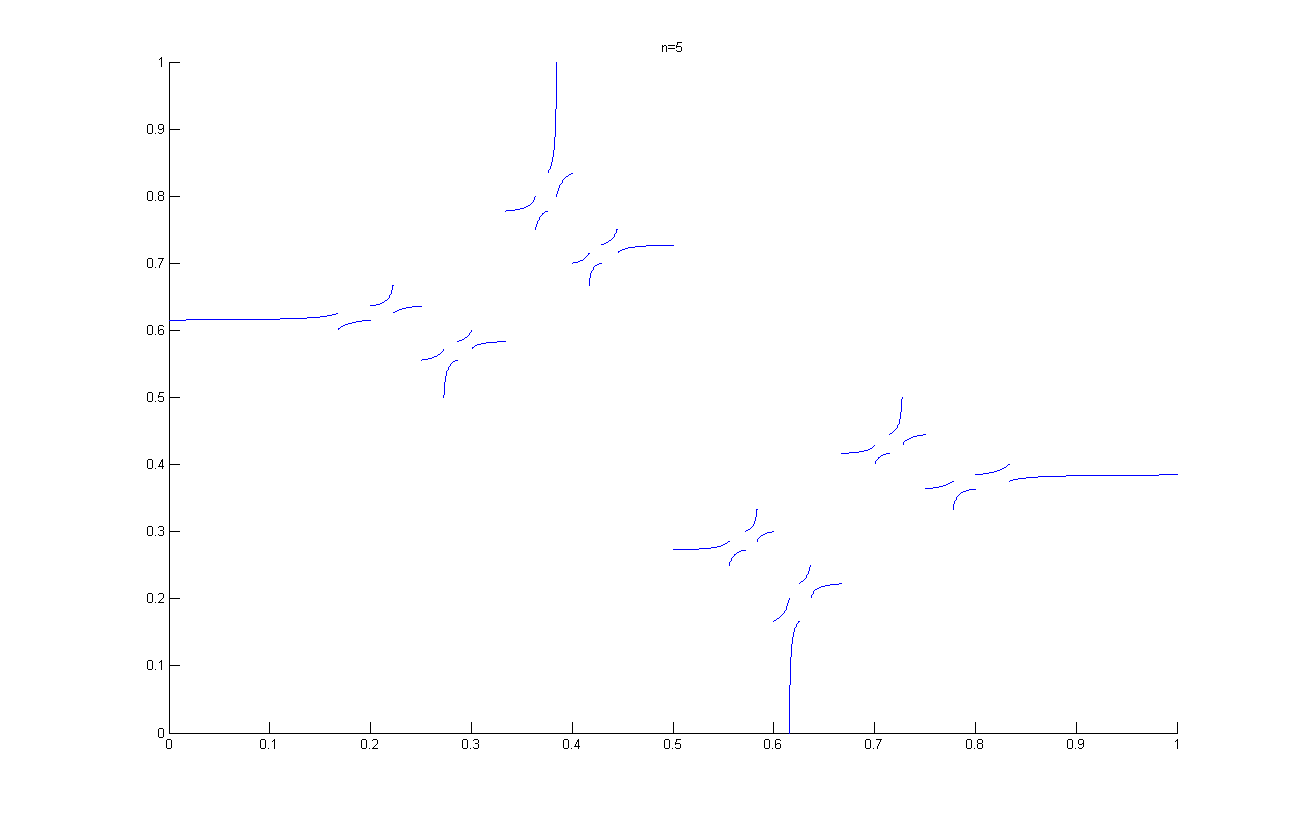}}\\
{\includegraphics[scale=.2]{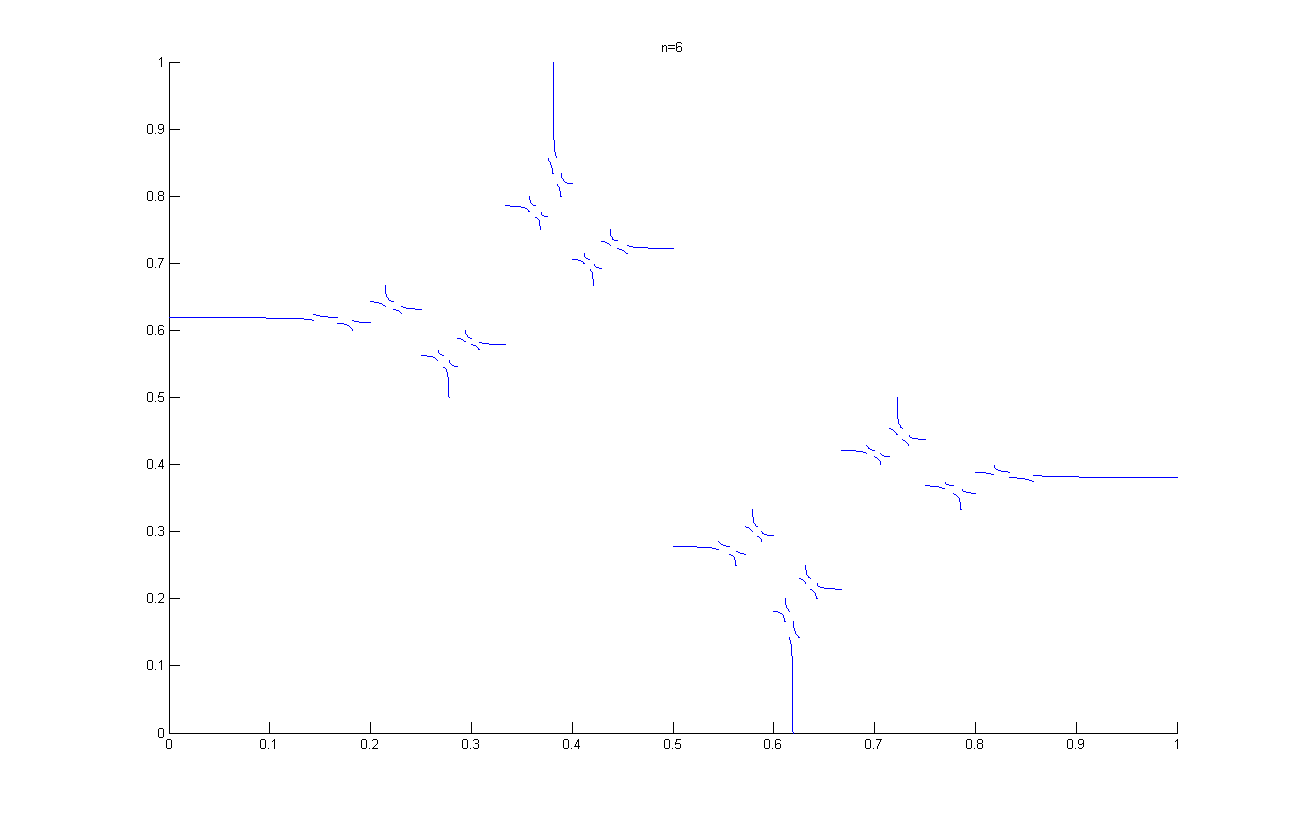}}
{\includegraphics[scale=.2]{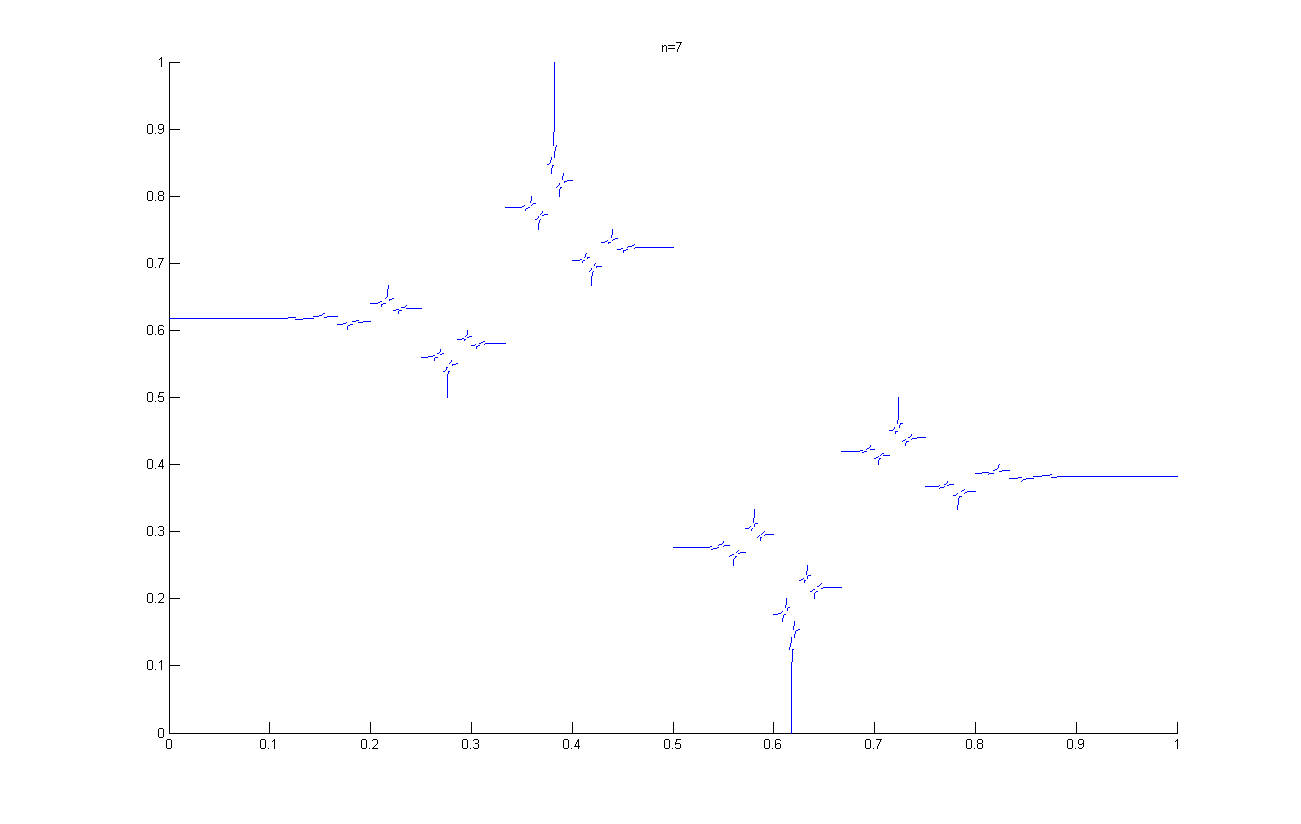}}
{\includegraphics[scale=.2]{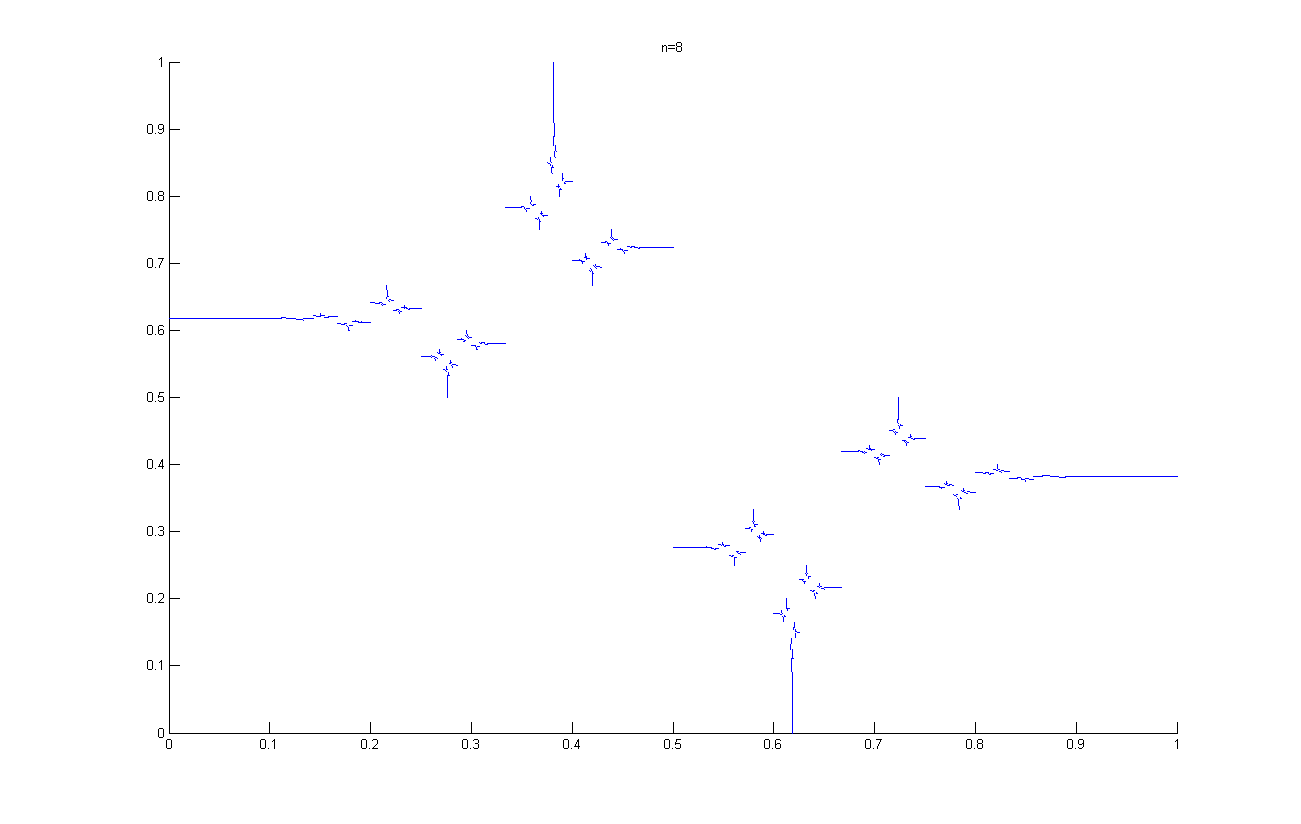}}
{\includegraphics[scale=.2]{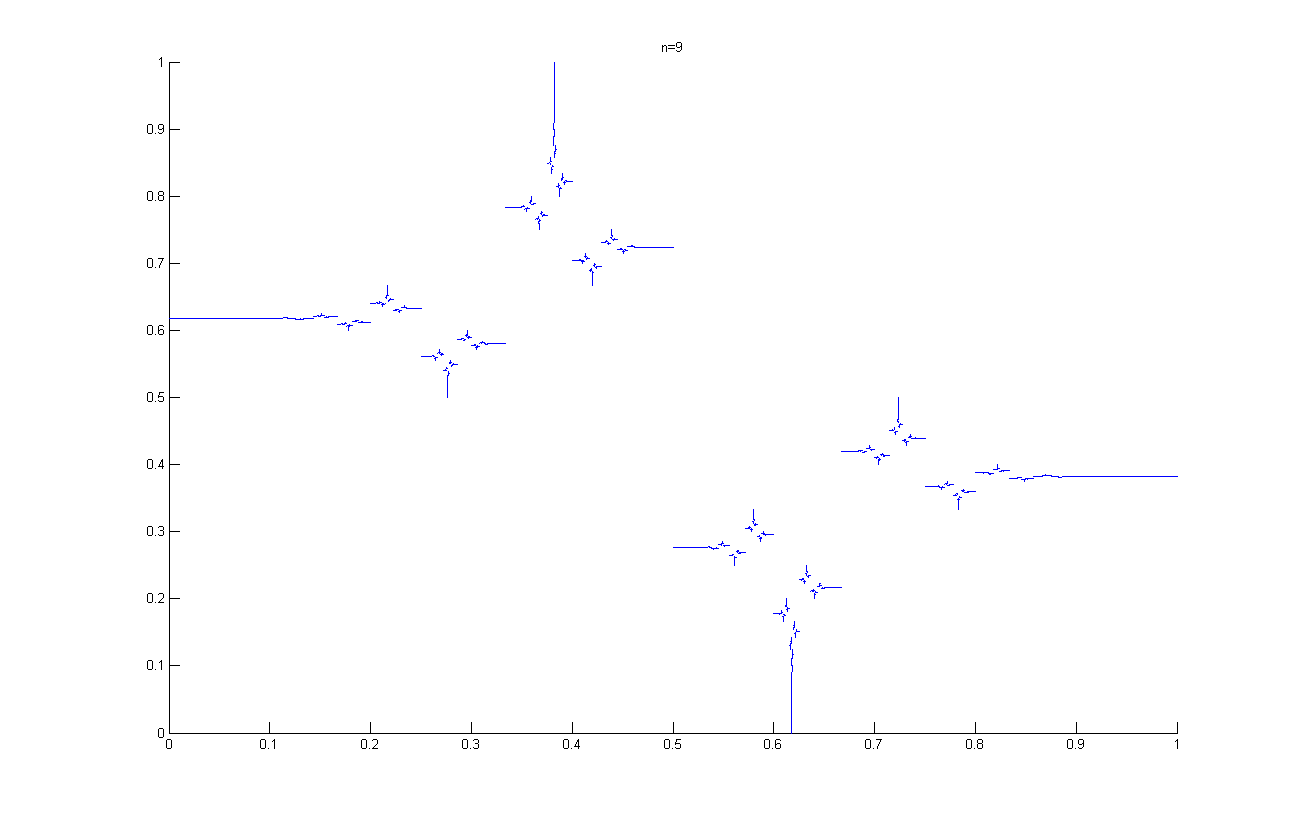}}
{\small {\bf Figure.} The plot of $\theta_{\F_n^*}$ on the interval $[0,1]$, for $n=2,3,...9$.}\\
\end{center}

\unut{It is of interest to further study the effects of automorphisms of $|\F|$ on the boundary circle and see how they conjugate the dynamical system of the Gauss map.
For example, any subgroup $\Gamma\subset \psl$ (or conjugacy class of a subgroup) corresponds to a set of vertices defining 
two automorphisms $\theta_\Gamma$ and $\sigma_\Gamma$. Any element $\gamma$ of the boundary $\partial\F$ also defines a shuffle-automorphism $\sigma_\gamma$  and a  twist-automorphism $\theta_\gamma$, constructed by using the set of vertices lying on the path connecting the base edge to $\gamma$. }

\unut{\paragraph{Modular graphs.} Since the modular group acts on the Farey tree, so does any subgroup $\Gamma$ of the modular group, and the resulting quotient graph $\F/\Gamma$ of the Farey tree $\F$ is called a {\it modular graph}, see \cite{kunming}. It is a bipartite graph with a ribbon structure. Conjugate subgroups give rise to isomorphic modular graphs. 
These are in fact special\nt{\footnote{Dessins describe the coverings of the thrice punctured sphere whereas modular graphs describe the coverings of the modular curve. Since the thrice-punctured sphere is a degree-6 covering of the modular curve (corresponding to the congruence subgroup of level 2), modular graphs are no weaker then dessins in their descriptive capacity.}}  kind of dessins, the ambient punctured surface of a modular graph carries a canonical arithmetic  structure. 

The modular group is not invariant under $\Jimm$, however, the subgroup 
$\langle R, SRS\rangle\simeq \Z/3\Z*\Z/3\Z$ is invariant and its subgroups are sent to each other under $\Jimm$. 
Hence $\Jimm$ acts on the corresponding modular graphs by 
$\F/\Gamma\to \F/\Jimm\Gamma$. This action shuffles (i.e. reverses the cyclic ordering of) every other trivalent vertex of the graph.
The genus of the ambient surface is not invariant under this action. On the other hand, its fundamental group is invariant. 
Jones and Thornton \cite{jones1983operations} described this action in terms of the dual graphs (i.e. triangulations).
The use of triangulations requires to cut the ambient surface into pieces, shuffle, and then glue the pieces back; which is somewhat hard to imagine for us. }

\sherh{Note that twists of modular graphs are not well-defined.}

\sherh{The class of finite modular trees with terminal points of type $\bullet$ only, is invariant under the $\Jimm$-action.
This action admits a simple description in this case (i.e. shuffling every other vertex). It might be of interest to compare the generalized Chebyshev polynomials to see if they are related. This is the only concrete connection with the Galois theory I could fathom.}

\unut{On the other hand, note that metrization of the graphs $\F/\Gamma$ parametrize the decorated Teichm\"uller space of Penner 
\cite{penner2012decorated} and we see that $\Jimm$ induces a certain duality on punctured Riemann surfaces. This duality do not respect the genus, but it does respect the rank of the fundamental group.}

\paragraph{Jimm on $\Q\setminus\{0,1\}$.}\label{jimmq}
In virtue of Lemma \ref{correspondences}, there is a 1-1 correspondence between the set of trivalent vertices and the set 
$\Q\setminus\{0,1\}$. Since any element of $Aut_I(\F)$ defines a bijection of $V_\bullet(\F)$, we see that every automorphism of $\F$ that fixes $I$, defines a unique bijection of $\Q\setminus\{0,1\}$. In particular, this is the case with $\Jimm_\F$. We denote the resulting
bijection with $\Jimm_\Q$. One has $\Jimm_\Q(1)=1$, and for $x>0$, its values can be computed by using the functional equations 
$\Jimm_\Q(1+x)=1+1/\Jimm_\Q(x)$ and $\Jimm_\Q(1/x)=1/\Jimm_\Q(x)$. It tends to $\Jimm_\R$ at irrational points. However, it must be emphasised that $\Jimm_\Q(x)$ is not the restriction of $\Jimm_\R$ to $\Q$; this latter function is by definition two-valued at rationals\footnote{It might have been convenient to declare the values of $\Jimm_\R$ at rational arguments to be given by $\Jimm_\Q(x)$, so that $\Jimm_\R$ would be a well-defined function everywhere (save 0 and $\infty$). However, we have chosen to not to follow this idea, for the sake of uniformity in definitions.}.
The involution $\Jimm_\Q$ conjugates the multiplication on $\widehat\Q\setminus\{0,\infty\}$ to an operation with 1 as its identity, and such that the inverse of $q$ is $1/q$. 

\subsection{Fun with functional equations}
It is possible to derive many equations from the functional equations for $\Jimm$, or express them in alternative forms.
In this section we record some of these equations. We leave the task of verifying to the reader. 
(As usual we drop the subscript $\R$ from $\Jimm_\R$ to increase the readability):

To start with, note that $\Jimm$ satisfies the following functional equations
$$
\Jimm\left(1-\frac{1}{x}\right)=1-\frac{1}{\Jimm(x)}, \quad 
\Jimm\left(\frac{x}{x+1}\right)+\Jimm\left(\frac{1}{x+1}\right)=1.
$$

There is the two-variable version of the functional equations
\smallskip\begin{center}\label{twovariable}
\fbox{\fbox{\begin{minipage}{10cm}
$$
\Jimm(x)=y  \iff \Jimm(y)=x
$$

$$
xy=1\iff \Jimm(x)\Jimm(y)=1
$$

$$
x+y=0\iff \Jimm(x)\Jimm(y)=-1
$$

$$
x+y=1\iff \Jimm(x)+\Jimm(y)=1
$$

$$
\frac{1}{x}+\frac{1}{y}=1\iff \frac{1}{\Jimm(x)}+\frac{1}{\Jimm(y)}=1
$$

\vspace{1mm}
\end{minipage}}}
\end{center}

\unut{\noindent
Functional equations can be expressed in terms of the involution
$\Jimm^*(x):=1/\Jimm(x)=U\Jimm(x)$: 

\noindent
\begin{center}
\fbox{\fbox{\begin{minipage}{10cm}

$$
\Jimm^*(\Jimm^*(x))=x 
$$

$$
xy=1\iff \Jimm^*(x)\Jimm^*(y)=1
$$

$$
x+y=0\iff \Jimm^*(x)\Jimm^*(y)=-1
$$

$$
x+y=1\iff \frac{1}{\Jimm^*(x)}+\frac{1}{\Jimm^*(y)}=1
$$

$$
\frac{1}{x}+\frac{1}{y}=1\iff {\Jimm^*(x)}+{\Jimm^*(y)}=1
$$

$$
\Jimm^*(x+1) ={1 \over 1+\Jimm^{*}(x)}
$$
\end{minipage}}}
\end{center}}
\noindent

\sherh{
Another way to write these equations is: 
$$
\Jimm^2=U^2=K^2=V^2=Id, \quad [\Jimm,U]=0, \quad [\Jimm,K]=0, \quad [\Jimm,V]=\Jimm-\frac{1}{\Jimm}
$$}

\sherh{Is there a way to modify/deform these functional equations to obtain new, interesting ones? 
One way to deform the equations might be as follows.
$$
\Jimm(1+x)=1+\alpha/x, \quad \Jimm(\beta/x)=\gamma/\Jimm(x)
$$
One should of check if these satisfy the group relations for $\pgl$. 
}
Finally, the functional equations can be expressed in terms of the function $\Psi(x):=\log |\Jimm(x)|$, where they look like cocycle relations, compare~\cite{zagier2000quelques}.

\noindent
\begin{center}
\fbox{\fbox{\begin{minipage}{10cm}

$$
\Psi(1/x)+\Psi(x)=0
$$

$$
\Psi(x)+\Psi(1+x)+\Psi(x/(1+x))=0
$$

$$
\exp(\Psi(1-x))+\exp({\Psi(x)})=1
$$

\end{minipage}}}
\end{center}
\noindent
(the involutive relation is left to the reader).

\unut{One is tempted to seek the solutions of (some of the) functional equations which are analytic in some region, the most appealing one being the equation $\Psi(1+x)=1+1/\Psi(x)$. Beware that by iterating this equation we get
\begin{equation}\label{limit}
\lim_{x\to +\infty} \Psi(x)=\Phi \mbox{ and } \lim_{x\to -\infty} \Psi(x)=-1/\Phi,
\end{equation}
which implies that there are no rational solutions. Here $\phi=(1+\sqrt{5})/2$ as usual.
We don't know if there exists an analytic solution in some strip around the real line. Another temptation is to extend 
$\Jimm$ to the upper half plane in some sense, but we have no idea how this must be done.

Since $\Jimm_\R$ satisfies this functional equation, its limit at $\pm \infty$ exists and these limits are given by 
(\ref{limit}). See also Lemma \ref{easyconsequence} (ii) below.}

\section{Analytical aspects of jimm}
\subsection{Fibonacci sequence and the Golden section}
The fabulous Fibonacci sequence is defined by the recurrence 
$$
F_0=0, \quad F_1=1,\mbox{ and } F_{n}=F_{n-1}+F_{n-2} \mbox{ for } n\in \Z,
$$
one has then $F_{-n}=(-1)^{n+1}F_n$. Here is a table of some small Fibonacci numbers:
$$
\begin{array}{|c|c|c|c|c|c|c|c|c|c|c|c|c|c|c|c|}
\hline
\dots& F_6&F_{-5}&F_{-4}&F_{-3}&F_{-2}&F_{-1}&F_0&F_1&F_2&F_3&F_4&F_5 &F_6&F_7\dots \\ \hline
\dots&-8&5&-3&2&-1&1&0&1&1&2&3&5&8&13\dots\\ \hline
\end{array}
$$
One has 
$$
\widetilde{T}(x)=1+\frac{1}{x}\implies
\widetilde{T}^n=\frac{F_{n+1}x +F_n}{F_nx+F_{n-1}} \quad (n\in \Z)
$$
and therefore
$$
\det(\widetilde{T}(x))=-1
\implies
\det(\widetilde{T}^n)=(-1)^n 
\implies
F_{n+1}F_{n-1}-F_n^2=(-1)^n \quad (n\in \Z)
$$
The transformation $\widetilde{T}$ has two fixed points, solutions of the equation
\begin{equation}\label{fixedpoint}
\widetilde{T}(x)=x \iff 1+\frac{1}{x}=x \iff  x^2-x-1=0 \iff  x=\frac{1\pm \sqrt{5}}{2}
\end{equation}
As we already mentioned, the positive root of this equation is called the golden section which we shall denote by $\Phi$. 
It is the limit
$
\Phi=\lim_{n\rightarrow \infty} {F_{n+1}}/{F_n}.
$
The negative root of Equation (\ref{fixedpoint}) equals $\Phi_*:=-1/\Phi$.
The numbers $F_n$ grow exponentially as $F_n=(\Phi^n-\Phi_*^n)/\sqrt{5}$ by Binet's formula. 
This fact will be used in the following form
\begin{equation}\label{binet}
F_n>{1}{\sqrt{5}}\Phi^n \implies F_n^{-1}<{\sqrt{5}}\Phi^{-n} 
\end{equation}

\sherhh{Summing Binet formulas with weight $z^n$ we obtain the generating function
$$
\frac{e^{\phi z}-e^{z\overline{\phi}}}{\phi-\overline{\phi}}=\sum_{n=0}^\infty F_n{z^n \over n!}
$$
Another beautiful generating function is
$$
\sum_{n=0}^\infty F_nz^n=\frac{1}{1-(z+z^2)}
$$
We also record the useful and well-known formula
$$
F_n=F_{k+1}F_{n-k}+F_kF_{n-k-1} \quad (n,k\in \Z)
$$}

The following lemma is an easy consequence of the preceding discussions.
\begin{lemma}\label{easyconsequence}
Let $x$ be an irrational number or an indeterminate. Then\\
(i) 
$$
\Jimm(1+x)=1+\frac{1}{\Jimm(x)} \iff \Jimm(Tx)=\widetilde{T}\Jimm(x).
$$
(ii) 
$$
\Jimm(n+x)=\widetilde{T}^n\Jimm(x)=\frac{F_{n+1}\Jimm(x) +F_n}{F_n \Jimm(x)+F_{n-1}} \quad (n\in \Z).
$$
(iii) 
$$
\Jimm([1,x])=\Jimm(1+\frac{1}{x})=1+\Jimm(x) \implies \Jimm([1_n,x])=n+x \quad (n\in \Z).
$$
(iv)
$$
\Jimm([n,x])=\Jimm(n+1/x)=\frac{F_{n+1}+F_n\Jimm(x) }{F_n +F_{n-1}\Jimm(x)} \quad (n\in \Z).
$$
\end{lemma}
\begin{lemma} One has $0\leq\Jimm(x)\leq 1$ if $0\leq x\leq 1$,
with $\Jimm(1/\Phi)=0$, $\Jimm(1/(1+\Phi))=1$, 
$\Jimm(1/(n+\Phi))=F_n/F_{n+1}$. Hence $\Jimm$ restricts to an involution of the unit interval.
More generally, $\Jimm$ maps Farey intervals to Farey intervals.
\end{lemma}
 \unut{The lemma below is also a routine observation.
\begin{lemma} If $x$ is given by the ``minus-continued fraction" 
$$
x=n_0-\cfrac{1}{n_1-\cfrac{1}{n_2-\cfrac{1}{\ldots}}}
=T^{n_0}\circ UV  \circ T^{n_1}\cdot UV  \dots,
$$
then $\Jimm(x)$ is given by the continued fraction
$$
\Jimm(x)=\widetilde{T}^{n_0}\circ V \circ \widetilde{T}^{n_1}\circ V \circ\widetilde{T}^{n_2}\circ \dots. $$
\end{lemma}
Finally, one has the following simple observation.
\begin{lemma}
The following is a $\pgl$-invariant of four-tuples of real numbers
$$
\bigl[ x:y:z:t \bigr]_\jimm:=
\bigl[ \Jimm(x): \Jimm(y): \Jimm(z): \Jimm(t) \bigr],
$$
where $\bigl[ x:y:z:t \bigr]$ denotes the cross ratio.
\end{lemma}}
\subsection{Continuity and jumps}
Since the involution $\Jimm_{\partial\F}: \partial\F\to \partial\F$
is a homemorphism and since the subspaces 
$$
\R\setminus \Q \simeq \partial\F\setminus \{\mbox{rational paths}\}
$$ 
are also homeomorphic, the function  $\Jimm_\R$ is continuous at irrational points. Let us record this fact.
\begin{theorem}\label{jimmcontini}
The function $\Jimm_\R$ on $\R\setminus \Q$ is continuous. 
\end{theorem}
Moreover, recall that there is a canonical ordering on $\partial\F$
inducing on $\widehat{\R}$ its canonical ordering compatible with its topology.
So the notion of lower and upper limits exists on $\partial\F$
and coincides with lower and upper limits on $\widehat{\R}$.
This shows that the two values that $\Jimm_\R$ assumes on rational arguments are nothing but the limits
$$
\Jimm(q)^-:=\lim_{x\to q^-} \Jimm(x) \mbox{ and }
\Jimm(q)^+:=\lim_{x\to q^+} \Jimm(x) 
$$
By choosing one of these values coherently, one can make $\Jimm$ an everywhere upper (or lower) continuous function. Note however that it will (partially) cease to satisfy the functional equations at rational arguments. One has, for irrational $r$ and rational $q$,
$$
\lim_{q\rightarrow r} \Jimm_\R(q)= \Jimm_\R(r),
$$
no matter how we choose the values of $\Jimm_\R(q)$.

The jump function 
$$
\delta_\jimm(q)=\Jimm(q)^+-\Jimm(q)^-
$$ 
can be considered as a (signed) measure of complexity of a rational number $q$. For irrational numbers it is 0. 
It splits the set $\Q$ into two pieces, according to the sign of $\delta_\jimm(q)$.

The following functional equations are easily deduced from those for $\Jimm$:

\smallskip\begin{center}
\fbox{\fbox{\begin{minipage}{10cm}
$$
\delta_\jimm\left(-1/x\right)=-\delta_\jimm(x) \iff \delta_\jimm\left(-x\right)=-\delta_\jimm\left(1/x\right)
$$

$$
\delta_\jimm\left(1-x\right)=\delta_\jimm(x) \iff \delta_\jimm\left(1+x\right)=\delta_\jimm(-x)
$$

$$
\delta_\jimm\left(-x\right)=-{\delta_\jimm(x) \over \Jimm^+(x)\Jimm^-(x)}
$$

\vspace{1mm}
\end{minipage}}}
\end{center}
Hence $|\delta_\jimm|$ is invariant under the infinite dihedral group $\langle -1/x, 1-x \rangle$. 
In order to compute some values of $\delta_\jimm(q)$ more explicitly, one has, for $k$ odd,
$$
q=[n_0, n_1,\dots n_k,\infty]=[n_0,n_1,\dots n_k-1,1,\infty]\implies
$$
$$
\Jimm(q)^+=\Jimm([n_0, n_1,\dots n_k,\infty]), \quad 
\Jimm(q)^-=\Jimm([n_0,n_1,\dots n_k-1,1,\infty])
$$
$$
\implies\delta_\jimm(q)=
\Jimm([n_0, n_1,\dots n_k,\infty])-
\Jimm([n_0,n_1,\dots n_k-1,1,\infty])
$$
(and for $k$ even, $\delta_\jimm(q)$ changes sign).

\begin{center}\label{jimmplot}
\noindent{\includegraphics[scale=.4]{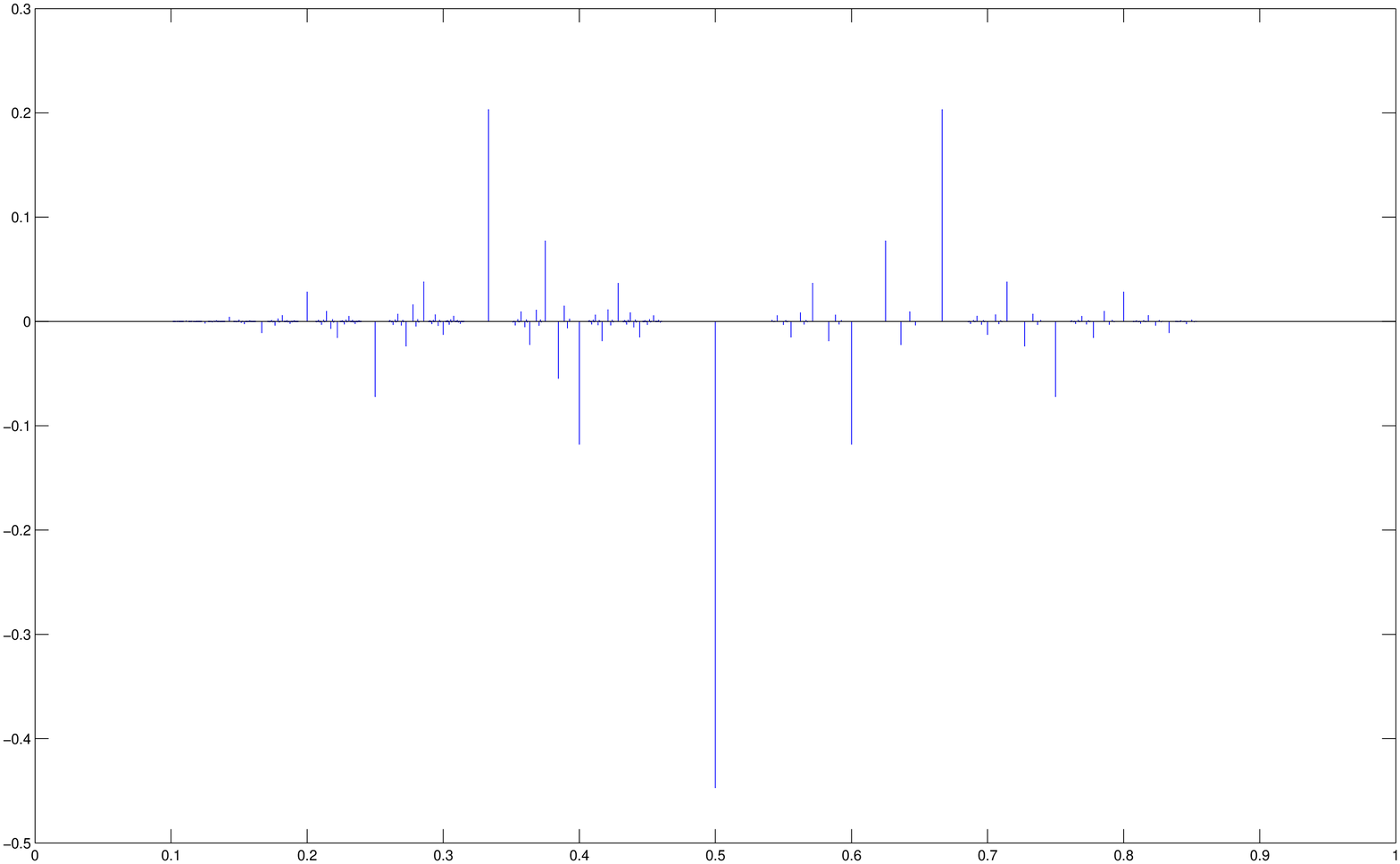}}\\
{\small {\bf Figure.} The plot of $\delta_\jimm$, centered around the point $x=1/2$. }\\
\end{center}

Suppose $q=n$ is an integer. Then
$$
\Jimm(n)^+=\frac{F_n\Phi+ F_{n+1}}{F_{n-1}\Phi+F_n}
$$
$$
\Jimm(n)^-=\frac{F_{n-1}\Phi+ F_{n}}{F_{n-2}\Phi+F_{n-1}}
$$
More conceptually, one has
$$
\Jimm[n,x]=\left[\begin{array}{ll} 
F_{k+1}& F_k\\ F_k& F_{k-1}\end{array}\right]
\Jimm[n-k,x] \mbox{ for } 0\leq k \leq n.
$$
Hence, 
$$
\Jimm[n,\infty]=\left[\begin{array}{ll} 
F_{n}& F_{n-1}\\ F_{n-1}& F_{n-2}\end{array}\right]
\Jimm[1,\infty], \quad
 \Jimm[n-1,1,\infty]=\left[\begin{array}{ll} 
F_{n}& F_{n-1}\\ F_{n-1}& F_{n-2}\end{array}\right]
\Jimm[0,1,\infty].
$$
Since 
$$
\Jimm[1,\infty]=1+1/\Jimm[0,\infty]=1+\Jimm[\infty]=1+\Phi=\Phi^2
$$
and 
$$
\Jimm[0,1,\infty]=1/\Jimm[1,\infty]=1/(1+\Phi)=\Phi^{-2}
$$
we 
$$
\delta_\jimm(n)=
\left[\begin{array}{ll} 
F_{n}& F_{n-1}\\ F_{n-1}& F_{n-2}\end{array}\right]
\Phi^{2}
 -
\left[\begin{array}{ll} 
F_{n}& F_{n-1}\\ F_{n-1}& F_{n-2}\end{array}\right]
\Phi^{-2}
$$
$$
=(-1)^{n+1}\frac{\Phi^{2}-\Phi^{-2}}{(F_{n-1}\Phi^{2}+F_{n-2})(F_{n-1}\Phi^{-2}+F_{n-2})}
$$
and we finally have
$$
\delta_\jimm(n)=
(-1)^{n+1}\frac{\sqrt{5}}{(F_{n-1}^2+3F_{n-1}F_{n-2}+F_{n-2}^2)}=
(-1)^{n+1}\frac{\sqrt{5}}{F_{n}^2+F_{n-1}F_{n-2}}
$$
We see that this rapidly tends to zero as $n\rightarrow\infty$.

\unut{Now, it should be true that $\Jimm$ is of bounded variation on any finite or infinite closed interval 
not containing $\Phi$ nor $-\Phi^{-1}$, and that the signed measure $d\Jimm$ is just the function $\delta_\jimm$, 
but we defer the proofs out of weariness. On the other hand, it is proved in Theorem~\ref{derivativevanish}
the derivative of $\Jimm$ exists vanishes almost everywhere. }

\sherh{
{\bf Variation of $\Jimm$}
One may also consider  the function $\int |d\Jimm|$.

$\Jimm$ defines a signed measure and it must be studied here.

(from wikipedia)
It can be proved that a real function f is of bounded variation in an interval if and only if it can be written as the difference $f = f_1-f_2$ of two non-decreasing functions: this result is known as the Jordan decomposition.
Through the Stieltjes integral, any function of bounded variation on a closed interval [a, b] defines a bounded linear functional on C([a, b]). In this special case,[3] the Riesz representation theorem states that every bounded linear functional arises uniquely in this way. The normalized positive functionals or probability measures correspond to positive non-decreasing lower semicontinuous functions.

Functions of bounded variation have been studied in connection with the set of discontinuities of functions and differentiability of real functions, and the following results are well-known. If f is a real function of bounded variation on an interval [a,b] then

f is continuous except at most on a countable set;

have only jump-type discontinuities 

f has one-sided limits everywhere (limits from the left everywhere in (a,b], and from the right everywhere in [a,b) ;

the derivative f'(x) exists almost everywhere (i.e. except for a set of measure zero).}

\unut{
Concerning the integrals of $\Jimm$, below is an easy consequence of the functional equation $\Jimm K=K\Jimm$.
\begin{lemma}\label{integral}
One has
$$
\int_0^1\Jimm(x)dx=\frac{1}{2}, \quad \int_0^t\Jimm(x)dx= t - 1/2 + \int_0^{1-t} \Jimm(x) dx
$$
\end{lemma}
By the change of variable $y=1/x$ this gives
$$
\int_1^\infty \frac{dy}{y^2\Jimm(y)}=\frac{1}{2}
$$
We failed to compute other integrals (moments, transforms, etc). 
On the other hand, as we remarked in the introduction we have developed some numerical algorithms and implemented them on a computer for the interested reader.}

\unut{\subsection{The graph of jimm}
Since $\Jimm_\R$ has a non-removable discontinuity at every rational point, it is not possible to draw its graph in the usual way. Below is a computer-produced box-graph of $\Jimm(x)$ for $x\geq 0$, the graph lies inside the smaller (and darker) boxes. 

This graph is found as follows.
If $x\in (1,2)$, then $x=[1,n_1,r]$ with $r\in (1,\infty)$, and there are two cases
$$x\in (1,\infty)=\bigcup_{n_0=1}^{\infty} (n_0,\infty)\iff x=[n_0,r] \mbox{ with } n_0\geq 1, \mbox{ and } r\in (1,\infty)$$
\begin{eqnarray*}
(n_0\geq 1), \quad x\in (n_0,n_0+1) =\bigcup_{n_1=1}^{\infty}\left[n_0+\frac{1}{n_1+1}, n_0+\frac{1}{n_1}\right]\\
x=[n_0,n_1, r] \mbox{ with }  n_1\geq 1, \quad 1<r<\infty \\
\implies \Jimm(x)=[\underbrace{1,1,\dots,1}_{n_0-1\, times},m_1,\underbrace{1,1,\dots,1}_{n_1-2 \, times},\dots,1,2,\dots ],
\end{eqnarray*}
where $m_1$=2 if $n_1>1$ and $m_1\geq 3$ if $n_1=1$. So in any case, $m_1\geq 2$, so that 
$$
(n_0\geq 1), \quad x\in (n_0,n_0+1) \implies \Jimm(x)\in [\underbrace{1,1,\dots,1}_{n_0-1\, times}, [2,\infty]]
$$
In other words,
$$
(n_0\geq 2), \implies \quad \Jimm((n_0,n_0+1))=\left[\frac{F_{n_0+3}}{F_{n_0+2}}, \frac{F_{n_0+1}}{F_{n_0}} \right]
$$
(intervals being inverted depending on the parity of $n_0$).
Now suppose that $x\in [2,3)$. Continuing in this manner gives the following graph:

\begin{center}\label{jimmplot}
\noindent{\includegraphics[scale=.5]{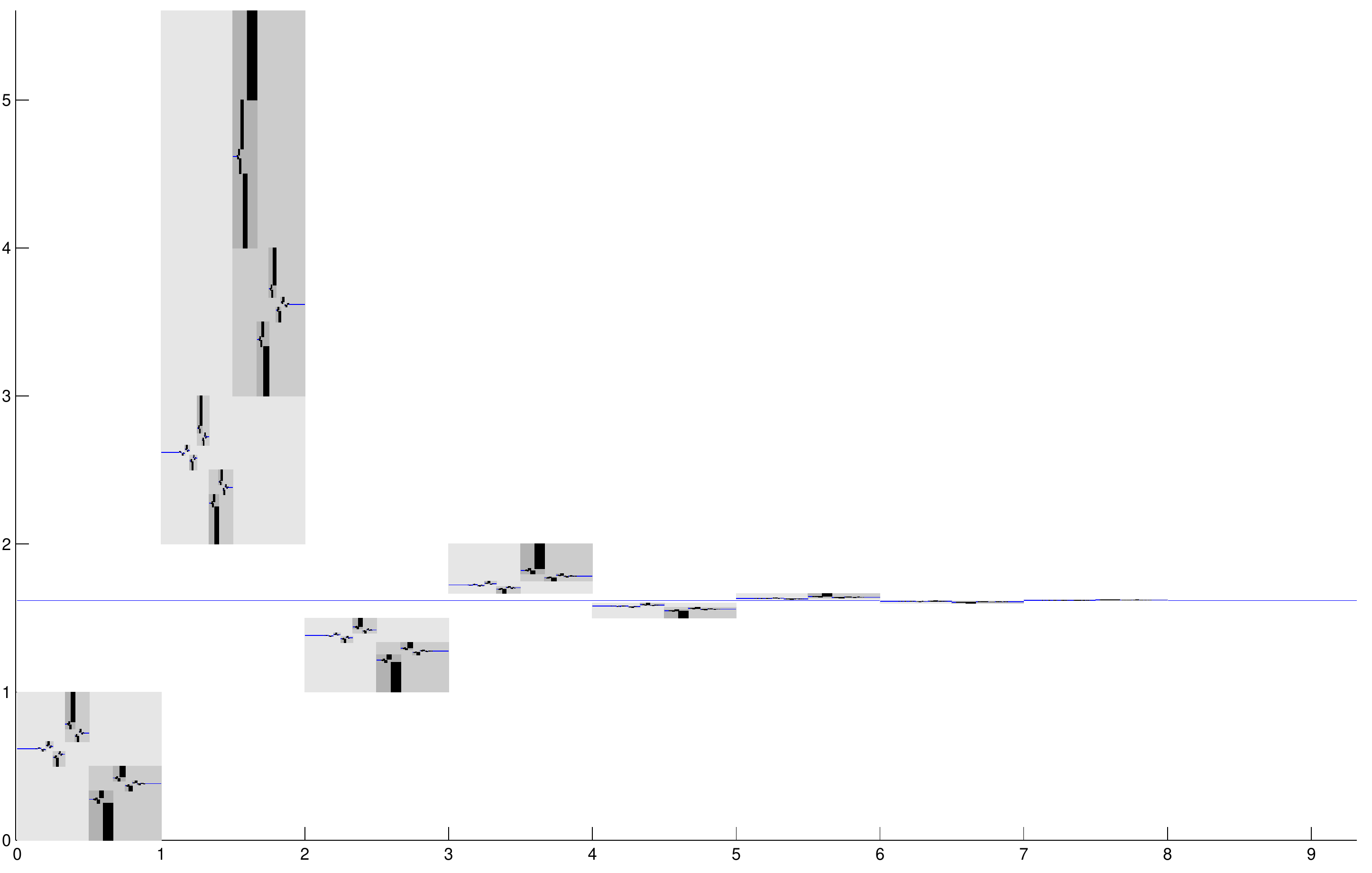}}\\
\end{center}

The graph on the negative sector can be found by using the functional equation
$\Jimm(-x)=-1/\Jimm(x)$.}

\section{Action on the quadratic irrationalities}
Since $\Jimm_\F$ sends eventually periodic paths to eventually periodic paths, 
$\Jimm_\R$ preserves the real-multiplication set. This is the content of our next result:

\begin{theorem}\label{orbits}
$\Jimm$ defines an involution of the set of real quadratic irrationals, (including rationals) and it respects the $\pgl$-orbits.
\end{theorem}

Beyond this theorem, we failed to detect any further arithmetic structure or formula relating the quadratic number to its $\Jimm$-transform.  
All we can do is to give some sporadic examples, which amounts to exhibit some special real quadratic numbers with known continued fraction expansions. We shall do this below. Before that, however, note that the equation below is solvable for every $n\in \Z$:
$$
\Jimm x=n+ x\implies  x=\Jimm\Jimm x=\widetilde{T}^n\Jimm x\implies  x=\widetilde{T}^n( x+n)
$$
$$
\frac{F_{n+1}( x+n)+F_n }{F_n( x+n) +F_{n-1}}= x=\frac{F_{n+1} x+F_n+nF_{n+1} }{F_n x^n +F_{n-1}+nF_n}
$$
$$
\implies F_n x^2+(F_{n-1}+nF_n-F_{n+1}) x -(F_n+nF_{n+1})
$$
More generally, the equation $\Jimm x=M x$ is solvable for $M\in\pgl$. Indeed one has
$$
\Jimm x=M x\implies 
 x=\Jimm\Jimm x=(\Jimm M)(\Jimm x)=
(\Jimm M) M x,
$$
and the solutions are the fixed points of $(\Jimm M) M$. 
These points are precisely the words represented by the infinite path
\begin{eqnarray}\label{fixed}
x=(\Jimm M) M(\Jimm M) M(\Jimm M) M(\Jimm M) M\dots
\end{eqnarray}
where we assume that both $M$ and $\Jimm M$ are expressed as words in $U$ and $T$.
Note that these words are precisely the points whose conjugacy classes remain stable under $\Jimm$.
$$
\Jimm x=M(\Jimm M) M(\Jimm M) M(\Jimm M) M(\Jimm M) M\dots
$$
One may also view these numbers $x$ as sort of eigenvectors of the $\Jimm$ operator:
$\Jimm(x)=M\Jimm(x)$. 

Hence the $\pgl$-orbits of the points in (\ref{fixed}) remain stable under $\Jimm$. Since the set of these orbits is the moduli space of real lattices, we deduce the result:
\begin{proposition}
The fixed points of the involution $\Jimm$ on the moduli space of real lattices are precisely the points (\ref{fixed}) .
\end{proposition}
For example, if $\Jimm M=\widetilde{T}^{k+1}$  then  $M=T^{k+1}$, and the point (\ref{fixed})  is nothing but the point 
$[\overline{1_k,k+2}]$ mentioned in the introduction. 

\bigskip\noindent
\paragraph{Example.} In  \cite{manfred} it is shown that 
$$
x=[0;{\overline{1_{n-1},a}}]=\frac{a}{2}\left(\sqrt{1+4\frac{aF_{n-1}+F_{n-2}}{a^2F_n}}-1\right)
$$
Hence
$
\Jimm(x) = [0;n,\overline{a-2,n+1}], and 
$
so 
$
\Jimm(x)=({n+\frac{1}{y}})^{-1},
$
where 
$$
y=\frac{(a-2)(n+1)+\sqrt{(a-2)^2(n+1)^2+4(a-2)(n+1)}}{2(n+1)}.
$$

\paragraph{Example.} (From Einsiedler \& Ward \cite{ward}, Pg.90, ex. 3.1.1)
This result of McMullen from~\cite{mcmullen} illustrates how $\Jimm$ behaves on one real quadratic number field.
$\Q(\sqrt{5})$ contains infinitely many elements with a uniform bound on their partial quotients, since 
$
[\overline{1_{k+1},4,5,1_k,3}]\in \Q(\sqrt{5}), \quad \forall k=1,2,\dots.
$
Routine calculations shows that the transforms 
$
\Jimm([\overline{1_{k+1},4,5,1_k,3}])=[k+2, \overline{2,2,3,k+2,1,k+3}]
$
lies in different quadratic number fields for different $k$'s. Hence, not only $\Jimm$ does not preserve the property of ``{belonging to a certain quadratic number field}", it sends elements from one quadratic number field to different quadratic number fields.
The following result provides even more examples of this nature: 
\begin{theorem} (McMullen~\cite{mcmullen})
For any $s>0$, the periodic continued fractions
$$
x_m=[\overline{(1,s)^m,1,s+1,s-1, (1,s)^m,1,s+1,s+3)}]
$$
lie in $\Q(\sqrt{s^2+4s})$ for any $m\geq 0$.
\end{theorem}
The fact recently shown by McMullen is that $\Q(\sqrt{d})$ contains infinitely elements with partial quotients bounded by some number $M_d$. He also raises the question: can one take $M_d=2$? 

\nt{Study the derivatives of $\Jimm$ at these points. In general, study the derivative of $\Jimm$ at quadratic irrationalities (and Markov numbers as stated elsewhere)}

\unut{
\paragraph{Example 4 (Powers of the golden ratio)}

\nt{I don't trust these results entirely, it wont do harm to do a check}
\begin{theorem}(Fishman and Miller \cite{fishman2013closed})
One has the following
$$
\Phi^k=
\begin{cases}
[\overline{F_{k+1}+F_{k-1}}],& \mbox{\it if k is even},\\
[F_{k+1}+F_{k-1}-1, \overline{1,F_{k+1}+F_{k-1}-2}]& \mbox{\it if k is odd}.
\end{cases}
$$
\end{theorem}
\begin{corollary}
The continued fraction of the kth power of the golden section is
$$
\Jimm(\Phi^k)=
\begin{cases}
[1_{F_{k+1}+F_{k-1}-1},\overline{2, 1_{F_{k+1}+F_{k-1}-2}}],& \mbox{\it if k is even},\\
[1_{F_{k+1}+F_{k-1}-2},\overline{3,1_{F_{k+1}+F_{k-1}-4}}]& \mbox{\it if k is odd}.
\end{cases}
$$
\nt{determine both sides as a quadratic surd}
\end{corollary}
The following generalization of this is also known:
\begin{theorem}(Fishman and Miller \cite{fishman2013closed})
Consider the recurrence relation 
$$
g_{n+1}=mg_n+lg_{n-1}, \mbox{ with } g_0=0, g_1=1, \quad, m\in \Z_{>0},\quad l=\pm1.
$$
Let
$$
\Phi_{m,l}:=\lim_{n\rightarrow\infty}\frac{g_n}{g_n-1}=
\frac{m\pm\sqrt{m^2+4l}}{2}.
$$
Then for any positive integer $k$, if $l=1$ one has 
$$
\Phi_{m,1}^k=
\begin{cases}
[\overline{g_{k+1}+g_{k-1}}],& \mbox{\it if k is odd},\\
[g_{k+1}+g_{k-1}-1, \overline{1,g_{k+1}+g_{k-1}-2}]& \mbox{\it if k is odd},
\end{cases}
$$
while if $l=1$ one has 
$$
\Phi_{m,-1}^k=
[g_{k+1}-g_{k-1}-1, \overline{1,g_{k+1}-g_{k-1}-2}] \quad \mbox{\it  if k is odd},
$$
\end{theorem}
We leave it to the reader to determine the $\Jimm$-transform of $\Phi_{m,l}^k$.}

\paragraph{Example.}  Suppose $\alpha=p+\sqrt{q}$ be a quadratic irrational with $p,q\in \Q$. Then $\alpha^*=1/\alpha$ provided
$$
|\alpha|^2=p^2-q=1 \implies q=p^2-1.
$$
Suppose $\alpha$ is of this form, i.e. $\alpha=p+\sqrt{p^2-1}$ and suppose $\Jimm(\alpha)=x+\sqrt{y}$. 
Then since 
$$
x-\sqrt{y}=\Jimm(\alpha)^*=\Jimm(\alpha^*)=\Jimm(1/\alpha)=1/\Jimm(\alpha)=1/(x+\sqrt{y})
\implies x^2-y=1,
$$
and we conclude that $\Jimm(\alpha)$ is again of the form $x+\sqrt{x^2-1}$. 

On the other hand, suppose $\Jimm(\sqrt{q})=x+\sqrt{y}$. Then 
$$
\Jimm(\sqrt{q})^*=\Jimm(\sqrt{q}^*)=\Jimm(-\sqrt{q})=-1/\Jimm(\sqrt{q})\implies
|\Jimm(\sqrt{q}^*)|=x^2-y=-1
$$
Hence, $\Jimm(\sqrt{q})$ is of the form $x+\sqrt{1+x^2}$. Conversely, 
if $\alpha$ is of the form $p+\sqrt{p^2+1}$, then $\Jimm(\alpha)$ is a quadratic surd.

\sherh{because
$$
a-2+\frac{1}{n+1+\frac{1}{y}}=y\implies a-2+\frac{y}{(n+1)y+1}=y \implies 
$$$$
(a-2)[(n+1)y+1]=(n+1)y^2\implies (n+1)y^2-(a-2)(n+1)y-(a-2)=0\implies
$$}

\sherh{
\bigskip\noindent
\paragraph{\bf Example 1} 
$$
\Jimm([\overline{n}])=[1_{n-1},\overline{2,1_{n-2}}] \implies
\Jimm(\frac{n+\sqrt{n^2+4}}{2})=\dots
$$

\bigskip\noindent
\paragraph{\bf Example 2} 
$$
\Jimm([\overline{n,m}])=[1_{n-1},\overline{2,1_{n-2},2,1_{m-2}}] \implies
\Jimm(\frac{xxx+\sqrt{xxx^2+4}}{2})=\dots
$$

In the same place it is indicated that the cf expansions of
$$
^k\!\sqrt{\frac{F_{n+k}}{F_n}}
$$
starts with at least $k$ ones's. Hence its jimm starts with at least $n$.

the page below is also helpful to get results of this kind
\url{http://math.arizona.edu/~ura-reports/993/miller.justin/Miller.html}

{\bf something that might be useful}
\begin{lemma} \cite{DiSario}
Let $F(n)=F_{F_n}$. Then 
$$
F(n)=\frac{(F(n-1))^2-(-1)^{F_n}(F(n-2))^2}{F(n-3)}
$$
\end{lemma}}

\sherh{
{\bf Integration}
In this chapter we do some hard analysis.
$$
[0,1]=\bigcup_{n=1}^{\infty} \frac{1}{n+[0,1]}=
\bigcup_{n=1}^{\infty} [\frac{1}{n+1}, \frac{1}{n}]
$$
$$
\int_0^1 \Jimm(x) dx= \sum_{n=1}^\infty \int_{1/(n+1)}^{1/n} 
\Jimm(x) dx 
$$
Set $x=1/(n+u)$  $\implies$
$$
=\sum_{n=1}^\infty  \int_{1}^{0} \Jimm(\frac{1}{n+u}) \frac{-du}{(n+u)^2} =
\sum_{n=1}^\infty  \int_{1}^{0} \frac{1}{\Jimm(n+u)} \frac{-du}{(n+u)^2}
$$
$$
 =
\sum_{n=1}^\infty  \int_{0}^{1} 
\frac{F_n\Jimm(u)+F_{n-1}}{F_{n+1}\Jimm(u)+F_n} \frac{du}{(n+u)^2}
 =
\sum_{n=1}^\infty  \frac{F_n}{F_{n+1}}\int_{0}^{1} 
\frac{\Jimm(u)+F_{n-1}/F_n}{\Jimm(u)+F_n/F_{n+1}} \frac{du}{(n+u)^2}
$$
$$
 =
\sum_{n=1}^\infty  \frac{F_n}{F_{n+1}}\int_{0}^{1} 
\left(1+\frac{F_{n-1}/F_n-F_n/F_{n+1}}{\Jimm(u)+F_n/F_{n+1}} \right)\frac{du}{(n+u)^2}
$$
$$
 =
\sum_{n=1}^\infty  \frac{F_n}{F_{n+1}}\int_{0}^{1} \frac{du}{(n+u)^2}
+
\sum_{n=1}^\infty  \frac{(-1)^n}{F_{n+1}^2}\int_{0}^{1} 
\frac{1}{\Jimm(u)+F_n/F_{n+1}} \frac{du}{(n+u)^2}
$$
$$
 =
\sum_{n=1}^\infty  \frac{F_n}{F_{n+1}}\frac{1}{n(n+1)}
+
\sum_{n=1}^\infty  \frac{(-1)^n}{F_{n+1}^2}\int_{0}^{1} 
\frac{1}{\Jimm(u)+F_n/F_{n+1}} \frac{du}{(n+u)^2}
$$
Integral of $\Jimm$.
$$
\int_0^1\Jimm(x)dx=?
$$

$$
[0,1]=\bigcup_{n=1}^{\infty} \frac{1}{n+\frac {1}{m+[0,1]}}=
\bigcup_{n,m=1}^{\infty} 
\left[\frac{1}{n+\frac{1}{m+1}}, \frac{1}{n+\frac{1}{m}}\right]
$$
$$
\int_0^1\Jimm(x)dx=
\sum_{n,m=1}^\infty \int_{1/(n+\frac {1}{m})}^{1/(n+\frac {1}{m+1})}\Jimm(x)dx
$$
Set $x=1/(n+\frac{1}{m+u})$ $\implies$ 
$\frac{m+u}{mn+nu+1}$ $\implies$ $dx=\frac{du}{(mn+nu+1)^2}$
$$
\sum_{n,m=1}^\infty \int_{{1}/(n+\frac {1}{m})}^{1/(n+\frac {1}{m+1})}\Jimm(x)dx=
\sum_{n,m=1}^\infty \int_{0}^{1}
\Jimm\bigl(\frac{1}{n+\frac {1}{m+u}}\bigr)\frac{du}{(mn+nu+1)^2}
$$
$$
\sum_{n=1}^\infty  \int_{0}^{1} 
\frac{F_n\frac{F_m\jimm(u)+F_{m-1}}{F_{m+1}\jimm(u)+F_m}
+F_{n-1}}{F_{n+1}\frac{F_m\jimm(u)+F_{m-1}}{F_{m+1}\jimm(u)+F_m}
+F_n} \frac{du}{(mn+nu+1)^2}
$$
$$
\sum_{n=1}^\infty  \int_{0}^{1} 
\frac{(F_nF_m+F_{n-1}F_{m+1})\Jimm(u)+(F_nF_{m-1}+F_mF_{n-1})}
{(F_{n+1}F_m+F_{m+1}F_n)\Jimm(u)+(F_{n+1}F_{m-1}+F_nF_m)}
\frac{du}{(mn+nu+1)^2}
$$
$$
\sum_{n=1}^\infty  
\frac{F_nF_m+F_{n-1}F_{m+1}}{F_{n+1}F_m+F_{m+1}F_n}
\int_{0}^{1} 
\frac{\Jimm(u)+(F_nF_{m-1}+F_mF_{n-1})/(F_nF_m+F_{n-1}F_{m+1})}
{\Jimm(u)+(F_{n+1}F_{m-1}+F_nF_m)/(F_{n+1}F_m+F_{m+1}F_n)}
\frac{du}{(mn+nu+1)^2}
$$
This is f.n hard.

$$
-\int_{0}^{1} \frac{du}{(mn+nu+1)^2}=\frac{1}{n(mn+nu+1)}]_{0}^{1} 
=\frac{1}{n(mn+1)}-\frac{1}{n(mn+n+1)}=
$$$$
\frac{1}{(mn+1)(mn+n+1)}
$$

$$
\sum_{n=1}^\infty  
\frac{F_nF_m+F_{n-1}F_{m+1}}{F_{n+1}F_m+F_{m+1}F_n}\frac{1}{(mn+1)(mn+n+1)}=?
$$
There are several possibilities for computing this integral:

First, by using the partition $\{1/n\}$. This involves infinitely iterated infinite sums and appears to be hard.

It seems that a better possibility is to use the Farey partition or the standard partition. 

Yet another possibility is to use an irrational partition. I mean the noble numbers. This way the sums might become rational sums.

Next: the signed measure $d\Jimm$. I believe this is entirely singular, i.e. it has no continuous part. Hence integrating it against some function must yield over a sum over the jumps.

Next:
$$
\waw(t):=\int_0^t\Jimm(x)dx=?
$$

Next: Moments
$$
\waw(s):=\int_0^1 x^s\Jimm(x)dx=?
$$
$$
\int_0^1 x^sd\Jimm=?
$$

Next: Fourier transforms
$$
\waw(t):=\int_0^1 \exp{2\pi ixt}\Jimm(x)dx=?
$$
Next: Fourier transforms
$$
\int_0^1 \exp{2\pi ixt}d\Jimm=?
$$
}

\section{Statistics}
By the Gauss-Kuzmin theorem, for almost all real numbers $x$, the frequency of an integer $k>0$ in the continued fraction expansion of $x$ converges to
$$
\frac{1}{\log 2}\log \left(1+{1\over k(k+2)}\right)
$$ 
(see~\cite{iosifescu2002metrical}). We say that $x$ obeys the Gauss-Kuzmin statistics if this is the case.

\sherh{WHAT FOLLOWS IS PLAINLY WRONG
\subsection{The initial string of 1's.}
One simple question we may answer by a straightforward application of  the Gauss-Kuzmin theorem is the following.
Suppose that $X$ is a random variable uniformly distributed over the interval $[0,1]$. What is the probability $p_k$ that $\Jimm(X)$ starts by 
$k+1$?
In other words, what is the length of the string that $T_\jimm$ forgets when applied to $X$?  
Let $q_k$ be the probability that $X$ belongs to $E_k:=[0,1_k, \{1\leq x\leq \infty\}]$.
Note that $\{E_k\}$ is a decreasing sequence of events, and the probability we are looking for is 
$$
p_k=\mathbf{Prob}(X\in E_{k}\setminus E_{k+1})=q_k-q_{k+1}.
$$ 
One has  $E_k=\{F_{k+1}/F_{k+2}\leq x\leq  F_{k}/F_{k+1}\}\implies$
$$
q_k=\mathbf{Prob}(X\in E_{k})=\frac{1}{\log 2}\int_{F_{k+1}/F_{k+2}}^{F_{k}/F_{k+1}} \frac{1}{1+x} dx
=\frac{1}{\log 2}\log \frac{F_{k+1}F_{k+3}}{F_{k+2}^2}
$$ 
So we have
$$
p_k=\frac{1}{\log 2}\log 
\frac{F_{k+1}}{F_{k+4}}
\frac{F_{k+3}^3}{F_{k+2}^3}
$$
One has $\sum_0^\infty p_k=1$ as  one may check.

%
%
%
%
}
\sherh{

{\bf Exercise}. Given two probability distributions $p(n)$, $q(n)$ on $\Z_{>0}$, find $x=[0, n_1, n_2, \dots]\in [0,1]$ 
with $\Jimm x=[0, m_1, m_2, \dots]$
such that 
$$
\lim_{k\to\infty} { \sharp \{1\leq i\leq k\, :\, n_i=n\} \over k}=p(n)
$$
\mbox{ and }
$$
\lim_{k\to \infty} { \sharp \{1\leq i\leq k\, :\, m_i=m\} \over k}=q(m)
$$
}
If $X\in [0,1]$ is a uniformly distributed random variable, then almost everywhere the arithmetic mean of its partial quotients tends to infinity, i.e. if $X =[0, n_1, n_2, \dots]$ then 
\begin{equation}
\lim_{k\to \infty} {n_1 +\cdots +n_k \over k} = \infty\qquad\hbox{(a.s.)}
\end{equation}
(see \cite{iosifescu2002metrical}). In other words, the set of numbers in the unit interval such that the above limit is infinite, is of full Lebesgue measure. Denote this set by $A$.

Now since the first $k$ partial fractions of $X$ give rise to at most $n_1 +\cdots +n_k -k$ partial fractions of $\Jimm(X)$ and 
at least $n_1 +\cdots +n_k -2k$ of these are 1's, one has 
$$
{n_1 +\cdots +n_k -k \over  n_1 +\cdots +n_k -2k} \to 
{\frac{n_1 +\cdots + n_k}{k} -1 \over  \frac{n_1 +\cdots + n_k}{k} -2} \to 1
$$
\begin{lemma}\label{density}
The density of $1$'s in the continued fraction expansion of $\Jimm(X)$ equals 1 a.e., and therefore the
the continued fraction averages of $\Jimm(X)$ tend to 1 a.e.
\end{lemma}

We conclude that $\Jimm(A)$ is a set of zero measure. 
This says that, even though the involution $\Jimm$ is a fundamental symmetry of 
$\R$, many properties of reals are strongly biased under it, unlike the other fundamental involutions $U$, $V$ and $K$.

\sherhh{What is the probability that the continued fraction of $\Jimm(x)$ is bounded? It seems 1.. There must be other questions of this type, for the experimental paper. Forex: if we simply ignore 1's in $\Jimm(x)$ (either by simply deleting or by applying the more delicate process of collapsing the continued fraction, i.e. subtracting  all 1's from the cf digits and removing the ensuing 0's), then what is the distribution of the continued fraction digits? 
}

\unut{
We leave open the questions such as: what is the distribution of the continued fraction entries of $\Jimm(X)$, if we ignore the 1's? }

\subsection{The derivative of jimm.}  
A function discontinuous on a dense subset of $[0,1]$ can not be differentiable everywhere on the residual set; it can be differentiable at most on a 
meager set (i.e. a countable union of nowhere dense sets), see Fort's paper \cite{fort}. But meager does not mean negligible: There exists meager sets of full Lebesgue measure, and it was also exhibited in the location cited a function discontinuous at rationals and yet differentiable on a set of full measure. It turns out that $\Jimm$ is a function of this kind.

To see why, recall that average values of continued fraction elements tend to infinity almost everywhere. Let $a=[0, n_1, n_2, \dots]$ be a number with this property:
\begin{equation}
\lim_{k\to \infty} {n_1 +\cdots +n_k \over k} = \infty.
\end{equation}
Then for every constant $M$, there is some $k$ with $n_1+\dots+n_k>kM$. But then the $\Jimm$- transform of the initial length-$k$ segment of $x$ is of length at least $kM-k$. Hence if $y$ is any number whose continued fraction expansion coincide with that of $x$ up to the place $k$, then the continued fraction $\Jimm(y)$ coincide with that of $\Jimm(x)$ at least up to the place $kM-k$. Since 
$kM-k$ is arbitrarily big compared to $k$, and since longer continued fractions give exponentially better approximations, we see that $\Jimm(y)$ a.e. is much closer to $\Jimm(x)$ then $y$ is to $x$. Hence the idea of the following theorem.

\begin{theorem}\label{derivativevanish}
The derivative  of $\Jimm(a)$  exists almost everywhere and vanishes almost everywhere.
\end{theorem}

To prove this, we need to show that for almost all $a$, 
$$
\lim_{x\to a} \frac{\Jimm(a)-\Jimm(x)}{a-x} =0.
$$

Assume that $x$ is irrational or equivalently its continued fraction expansion is non-terminating.

Let $x\in [0,1]$ with $0<|x-a|<\delta$ for some $\delta$. Then there is a number $k=k_\delta$, such that the continued fractions of 
$a$ and $x$ coincide up to the $k$th element. 
Hence 
$
x=[0,n_1,n_2, \dots n_k, m_{k+1}, \dots] \mbox{ with } m_{k+1}\neq n_{k+1}.
$
Note that this latter condition also guarantees that $0<|x-a|$. 
Now let 
$$
M_k(z):=[n_1, n_2, \dots,n_{k-1}, n_k+z]=
\frac{\alpha_kz+\beta_k}{\gamma_kz+\theta_k}
$$
and put 
$
a_k:=[0,n_{k+1}, n_{k+2}, \dots], \quad x_k:=[0,m_{k+1}, m_{k+2}, \dots].
$
Then one has $0<a_k<1$ (with strict inequality since $a$ is irrational) and $0\leq x_k<1$ for every $k=1,2,\dots$. One has
$$
a=M_k(a_k),\quad x=M_k(x_k)  \mbox{ and } \det(M_k)=(-1)^k.
$$
\begin{lemma} Let $a:=[0,n_0, n_1, \dots]$ and suppose that the continued fractions of $a$ and $x$ coincide up to the place $k$ (but not $k+1$), where $x\in [0,1]$. Put $N_k:=\sum_{i=1}^k n_i$, and $\mu_k:=N/k$. Then 
$$
|a-x|>\frac{1}{24} (2\mu_{k+3})^{-2(k+3)}
$$
\end{lemma}
{\it Proof.} One has 
$$
|a-x|=\left|M_k(a_k)-M_k(x_k)\right|=
\left|\frac{\alpha_ka_k+\beta_k}{\gamma_ka_k+\theta_k}-\frac{\alpha_kx_k+\beta_k}{\gamma_kx_k+\theta_k}\right|
=\frac{1}{(\gamma_ka_k+\theta_k)(\gamma_k x_k+\theta_k})
$$
Since 
$$
M_{i+1}(z)=M_i\left(\frac{1}{n_{i+1}+z}\right)=\frac{\beta_i z+(\alpha_i+n_{i+1}\beta_i)}{\theta_i z+(\gamma+n_{i+1}\theta_i)},
$$
one has $\gamma_{i+1}=\theta_i$ and $\theta_{i+1}=\gamma_i+n_{i+1}\theta_i$.
Hence $\theta_{i+1}>\gamma_{i+1} \implies \theta_{i+1}> \theta_i(1+n_{i+1})$.
This implies 
\begin{eqnarray*}
\theta_i<(1+n_1)(1+n_2)\dots (1+n_i),\\
\gamma_i<(1+n_1)(1+n_2)\dots (1+n_{i-1}).
\end{eqnarray*}
Since $0 \leq a_k, x_k <1 $, this implies
\begin{eqnarray*}
\gamma_ka_k+\theta_k<\gamma_k+\theta_k<2(1+n_1)(1+n_2)\dots (1+n_k),\\
\gamma_kx_k+\theta_k<\gamma_k+\theta_k<2(1+n_1)(1+n_2)\dots (1+n_k).
\end{eqnarray*}
Hence, we get 
$$
|a-x|>\frac{|a_k-x_k|}{4(1+n_1)^2(1+n_2)^2\dots (1+n_k)^2}
$$
To estimate $|a_k-x_k|$, consider 
$$
a_k-x_k=\frac{1}{n_{k+1}+a_{k+1}}-\frac{1}{m_{k+1}+x_{k+1}}=
\frac{m_{k+1}-n_{k+1}+x_{k+1}-a_{k+1}}{(n_{k+1}+a_{k+1})(m_{k+1}+x_{k+1})}
$$
$$
>\frac{m_{k+1}-n_{k+1}+x_{k+1}-a_{k+1}}{(1+n_{k+1})(1+m_{k+1})}.
$$
Now, if $m_{k+1}<n_{k+1}$ then set $m_{k+1}=n_{k+1}-t$ with $t\geq 1$. Then one has 
$$
|a_k-x_k|>
\frac{|-t+x_{k+1}-a_{k+1}|}{(1+n_{k+1})(1+n_{k+1}-t)}>\frac{a_{k+1}}{(1+n_{k+1})^2}
$$
On the other hand, if $3n_{k+1}\geq m_{k+1}>n_{k+1}$ then 
$$
|a_k-x_k|>
\frac{|1+x_{k+1}-a_{k+1}|}{(1+n_{k+1})(1+3n_{k+1})}>\frac{1-a_{k+1}}{3(1+n_{k+1})^2},
$$
and if $m_{k+1}>3n_{k+1}$ then
$$
|a_k-x_k|=
\frac{1-\frac{n_{k+1}}{m_{k+1}}+\frac{x_{k+1}}{m_{k+1}}-\frac{a_{k+1}}{m_{k+1}}}%
{(1+n_{k+1})(1+\frac{1}{m_{k+1}})}>\frac{1}{6(1+n_{k+1})}
$$
So one has 
$$
|a_k-x_k|>
\frac{a_{k+1}(1-a_{k+1})}{6(1+n_{k+1})^2},
$$
which gives the estimation from below
$$
|a-x|>\frac{a_{k+1}(1-a_{k+1})}{24(1+n_1)^2(1+n_2)^2\dots (1+n_k)^2(1+n_{k+1})^2},
$$
estimation obtained under the assumption that the continued fraction expansions of $x$ and $a$ coincide up until the 
$k$th term and differ for the $k+1$th term.

Now we have the crude estimate
$$
{\cfrac{1}{n_{k+2}+\cfrac{1}{n_{k+3}+1}}}>a_{k+1}>{1 \over 1+n_{k+2}}
\implies
{a_{k+1}(1-a_{k+1})}>{1 \over (1+n_{k+2})^2 } {1 \over (1+n_{k+3})^2}.
$$
which gives
$$
|a-x|>\frac{1}{24(1+n_1)^2(1+n_2)^2\dots (1+n_{k+2})^2(1+n_{k+3})^2},
$$
Now put $N_k:=\sum_{i=1}^k n_i$, and $\mu_k:=N/k$. Then 
$$
(1+n_1)^2(1+n_2)^2\dots (1+n_{k})^2 \leq (1+\mu_k)^{2k}\leq  (2\mu_k)^{2k}
$$
The last inequality follows from the fact that $\mu_k\geq 1$ for all $k$, since $n_i\geq 1$ for all $i$. 
We finally obtain the estimate
$$
|a-x|>\frac{1}{24} (2\mu_{k+3})^{-2(k+3)} =\frac{1}{24}\exp\{-2(k+3) \log 2\mu_{k+3}\}
$$
\hfill $\Box$

On the other hand, if the c.f. expansions of $a$ and $x$ coincide up to the $k=k(x)$th place, 
then the c.f. expansions of $\Jimm(a)$  and $\Jimm(x_i)$ coincide up to the place $N_k$, and by~(\ref{binet}) we have
$$
|\Jimm(a)-\Jimm(x)|<F_{N_k}^{-2}<{\sqrt{5}}\phi^{-2N_k}={\sqrt{5}} \exp\{-2k\mu_k \log\phi\}
$$
(This estimate should be close to optimal (a.e.), since the density of 1's in the c.f. expansion of $\Jimm(a)$ is 1 (a.e.) by Lemma~\ref{density}.)
This gives
$$
\left|\frac{\Jimm(a)-\Jimm(x)}{a-x}\right|<
{24\sqrt{5}}\exp k \{2(1+3/k) \log 2\mu_{k+3}-2\mu_k \log\phi \}\implies
$$
$$
\qquad\qquad\qquad \left|\frac{\Jimm(a)-\Jimm(x)}{a-x}\right|<
A\exp \Bigl\{2k \log \phi  
\bigl(
B  \log 2\mu_{k+3}-\mu_k
\bigr) 
\bigr\}
$$
where $A$ is some absolute constant and $B=(1+3/k)/ \log\phi $ can be taken arbitrarily close to $1/\log\phi<2.08$ 
by assuming $k$ is big enough.

We see immediately that, if $a=[0,n,n,n,n,\dots]$ then $\mu_k$ is constant $=n$, and if $n$ is taken big enough so that 
$2.08\log 2n-n<0$, then the derivative exists and is zero. This is true for $n>4$. We don't claim that our estimations are optimal in this respect, however. 

On the other hand, since $\mu_k\to\infty$ almost surely, we see that $B  \log 2\mu_{k+3}-\mu_k<0$ for $k$ sufficiently big
and the derivative exists and vanishes. This is because by choosing a sufficiently small neighborhood $\{|x-a|<\delta\}$, we 
can guarantee that  $k=k(x)$ is always greater than a given number for any $x$ in this neighborhood. This concludes the proof of the theorem.

Note that, if $\mu_k\to\infty$ then the average of continued fraction elements of $\Jimm(a)$ tends to 1, 
and $\Jimm$ is not differentiable at $\Jimm(a)$. In other words, $\Jimm$ is almost surely not differentiable at $\Jimm(a)$.
In the same vein, the derivative of $\Jimm$ at $a=[0,n,n,n,n,\dots]$ vanish for $n>4$, and we see that 
$\Jimm$ is not differentiable at $\Jimm(a)=[0,1_{n-1}, \overline{2,n-2}]$ or at best it will be of infinite slope at this point.

\unut{
It might be of interest to study the derivatives at quadratic irrationalities in general. 
An exciting possibility would be to establish a connection with the derivative of $\Jimm$ and the Markov irrationalities, \cite{aigner2013markov}.
We don't know if the derivative may attain or not a finite non-zero value at some point.}

\sherh{Concerning the Markov theory, there is the following fact. The set $\{\alpha: L(\alpha)<3\}$ in the Lagrange spectrum consists of numbers whose continued fractions consists exclusively of 1's and 2's. Now this is hard to describe as a subset of the reals. However, the $\Jimm$-transform of this set is simply the set of numbers whose continued fractions do not contain 1. This latter set is easy to draw and understand. In every Farey interval, we don't take one half of the interval.}

\sherhh{
\bigskip\noindent
{\bf Exercice 1.} Show that the derivative at $\sqrt{2}$ and $1/\phi$ does not exist. Construct other examples with this property.
Hint. One may try things like $[\overline{n, 1_{n^2}}]$. Here $\mu_n\to 1$. Is it true that if $\mu_n\to 1$, then the derivative does not exist?

\smallskip\noindent
{\bf Exercice 2.} Is it true that, if the derivative of $\Jimm$ exists and has a non-zero value $r$ at a point $a$, then the derivative of $\Jimm$ at the point $\Jimm(a)$ exists and equals $1/r$?}

\sherh{
\subsection{Analytic properties of other tree automorphisms}
Let $\phi\in Aut_I(|\F|)$ be an automorphism. 
\begin{itemize}
\item What can be said about the continuity of the induced map on $\R$ in general?
\item What can be said about the differentiability of the induced map? 
\item Given a point on the boundary, can you construct an automorphism which admits some derivatives at that point? 
Can you arrange the value of the derivatives? 
\end{itemize}
{\bf Conjecture.}  Let $\tau_\mu$ be a shuffle (or a twist). If its locus $\mu$ is invariant under $Aut_I(|\F|)$, then $\tau_\mu$ is similar to $\Jimm$, in that it is continuous on irrationals, differentiable a.e. with null derivative a.e. Note that the $Aut_I(|\F|)$-invariance of $\mu$ means that if a vertex belongs to $\mu$, then every vertex of the same distance to the base also belongs to $\mu$. In other words, $\mu$ is ``{radially symmetric}". Let us call the corresponding automorphisms ``{radially symmetric}". 
Then  I suspect that the center of $Aut_I(|\F|)$ consists of radially symmetric automorphisms.
}

\unut{\section{Concluding remarks.} \label{concluding}
\subsection{The Minkowski measure and the Denjoy measures.} 
The reader with a taste for singular functions such as $\Jimm$ will inevitably ask how it relates to Minkowski's question mark function. 
In our setting, this latter function is best understood as the cumulative distribution function of the  unique $Aut_I(|\F|)$-invariant measure on the boundary $\partial \F$.

More generally\footnote{In a yet more general setting, \nt{which comprises the Lebesgue measure}, 
take any function $P:V_\bullet(\F)\to [0,1]$ and 
let the random walker, upon arriving to a vertex, choose to turn right with probability $P(v)$ and to turn left with probability $1-P(v)$. Then this introduces a probability measure $\mu_P$ on $\partial\F$. A special class of measures is obtained, by requiring $P$ to be only a function of the distance to the starting edge. },  let $p\in (0,1)$ be a real number. 
Suppose that a random walker on $\F$ starting at the edge $I$ decides with probability $1/2$ to move through the positive sector of $\F$ and then each time he arrives to a bifurcation, he chooses to turn right with probability $p$ and he chooses to turn left with probability $q:=1-p$. Similarly, he decides, with probability $1/2$ to move through the negative sector and then proceeds in the same manner at bifurcations. The probability that the walker ends up in the 
Farey interval ${\mathcal O}_{e'}$ is the product of probabilities of choices he makes to go from the base edge $I$ to the edge $e'$. 

\nt{If we drop the condition that the probabilities add up to one at junctures, then it seems that we can obtain all (probability or not) measures on the circle, including the usual Lebesgue measure. The set of all measures thus form an infinite dimensional space akin to Penner's universal Teichm\"uller space. Of course, the action of Thompson's group might be relevant here, as a sort of universal mapping class group.} 

\nt{Since the Lebesgue measure is characterized by its invariance properties under translations and it also have rescaling property with respect to multiplication, and since this measure can be obtained by the above procedure, one is tempted to ask if it might be possible to define the addition and multiplication by using the Lebesgue measure}

Note that, unlike the classical random walkers, ours is a non-backtracking walker. He is not allowed to retrace his steps (so this is not a Markov process) \cite{ortner2008online}.
Since by definition the topology of the boundary is generated by the Farey intervals, this puts a measure on the Borel algebra of $\partial\F_I$ and since the quotient map $\partial\F_I\to S_I^1$ is measurable with respect to Borel algebras, we obtain a probability measure $\mu_{p}$ on $S_I^1$, the {\it Denjoy measure} first introduced and studied by Denjoy \cite{denjoy1938fonction}. Since $0<p<1$, the set of rationals is of measure 0, so that this passage to the quotient has no detectable effect on the measure space properties. Let ${\mathbf F}_{p}(x)$ be the cumulative distribution function of $\mu_{p}$. Then ${\mathbf F}_{1/2}(x)$ is precisely the right extension of the Minkowski question mark function to $\widehat{\R}$ and provides a canonical homeomorphism $S_I^1\to [0,1]/0\!\sim\!1$ (recall that $\widehat{\R}$ is identified with $S_I^1$ by the continued fraction map $cfm$ and we may consider $\mu_p$ as a measure on $\widehat{\R}$ via this identification). As such, it conjugates the modular group to a group of homeomorphisms of $[0,1]/0\!\sim\!1$. The functional equations for the question mark function (see Alkauskas' papers \cite{alkauskas2011semi}, \cite{alkauskas2010minkowski}) then becomes a mere expression of this action. In fact, the conjugate representation of the modular group is a subgroup of Thompson's group $T$ presented as the orientation preserving piecewise linear homeomorphisms of  $[0,1]/0\!\sim\!1$ with break points at dyadic rationals and such that the slopes of linear pieces are all powers of 2, see \cite{imbert1997isomorphisme}. Now ${\mathbf F}_{1/2}(x)$ also conjugates $\Jimm$ to an involution of $[0,1]/0\!\sim\!1$ (with jumps at dyadic rationals), and there are functional equations relating the action of this conjugate involution and the conjugate modular group action, similar to the ones we gave in the beginning of the paper.
There is a possibility that this conjugate $\Jimm$ by the Minkowski function (or more generally by the functions ${\mathbf F}_{p}$) may interact in interesting ways with other singular functions such as the Takagi function\cite{lagarias2011takagi}. See also the papers by Bonano and Isola \cite{bonanno2008orderings}, \cite{bonanno2007spectral} for some more variations on this theme. 

\sherh{On the other hand, the conjugate modular group action is not very convenient (i.e. it is piecewise linear but not linear) which makes 
the conjugate-$\Jimm$ somewhat less amenable for study. E.g. one has $?S?(x)=1/2-x (mod 1)$ and 
$$
?T?(x)=\begin{cases}
3/2x& 0\leq x\leq 1,\\
(1+x)/2 & 1/2\leq x\leq 1
\end{cases}
$$
{\bf Exercice.} Check this. Find also $?V?$. Then find the functional equations for $?\Jimm ?$. \\
(in fact, I did some computations to write $\Jimm$ analytically. They are in one of the latex files)

{\bf Problem.} Study the derivative of the $?$-conjugate $\Jimm$.}

\sherh{One should start by conjugating the modular group and Thompson's groups by the functions ${\mathbf F}_{p}$ and see how do they look like. Even more generally, one could conjugate with arbitrary boundary measures ${\mathbf F}$. What if we conjugate with the Lebesgue measure?}

\sherh{Another special class of measures are obtained as shifts (i.e. subgroups-periodic measures.)}

The measures $\mu_p$ are purely singular \cite{Viader1998} with respect to the Lebesgue measure, i.e. their c.d.f's ${\mathbf F}_p$ have a vanishing derivative almost everywhere, a property shared by the involution $\Jimm$. 
In a different vein, it can be readily verified that the Minkowski measure\footnote{The measure   $\mu_{1/2}$ have also been called the 
{\it Lebesgue measure} of $\partial\F$  (\cite{cohen2004distributions}) and also the {\it Farey measure} \cite{kessebohmer2008fractal}.}
 $\mu_{1/2}$ is an invariant measure for both the Gauss-Kuzming-Wirsing operator 
${\mathscr L}^G_2$ (see the insightful paper of  Vepstas \cite{vepstas2008minkowski}) and for the dual Gauss-Kuzming-Wirsing operator ${\mathscr L}^{\jimm}_2$.

It seems that $\Jimm$ and ${\mathbf F}_{1/2}(x)$ recognize each other in another context: the integrals $\int \Jimm \, d{\mathbf F}_{p}(x)$ may well be finite.

\subsection{Other trees and other continued fraction maps.} In fact, our story is about the group 
$\Z/2\Z*\Z/3\Z=\langle a, b \, |\, a^2=b^3=1\rangle$ and a certain representation of this group as a subgroup of ${\mathrm{PSL}}_2(\R)$, namely the representation whose image is the modular group.
Every representation gives rise to a sort of continued  fraction map, with particularly nice properties if the element $ab$ is sent to a parabolic element. There are variations on this theme concerning the trees related to the free products 
$\Z/p\Z*\Z/q\Z\dots $ and their representations in ${\mathrm{PSL}}_2(\R)$. The stories are particularly nice if the targets of the representations lies inside ${\mathrm{PGL}}_2(\Q)$. Hecke groups (see \cite{koruoglu2010},\cite{koruoglu}) are also promising in this respect. These stories involve various $\Jimm$-like automorphisms and boundary measures.\\

\sherh{There is a $\Jimm$-function for any element of the Teichm\"uller space of the modular curve. In other words, for any representation of the group $\Z_2*\Z_3$ in PSL$_2(\R)$, there is a $\Jimm$ function satisfying the functional equations. The most interesting thing is to start with the representations in PGL$_2(\Q)$.}

\subsection{Covariant functions}
As we already mentioned in the introduction, the functional equations satisfied by $\Jimm_\R$ says that it lax-covariant with the action of $\pgl$ on $\widehat\R$. Is it possible to find some analytic function on the upper half plane satisfying those functional equations (if necessary after modifying them by the complex conjugation)? 

On the other hand, there are some analytic $\psl$-covariant functions on the upper half plane, discovered by 
Brady \cite{brady1971meromorphic} and Smart \cite{smart1972meromorphic}, see the papers by Basraoui and Sebbar  for a recent account
(they use the term ``{equivariant}" instead of ``{covariant}")
\cite{el2012rational} \cite{saber2014equivariant} \cite{sebbar2012equivariant}. The fact is that, if $h$ is an equivariant function (i.e. if $h$ satisfies the functional equations of $\Jimm$), then its Schwartzian is a weight-4 modular form. In the opposite direction,  if $f$ is a weight-$k$ modular form, then the function 
$$
h(z)=z+k\frac{f(z)}{f'(z)}
$$
is $\psl$-covariant. The $\kappa$ function of Kaneko and Yoshida \cite{kaneko2003kappa} is covariant for an infinite subgroup of infinite index $\psl$ (they use the word ``{covariant}"). The Rogers-Ramanujan continued fraction $r(\tau)$ have some weak covariance properties, see \cite{duke2005continued}. Note that the identity function and the complex conjugation are $\psl$- covariant, too.

We  tried to concoct from $\Jimm_\R$ a $\pgl$-invariant (singular) function on $\widehat\R$, but we failed.

\subsection{Concepts similar to Jimm in Engineering}

Digital communications implement some coding techniques to represent a data, and some signal conditioning transforms to obtain some desired properties on these representations. In analog communication when a carrier wave is used to transmit the data, as in the case of RF communications, a modulation of the carrier wave can be used to encode data; a sinusoidal carrier wave can be modulated in terms of its frequency, amplitude or phase. Be it RF, light, sound or electrical signals, the natural media used to transmit the digital information being analog, some discretizing mapping from digitally represented information to analog states is necessary.

Phase shift keying (PSK) is a modulation method to represent data in terms of changes of a sinusoidal carrier wave. When a modulation method assumes discrete values on modulated attribute of carrier wave, as in the case of digital communication, it is called ``{keying}", hence the name phase shift keying. Binary PSK assumes two distinct states of 180 degrees separated phases with equal amplitude and frequency; simplest binary to BPSK mapping is to represent a 0 bit in digital as a specific phase and, a 1 bit as the inverted phase. All signal are subject to mostly arbitrary phase shift during transmission (except in very controlled media like length balanced wires or PCB traces), unless the source and the receiver has another information source to decide which of two phases is the neutral one (i.e. 0 bit), PSK modulation suffers from what's called phase ambiguity. Phase ambiguity can be addressed by supplying a non-modulated reference signal to compare the phase against (like a clock signal in digital electronics) which is not always practical; or by use of another transform on binary data called differential encoding. If a digital data is differentially encoded before PSK modulation, the each bit is not mapped to a specific phase of PSK states, but it is mapped to a specific change in phase. A differentially encoded data has the property called polarity insensitivity, which resolves the phase ambiguity issue (i.e. phone lines work no matter in which order the two copper wires are connected). Differential Manchester Encoding is a well known example of differential encoding which is used as the default encoding scheme of some versions of Ethernet (IEEE 802.3) computer network communications. Manchester encoding can also be thought of a digital counterpart of BPSK, it uses the digital clock signal to modulate the data signal via XOR operation; the resulting signal is same as the clock signal while the data bit (which is optionally differentially encoded before) is 0, and the inverse of clock signal if data bit is 1.

Jimm transformation of a number can be thought of being similar to BPSK; a number represented as zero and ones encoding the path on the farey tree, can be transformed using XOR operation to a modulation with a series of $(01)^\omega = 010101..$. which serves like a carrier signal or digital clock. We face the same phase ambiguity issue as in the case of BPSK, but in this case the problem can be solved by means of convention, such as by mapping the left turns to 0 and right turns to 1, and not the other way, for farey tree encoding, and the carrier signal to be $(01)^\omega $ and not to $(10)^\omega $ thus fixing the neural phase to be carrier series starting with 0. There is no invariant signal (number) under this transformation, as carrier signal contains 1s and XORing with 1 is equivalent to negation which is an unary operator with no fixed point; but the signal (number) closest to be minimally changing under this transformation is $(0011)^\omega $ giving 
$$(0011)^\omega  \xor (01)^\omega  = 011(0011)^\omega $$ which is the same signal 180 degrees out of phase, and the numbers related to $(0011)^\omega $ and $011(0011)^\omega $ are respectively $1+\sqrt2$ and $\sqrt2$.}

\bigskip\noindent
{\bf Acknowledgements.} 
We are grateful to 
Roelof Bruggeman,
Alex Degtyarev, 
Gareth Jones,
Masanobu Kaneko,
Jeffrey Lagarias,
Robert Langlands,
Özer Öztürk,
Recep \c Sahin,
Susumu Tanab\'e,
Michel Waldschmidt,
Masaaki Yoshida,
Mikhael Zaidenberg and
Ayberk Zeytin for various suggestions, critiques and encouragement. 
We are thankful to the participants of the Ankara-İstanbul Algebraic Geometry and Number Theory Meetings who suffered the first exposition of the content of this paper.

This research have been supported by the Galatasaray University Research Grant 14.504.001 
and the T\"UB\.ITAK 1001 Grant No 110T690. 

\small

\bibliographystyle{amsplain}
\bibliography{references2}

\def\cydot{\leavevmode\raise.4ex\hbox{.}}
\providecommand{\bysame}{\leavevmode\hbox to3em{\hrulefill}\thinspace}
\providecommand{\MR}{\relax\ifhmode\unskip\space\fi MR }
\providecommand{\MRhref}[2]{%
  \href{http://www.ams.org/mathscinet-getitem?mr=#1}{#2}
}
\providecommand{\href}[2]{#2}
\begin{thebibliography}{10}

\bibitem{Bugeaud}
Boris Adamczewski and Yann Bugeaud, \emph{Nombres réels de complexité
  sous-linéaire: mesures d'irrationalité et de transcendance}, Journal für
  die reine und angewandte Mathematik (Crelles Journal) \textbf{2011} (2011),
  no.~658, 65--98.

\bibitem{aigner2013markov}
Martin Aigner, \emph{Markov's theorem and 100 years of the uniqueness
  conjecture: a mathematical journey from irrational numbers to perfect
  matchings}, Springer Science \& Business Media, 2013.

\bibitem{alkauskas2010minkowski}
Giedrius Alkauskas, \emph{The minkowski question mark function: explicit series
  for the dyadic period function and moments}, Mathematics of Computation
  \textbf{79} (2010), no.~269, 383--418.

\bibitem{alkauskas2011semi}
\bysame, \emph{Semi-regular continued fractions and an exact formula for the
  moments of the minkowski question mark function}, The Ramanujan Journal
  \textbf{25} (2011), no.~3, 359--367.

\bibitem{bonanno2007spectral}
Claudio Bonanno, Sandro Graffi, and Stefano Isola, \emph{Spectral analysis of
  transfer operators associated to farey fractions}, arXiv preprint
  arXiv:0708.0686 (2007).

\bibitem{bonanno2008orderings}
Claudio Bonanno and Stefano Isola, \emph{Orderings of the rationals and
  dynamical systems}, arXiv preprint arXiv:0805.2178 (2008).

\bibitem{bonanno2014thermodynamic}
\bysame, \emph{A thermodynamic approach to two-variable ruelle and selberg zeta
  functions via the farey map}, Nonlinearity \textbf{27} (2014), no.~5, 897.

\bibitem{brady1971meromorphic}
Michael~M Brady, \emph{Meromorphic solutions of a system of functional
  equations involving the modular group}, Proceedings of the American
  Mathematical Society \textbf{30} (1971), no.~2, 271--277.

\bibitem{cohen2004distributions}
Joel~M Cohen, Flavia Colonna, and David Singman, \emph{Distributions and
  measures on the boundary of a tree}, Journal of mathematical analysis and
  applications \textbf{293} (2004), no.~1, 89--107.

\bibitem{denjoy1938fonction}
Arnaud Denjoy, \emph{Sur une fonction r{\'e}elle de minkowski}, J. Math. Pures
  Appl \textbf{17} (1938), no.~9, 105.

\bibitem{dicks}
Warren Dicks and Martin~John Dunwoody, \emph{Groups acting on graphs}, vol.~17,
  Cambridge University Press, 1989.

\bibitem{djokovic}
Dragomir~Ž Djoković and Gary~L Miller, \emph{Regular groups of automorphisms
  of cubic graphs}, Journal of Combinatorial Theory, Series B \textbf{29}
  (1980), no.~2, 195--230.

\bibitem{duke2005continued}
William Duke, \emph{Continued fractions and modular functions}, Bulletin of the
  American Mathematical Society \textbf{42} (2005), no.~2, 137--162.

\bibitem{dyer1978automorphism}
Joan~L Dyer, \emph{Automorphism sequences of integer unimodular groups},
  Illinois Journal of Mathematics \textbf{22} (1978), no.~1, 1--30.

\bibitem{ward}
Manfred Einsiedler and Thomas Ward, \emph{Ergodic {Theory}: with a view towards
  {Number} {Theory}}, vol. 259, Springer Verlag, 2010.

\bibitem{el2012rational}
Abdelkrim El~Basraoui and Abdellah Sebbar, \emph{Rational equivariant forms},
  International Journal of Number Theory \textbf{8} (2012), no.~04, 963--981.

\bibitem{fishman2013closed}
Daniel Fishman and Steven~J Miller, \emph{Closed form continued fraction
  expansions of special quadratic irrationals}, ISRN Combinatorics
  \textbf{2013} (2013).

\bibitem{fort}
MK~Fort, \emph{A theorem concerning functions discontinuous on a dense set},
  American Mathematical Monthly (1951), 408--410.

\bibitem{hua1952automorphisms}
LK~Hua and I~Reiner, \emph{Automorphisms of the unimodular group}, Transactions
  of the American Mathematical Society \textbf{71} (1951), no.~3, 331--348.

\bibitem{imbert1997isomorphisme}
Michel Imbert, \emph{Sur l'isomorphisme du groupe de richard thompson avec le
  groupe de ptol{\'e}m{\'e}e}, LONDON MATHEMATICAL SOCIETY LECTURE NOTE SERIES
  (1997), 313--324.

\bibitem{iosifescu2002metrical}
Marius Iosifescu and Cor Kraaikamp, \emph{Metrical theory of continued
  fractions}, vol. 547, Springer Science \& Business Media, 2002.

\bibitem{isola2014continued}
Stefano Isola, \emph{Continued fractions and dynamics}, Applied Mathematics
  \textbf{5} (2014), no.~07, 1067.

\bibitem{jones1983operations}
Gareth~A Jones and John~S Thornton, \emph{Operations on maps, and outer
  automorphisms}, Journal of Combinatorial Theory, Series B \textbf{35} (1983),
  no.~2, 93--103.

\bibitem{jones1986automorphisms}
\bysame, \emph{Automorphisms and congruence subgroups of the extended modular
  group}, Journal of the London Mathematical Society \textbf{2} (1986), no.~1,
  26--40.

\bibitem{kamano2013analytic}
Ken Kamano, \emph{Analytic continuation of the lucas zeta and l-functions},
  Indagationes Mathematicae \textbf{24} (2013), no.~3, 637--646.

\bibitem{kaneko2003kappa}
Masanobu Kaneko and Masaaki Yoshida, \emph{The kappa function}, International
  Journal of Mathematics \textbf{14} (2003), no.~09, 1003--1013.

\bibitem{kessebohmer2008fractal}
Marc Kesseb{\"o}hmer and Bernd~O Stratmann, \emph{Fractal analysis for sets of
  non-differentiability of minkowski's question mark function}, Journal of
  Number Theory \textbf{128} (2008), no.~9, 2663--2686.

\bibitem{koruoglu2010}
Ozden Koruoglu and Recep Sahin, \emph{Generalized fibonacci sequences related
  to the extended hecke groups and an application to the extended modular
  group}, Turk J Math \textbf{34} (2010), 325--332.

\bibitem{koruoglu}
Özden Koruoğlu, Recep Şahin, and Sebahattin İkikardeş, \emph{Trace
  {Classes} and {Fixed} {Points} for the {Extended} {Modular} {Group}}, Turkish
  Journal of Mathematics \textbf{32} (2008), no.~1, 11--19.

\bibitem{lagarias2011takagi}
Jeffrey~C Lagarias, \emph{The takagi function and its properties}, arXiv
  preprint arXiv:1112.4205 (2011).

\bibitem{lewis2001period}
J~Lewis and Don Zagier, \emph{Period functions for maass wave forms. i}, Annals
  of Mathematics \textbf{153} (2001), no.~1, 191--258.

\bibitem{marcolli}
Yuri Manin and Matilde Marcolli, \emph{Modular shadows and the {Levy}-{Mellin}
  infinity-adic transform}, arXiv preprint math/0703718 (2007).

\bibitem{mayer1990thermodynamic}
Dieter~H Mayer, \emph{On the thermodynamic formalism for the gauss map},
  Communications in mathematical physics \textbf{130} (1990), no.~2, 311--333.

\bibitem{mcmullen}
Curtis~T McMullen, \emph{Uniformly {Diophantine} numbers in a fixed real
  quadratic field}, Compositio Mathematica \textbf{145} (2009), no.~04,
  827--844.

\bibitem{murty2013fibonacci}
M~Ram Murty, \emph{Fibonacci zeta function}, Automorphic Representations and
  L-Functions, TIFR Conference Proceedings, edited by D. Prasad, CS Rajan, A.
  Sankaranarayanan, J. Sengupta, Hindustan Book Agency, New Delhi, India, 2013.

\bibitem{mushtaq1}
Q~Mushtaq, \emph{Modular group acting on real quadratic fields}, Bulletin of
  the Australian Mathematical Society \textbf{37} (1988), no.~02, 303--309.

\bibitem{mushtaq2}
Qaiser Mushtaq, \emph{Coset diagrams for an action of the extended modular
  group on the projective line over a finite field}, Indian J. Pure Appl. Math
  \textbf{20} (1989), no.~8, 747--754.

\bibitem{similar}
Volodymyr Nekrashevych, \emph{Self-similar groups}, no. 117, American
  Mathematical Soc., 2005.

\bibitem{northshield}
S~Northshield, \emph{Circle boundaries of planar graphs}, Potential Analysis
  \textbf{2} (1993), no.~4, 299--314.

\bibitem{tugce}
Tuğçe Çolak, A~Muhammed Uludağ, and Ayberk Zeytin, \emph{Boundaries of
  {Planar} {Trees}}, (to appear).

\bibitem{ortner2008online}
Ronald Ortner, \emph{Online regret bounds for markov decision processes with
  deterministic transitions}, Algorithmic Learning Theory, Springer, 2008,
  pp.~123--137.

\bibitem{penner2012decorated}
Robert~C Penner, \emph{Decorated teichm{\"u}ller theory}, European Mathematical
  Society, 2012.

\bibitem{saber2014equivariant}
Hicham Saber, \emph{Equivariant functions for the m$\backslash$"$\{$o$\}$ bius
  subgroups and applications}, arXiv preprint arXiv:1412.8100 (2014).

\bibitem{manfred}
Manfred~Robert Schroeder, \emph{Fractals, {Chaos} and {Power} {Laws}}.

\bibitem{sebbar2012equivariant}
Abdellah Sebbar and Ahmed Sebbar, \emph{Equivariant functions and integrals of
  elliptic functions}, Geometriae Dedicata \textbf{160} (2012), no.~1,
  373--414.

\bibitem{smart1972meromorphic}
John~Roderick Smart, \emph{On meromorphic functions commuting with elements of
  a function group}, Proceedings of the American Mathematical Society (1972),
  343--348.

\bibitem{kunming}
A.~M. Uluda{\u{g}} and H.~Ayral, \emph{Modular group and its actions}, Handbook
  of group actions (S-T.~Yau L.~Ji, A.~Papadopoulos, ed.), International Press,
  2014.

\bibitem{dynamicaljimm}
A~Muhammed Uludag and Hakan Ayral, \emph{Dynamical aspects of the involution
  jimm}, (to appear).

\bibitem{uludag2012binary}
A~Muhammed Uludag, A.~Zeytin, and M.~Durmus, \emph{Binary quadratic forms as
  dessins}, preprint (2012).

\bibitem{vallee2006euclideans}
Brigitte Vall{\'e}e, \emph{Euclideans dynamics}, Discrete and Continuous
  Dynamical Systems series S (2006), 281--352.

\bibitem{vepstas2008minkowski}
Linas Vepstas, \emph{On the minkowski measure. arxiv}, arXiv preprint
  arXiv:0810.1265 (2008).

\bibitem{Viader1998}
Pelegrí Viader, Jaume Paradís, and Lluís Bibiloni, \emph{A new light on
  minkowski's ?(x) function}, Journal of Number Theory \textbf{73} (1998),
  no.~2, 212 -- 227.

\bibitem{zagier2000quelques}
Don Zagier, \emph{Quelques cons{\'e}quences surprenantes de la cohomologie de
  sl (2, z)}, Le{\c{c}}ons de Math{\'e}matiques d’aujourd’hui, Cassini,
  Paris (2000), 99--123.

\bibitem{zagier2001new}
\bysame, \emph{New points of view on the selberg zeta function}, Proceedings of
  the Japanese-German Seminar “Explicit structures of modular forms and zeta
  functions” Hakuba, 2001.

\end{thebibliography}

\noindent
{\bf Appendix.} Latex preamble commands for a latin-compatible typography of $\Jimm$
\begin{verbatim}
\usepackage{arabtex}
\usepackage{adjustbox}
%Upper-case size jimm
\newcommand{\Jimm}{{\adjustbox{scale=.8, raise=0.7ex,%
trim=0px 0px 0px 7px, padding=0ex 0ex 0ex 0ex}{\RL{j}}}}
 %lower-case size jimm, to use as a substcript
\newcommand{\jimm}{{\adjustbox{scale=.5, raise=0.6ex,%
trim=0px 0px 0px 7px, padding=0ex 0ex 0ex 0ex}{\RL{j}}}}
\end{verbatim}

 \end{document}